%% file: main.tex
\documentclass[]{article}

\usepackage[utf8]{inputenc}
\usepackage[T1]{fontenc}
\usepackage[hmarginratio=1:1,top=32mm,columnsep=20pt]{geometry} 
\usepackage[hang, small,labelfont=bf,up,textfont=it,up]{caption} 
\usepackage{booktabs} 
\usepackage{abstract} 

\usepackage{titlesec}
\titleformat{\section}[block]{\large\scshape\centering}{\thesection.}{1em}{} 
\titleformat{\subsection}[block]{\large}{\thesubsection.}{1em}{} 

\usepackage{amssymb, epsfig,amssymb, latexsym}
\usepackage{bm}
\usepackage{amsfonts,psfrag,amsmath,bbm,color,url}
\usepackage{graphicx}
\graphicspath{ {./img/} }
\usepackage{subfig}

\usepackage{amsfonts}
 \usepackage{amsmath}
 \usepackage{xcolor}
\usepackage{tikz}
\usepackage{pgfplots}
\usepackage{comment}
\usepackage{bm}
\usepackage{mathabx}
\DeclareMathAlphabet{\mathbbx}{U}{bbold}{m}{n}

\usepackage{cases}
\usepackage[output-decimal-marker={,}]{siunitx}
\usepackage{lineno,hyperref}
\definecolor{vargreen}{rgb}{0.0, 0.5, 0.0}
\definecolor{navyblue}{rgb}{0.0, 0.0, 0.5}
\definecolor{mediumorchid}{rgb}{0.73, 0.33, 0.83}
\definecolor{crimson}{rgb}{0, 0, 0}
\definecolor{lightseagreen}{rgb}{0.13, 0.7, 0.67}
\definecolor{royalblue}{rgb}{0.25, 0.41, 0.88}
\newcommand{\RA}[1]{{\color{crimson}#1}}
\definecolor{hotpink}{rgb}{1.0, 0.41, 0.71}
\definecolor{magenta}{rgb}{1.0, 0.0, 1.0}
\definecolor{goldenrod}{rgb}{0.85, 0.65, 0.13}

\newcommand{\TI}[1]{{\color{blue}TI: #1}}
\definecolor{plum(traditional)}{rgb}{0.56, 0.27, 0.52}
\usepackage{enumitem}

\newcommand{\mylabel}[2]{#2\def\@currentlabel{#2}\label{#1}}
\usetikzlibrary{shapes}
\usetikzlibrary{plotmarks}
\usetikzlibrary{shapes.geometric}
\usepackage[export]{adjustbox}
\usepackage{mdframed}
\usepackage{array}

\newtheorem{remark}{Remark}[section]

\def\PP{{{\rm l}\kern - .15em {\rm P} }}
\def\PN2{{\PP_{N}-\PP_{N-2}}}

\newcommand{\btau}{\boldsymbol{\tau}}

\newcommand{\ba}{\boldsymbol{a}}

\newcommand{\bb}{\boldsymbol{b}}

\newcommand{\bff}{\boldsymbol{f}}

\newcommand{\bu}{\boldsymbol{u}}

\newcommand{\deleted}[1]{{}}

\pgfplotsset{compat=1.17}
\usepackage[mode=buildnew]{standalone}
\title{\textbf{Pressure Data-Driven Variational Multiscale Reduced Order Models}}

\date{ }
\pgfplotsset{compat=1.17}
\author{Anna Ivagnes  \\ \small SISSA, International School for Advanced Studies, \\ \small Mathematics Area, mathLab, Trieste, Italy. \\ \small  \href{mailto:aivagnes@sissa.it}{aivagnes@sissa.it} \normalsize \and Giovanni Stabile \\ \small SISSA, International School for Advanced Studies, \\ \small Mathematics Area, mathLab, Trieste, Italy. \\ \small  \href{mailto:gstabile@sissa.it}{gstabile@sissa.it} \normalsize  \and Andrea Mola  \\\small  Multi-scale Analysis of Materials Unit, \\ \small Scuola IMT Alti Studi, Lucca, Italy. \\ \small \href{mailto:andrea.mola@imtlucca.it}{andrea.mola@imtlucca.it} \normalsize \and Traian Iliescu  \\\small 
Department of Mathematics, Virginia Tech, \\ \small Blacksburg, VA, USA. \\ \small \href{mailto:iliescu@vt.edu}{iliescu@vt.edu} \normalsize \and Gianluigi Rozza  \\ \small  SISSA, International School for Advanced Studies, \\\small  Mathematics Area, mathLab, Trieste, Italy. \\ \small \href{mailto:grozza@sissa.it}{grozza@sissa.it}}

\begin{document}
\maketitle
\begin{abstract}
In this paper, we develop data-driven closure/correction terms to increase the pressure and velocity accuracy of reduced order models (ROMs) for fluid flows.
Specifically, we propose the first pressure-based data-driven variational multiscale ROM, in which we use the available data to construct closure/correction terms for both the momentum equation and the continuity equation. 
Our numerical investigation of the two-dimensional flow past a circular cylinder at 
$Re = \num[group-separator={,}]{50000}$ in the marginally-resolved regime shows that the novel pressure data-driven variational multiscale ROM yields significantly more accurate velocity and pressure approximations than the standard ROM and, more importantly, than the original data-driven variational multiscale ROM (i.e., without pressure components).
In particular, our numerical results show that adding the closure/correction term in the momentum equation significantly improves both the velocity and the pressure approximations, whereas adding the closure/correction term in the continuity equation improves only the pressure approximation.
\end{abstract}



\input{sections_paper1/section_introduction}
\input{sections_paper1/section_fom}
\input{sections_paper1/section_grom}
\input{sections_paper1/section_datasup}
\input{sections_paper1/section_datappe}

\input{sections_paper1/section_numresults}
\input{sections_paper1/section_conclusions}
\input{sections_paper1/section_acknowledgements}
\bibliographystyle{plain}
\bibliography{main}
\end{document}

%% file: sections_paper1/section_introduction.tex
\section{Introduction}


Reduced order models (ROMs) have demonstrated to be a key methodology to accelerate the resolution of numerical problems governed by partial differential equations \cite{handbook,hesthaven2015certified,quarteroni2015reduced,benner2015survey} 
for both linear \cite{rozza2008reduced} and nonlinear 
\cite{carlberg2013gnat,wilcoxNLROM2019,parish2020adjoint} problems. In this article,  we 
focus on projection-based ROMs for fluids \cite{stabile2019reduced,grimberg2020stability,noack2005need,iollo2000stability,reyes2020projection,rebollo2017certified,iliescu2014proper, rozza2013reduced}, 
which have been used in a large variety of applications 
in industry, geosciences, and biomedicine.
Moreover, we exclusively focus on reduced basis Galerkin ROMs 
in which the reduced basis space 
is constructed by 
using the proper orthogonal decomposition 
(POD) \cite{Sir87abc}.
These ROMs are built starting from the general Galerkin framework:
start with a set of basis functions (modes), $\{ \varphi_{1}, \ldots, \varphi_{r} \}$, express the unknown solution as a linear combination of these modes, $u(x,t) = \sum_{i=1}^{r} a_i(t) \, \varphi_i(x)$, and project the equations onto the space spanned by these modes.
The resulting Galerkin ROM (G-ROM) is a system of equations in which the unknowns are the coefficients in the linear combination ansatz used above:
\begin{eqnarray}
    \ba_t = \bff(\ba),
    \label{eqn:g-rom}  
\end{eqnarray}
where $\ba(t) := (a_1(t), \ldots, a_r(t))$.
The main difference between ROMs and the classical Galerkin methods (of which the quintessential example is the finite element method (FEM)), is that the ROM basis is a data-driven basis, i.e., a basis constructed from the available numerical or experimental data, whereas the classical Galerkin methods do not generally use data to build the basis. 

There are 
significant challenges in the development of ROMs for realistic fluid flow applications.
For example, 
finding a practical compromise between 
computational effort and 
accuracy in 
the reduced order modeling of realistic flows still remains an elusive goal.  Furthermore, 
ROMs for computational fluid dynamics (CFD) suffer from different types of stability issues, 
which are mainly of two different types: (i) instabilities associated with pressure recovery and non-compliant inf-sup spaces \cite{stabile2018finite,ballarin2015supremizer,caiazzo2014numerical,gerner2012, ali2020stabilized}, and (ii) instabilities due to 
the convection-dominated regime and turbulence 
\cite{farhat2015structure}.
Finally, one of the most important challenges in reduced order modeling of realistic flows is the {\it under-resolved regime}, i.e., using fewer ROM modes than the number required to accurately approximate the dynamics of the given 
flow.
In classical CFD, the under-resolved regime is generally associated with spatial and temporal meshes that are too coarse to represent the underlying dynamics.
For example, in the numerical simulation of turbulent flows (which are multiscale, chaotic phenomena), a standard ROM approximation could require hundreds or even thousands of ROM modes~\cite{ahmed2021closures}, which would significantly increase the ROM computational cost.
To ensure a low ROM computational cost, ROMs are generally built with relatively few basis functions.
In that case, ROMs are used in the under-resolved regime.

We note that there is also another regime, called the {\it marginally-resolved regime}, which is an intermediate regime, between the under-resolved regime and the fully-resolved regime (i.e., when the number of modes is enough to represent the underlying dynamics).
In the marginally-resolved regime, the number of ROM basis functions is enough to represent the main features of the underlying dynamics, but the standard G-ROM yields inaccurate approximations.
One example of ROM in a marginally-resolved regime is a low-dimensional ROM in the numerical simulation of a laminar flow (i.e., at a low Reynolds number).

In the marginally-resolved regime and, especially, in the under-resolved regime, one popular approach to increase the accuracy of the standard G-ROM~\eqref{eqn:g-rom} is to add an extra term, $\btau$, to the right-hand side of~\eqref{eqn:g-rom}:
\begin{eqnarray}
    \ba_t 
    = \bff(\ba)
    + \btau(\ba).
    \label{eqn:g-rom+closure}  
\end{eqnarray}
In the under-resolved regime, the new term in~\eqref{eqn:g-rom+closure}, $\btau$, is called the {\it ROM closure term}.
The same terminology is used in classical large eddy simulation (LES)~\cite{BIL05,sagaut2006large}.
The role of the ROM closure term is to model the effect of the discarded ROM modes, $\{ \varphi_{r+1}, \ldots \}$, on the ROM dynamics.
In the marginally-resolved regime, the new term in~\eqref{eqn:g-rom+closure}, $\btau$, is called the {\it ROM correction term}~\cite{mou2020data} and its role is to increase the ROM accuracy.

There are three main approaches to model the closure/correction term $\btau$ in~\eqref{eqn:g-rom+closure}:
(i) Functional modeling, in which physical insight is used to build a model for $\btau$.
(ii) Structural modeling, in which mathematical methods (e.g., expansions and asymptotics) are used to build a model for $\btau$.
(iii) Data-driven modeling, in which available numerical or experimental data is used to build a model for $\btau$.

In this paper, we exclusively focus on data-driven modeling of the closure/correction term $\btau$ in~\eqref{eqn:g-rom+closure}.
Specifically, we investigate the data-driven variational multiscale ROM (DD-VMS-ROM) framework put forth in~\cite{mou2021data,xie2018data}.
The classical VMS framework was pioneered by Hughes and his collaborators more than two decades ago as a means to model multiscale phenomena in science and engineering. 
The VMS methodology has become popular in classical Galerkin methods, such as the FEM.
The VMS framework has also been applied in a ROM setting.
Functional (physical) VMS-ROMs were proposed in~\cite{bergmann2009enablers,wang2012proper,iliescu2013variational,iliescu2014variational}.
The DD-VMS-ROMs, which are the focus of this paper, are built around a completely different principle:
instead of physical insight, they use available data to construct the closure/correction term $\btau$ in~\eqref{eqn:g-rom+closure}.
(For related, yet different approaches, see, e.g.,~\cite{baiges2015reduced,selten1997statistical} in a ROM setting and~\cite{pradhan2020variational,pradhan2021variational} in a 
FEM setting.)
Specifically, in the offline stage, we first postulate an ansatz (i.e., a model form) for the closure/correction term $\btau$, and then solve a least squares problem to find the ansatz parameters that yield the best fit between the exact closure/correction term $\btau$ computed with the full order model (FOM) data and the ansatz.
The DD-VMS-ROM has been successfully used in the numerical simulation of a 2D flow past a cylinder, the quasi-geostrophic equations, and the 3D turbulent channel flow at $Re_{\tau} = 395$.


We emphasize that all DD-VMS-ROMs have been developed exclusively for the velocity correction, i.e., introducing a correction/closure term for term $\mathbf{f}(\mathbf{a})$ in \eqref{eqn:g-rom+closure}. In 
\cite{xie2018data, mohebujjaman2019physically, mou2020data}, it was shown that this approach leads to improvements in the velocity approximation accuracy.
In this paper, we propose the first DD-VMS-ROM for the improvement of {\it pressure} approximation accuracy.
We note that the importance of pressure approximation in reduced order modeling has been recognized early on (see, e.g., the pioneering work of Veroy-Grepl and Rozza~\cite{rozza2007stability} in the reduced basis method (RBM) community, and more recent work in the POD community~\cite{caiazzo2014numerical,decaria2020artificial,kean2020error,noack2005need,rubino2020numerical}). 
For example, a ROM pressure approximation is needed to compute important engineering quantities (e.g., lift and drag) or when the snapshots used to construct the ROM basis are not weakly divergence-free (e.g., when the mass conservation equation is only approximately enforced in the FOM~\cite{caiazzo2014numerical,decaria2020artificial}).
Thus, it is important to investigate whether the accuracy of the standard ROM pressure approximation can be increased.
This is precisely the main aim of this paper.

To present the novel pressure DD-VMS-ROM model, we first note that the standard G-ROM~\eqref{eqn:g-rom} is modified to include a pressure approximation.
\begin{equation}
    \begin{cases}
        \ba_t 
        = \bff(\ba,\bb), & \\
        {\bf g}(\ba)
        = {\bf 0} , & 
    \end{cases}
    \label{eqn:p-g-rom}  
\end{equation}
where $\bb(t) := (b_1(t), \ldots, b_q(t))$ is the vector of unknown coefficients in the ROM pressure approximation.
Furthermore, there are now two types of ROM closure/correction terms instead of one (as in~\eqref{eqn:g-rom+closure}):
\begin{equation}
    \begin{cases}
        \ba_t 
        = \bff(\ba,\bb)
        + \btau_{\bu}(\ba, \bb), & \\
        {\bf g}(\ba)
        + \btau_{p}(\ba, \bb)
        = {\bf 0} . & 
    \end{cases}
    \label{eqn:p-g-rom+closure}  
\end{equation}
To construct models for the closure/correction terms $\btau_{\bu}$ and $\btau_{p}$ in~\eqref{eqn:p-g-rom+closure}, we use again the available numerical or experimental data.
Specifically, in the offline stage, we use the FOM data to compute the true (accurate) closure/correction terms $\btau_{\bu}$ and $\btau_{p}$, we postulate model forms (i.e., ansatzes) for these two terms, and then we solve one or two least square problems to find the model form parameters that ensure the best fit between the true models and the ansatzes.
Our numerical investigation shows that the resulting model, which we call the pressure DD-VMS-ROM, yields significantly more accurate velocity and pressure approximations than the standard G-ROM~\eqref{eqn:p-g-rom} and, more importantly, than the original DD-VMS-ROM~\cite{mou2021data,xie2018data}.
In particular, our numerical results show that 
adding the term $\btau_{\bu}$ in the momentum equation significantly improves both the velocity and the pressure approximations, whereas the $\btau_{p}$ contribution improves only the pressure approximation.

The rest of the paper is organized as follows: In section~\ref{sec:g-rom}, we outline the construction of the standard Galerkin ROM, presenting two different approaches to construct the ROM basis, both based on the 
POD: 
the first approach (section~\ref{sup}) is based on supremizers (i.e., ROM modes that are added to the standard ROM velocity modes in order to fulfill the inf-sup condition), whereas the second approach (section~\ref{ppe_sec}) is based on solving a pressure Poisson equation (PPE). The PPE approach is 
used to improve the accuracy of the pressure approximation. In sections \ref{datasup} and \ref{datappe}, 
the novel pressure DD-VMS-ROM framework is constructed for the supremizer and PPE approaches, respectively.
Section \ref{results} 
presents a numerical investigation of the new pressure DD-VMS-ROM framework in the simulation of a two-dimensional flow past a circular cylinder at a Reynolds number $Re=\num[group-separator={,}]{50000}$.

\bigskip




%% file: sections_paper1/section_fom.tex
\section{Full Order Problem}
As a mathematical model, we use 
the Navier-Stokes equations (NSE) for incompressible flows. 
We use the following notation: 
$\Omega \in \mathbb{R}^d$ with 
$d=2$ or $3$ is the fluid domain, $\Gamma$ its boundary, $t \in [0,T]$ 
the time, $\mathbf{u}=\mathbf{u}(\mathbf{x},t)$ 
the flow velocity vector field, 
$p=p(\mathbf{x},t)$ 
the normalized pressure scalar field divided by the fluid density, and $\nu$ 
the fluid kinematic viscosity. The strong form of the NSE is the following:
\begin{subnumcases}{\label{NSE}}
\frac{\partial \mathbf{u}}{\partial t}=-\nabla \cdot (\mathbf{u} \otimes \mathbf{u})+\nabla \cdot \RA{ (}\nu (\nabla \mathbf{u} + (\nabla \mathbf{u})^T ) \RA{)}-\nabla p & in $\Omega \times [0,T]\, , $ \label{mom_NSE}\\
\nabla \cdot \mathbf{u} = \mathbf{0} & in $\Omega \times [0,T]\, ,$ \label{cont_NSE}\\
+ \text{ boundary conditions }& on  $\Gamma \times [0,T]\, ,$ \label{bound_NSE}\\
+ \text{ initial conditions }& in  $(\Omega,0)\, .$ \label{init_NSE}
\end{subnumcases}

In the computation of the full order solutions of 
the NSE in \eqref{NSE}, we 
employ a finite volume discretization 
and \emph{OpenFOAM}. The finite volume method 
\cite{moukalled2016finite} 
is a mathematical technique that converts the partial differential equations (the NSE in our case) defined on differential volumes 
to algebraic equations defined on finite volumes.
The first preliminary step of this method is a polyhedral discretization of the domain, in order to define the finite control volumes. The second step is to integrate the NSE over each control volume of the domain of interest. The divergence theorem is used to convert the volume integrals 
to surface integrals. Those integrals are finally discretized as sums of the fluxes at the boundary faces of each control volume.  

At the full order level, 
the Reynolds--averaged Navier--Stokes (RANS) approach is used in the numerical simulation of the turbulent flow.
The RANS approach is based on the 
Reynolds decomposition 
\cite{reynolds1895iv}, in which each flow field is expressed as the sum of its mean and its fluctuating part. The RANS equations are obtained by taking the time-average of the 
NSE in \eqref{NSE} and adding a closure 
term to model the 
Reynolds stress tensor. 
The closure problem is solved by using an eddy viscosity 
model, which 
is based on the Boussinesq hypothesis. In this paper, the SST $\kappa-\omega$ model is used. This model is based on the resolution of two additional transport equations to describe the kinetic energy $\kappa$ and the specific turbulent dissipation rate $\omega$. The 
$\kappa-\omega$ model is presented in its standard form in \cite{kolmogorov1941equations}, whereas the SST formulation is developed in \cite{menter1994two}. For more details on the SST $\kappa-\omega$ model equations, we refer to \cite{hijazi2020data}.

%% file: sections_paper1/section_grom.tex
\section{The POD-Galerkin ROM} 
    \label{sec:g-rom}
The online phase consists 
of the time integration of the G-ROM that is 
obtained by performing a Galerkin projection of the 
NSE onto the POD space. In order to construct a ROM 
that can approximate both the velocity and pressure field, it is required to 
use stabilization approaches \cite{hijazi2020data}. In this work, we 
explore two different stabilization strategies: 
\begin{itemize}
    \item the SUP-ROM: a reduced order method based on a reduced version of system \eqref{NSE}, in which additional modes for the velocity space, named \emph{supremizer} modes, are introduced in order to fulfill the \emph{inf-sup} condition \cite{rozza2007stability, ballarin2015supremizer, stabile2018finite};
    \item the PPE-ROM: a reduced order method in which the continuity equation \eqref{cont_NSE} is replaced by the pressure Poisson equation \cite{akhtar2009stability,Stabile2017CAIM,stabile2018finite}.
\end{itemize}
In what follows, for the sake of completeness, we 
briefly recall the reduced order systems obtained by the two different formulations. 
\subsection{The SUP-ROM approach}
\label{sup}
The first approach we followed is the SUP-ROM approach, first introduced in \cite{rozza2007stability} and explored in \cite{ballarin2015supremizer, stabile2018finite}, in which additional velocity modes are introduced.
The method consists in computing the FOM solutions for different time instants $\{t_j\}_{j=1}^{N_T}$; each of the 
full order solutions is called FOM snapshot. 

The 
POD is applied 
to the FOM snapshot matrices 
\begin{equation*}
\mathcal{S}_u=\{\mathbf{u}(\mathbf{x},t_1),...,\mathbf{u}(\mathbf{x},t_{N_T})\}  \in \mathbb{R}^{N_u^h \times N_T}, \quad \mathcal{S}_p=\{p(\mathbf{x},t_1),...,p(\mathbf{x},t_{N_T})\}  \in \mathbb{R}^{N_p^h \times N_T},
\end{equation*}
where $N_u^h$ and $N_p^h$ are the numbers of degrees of freedom for the velocity and pressure fields. \RA{We specify that, since the velocity is a vector field and the pressure is a scalar field, in our case $N_u^h=3 \times N_{cells}$ 
and $N_p^h=N_{cells}$, where $N_{cells}$ is the total number of cells of the computational grid.}
After a POD modal decomposition is performed, the velocity and pressure POD spaces are assembled as follows:
\begin{equation}
\mathbb{V}^u_{\text{POD}}=\mbox{span}\{[\boldsymbol{\phi}_i]_{i=1}^{N_u}\},\quad 
\mathbb{V}^p_{\text{POD}}=\mbox{span}\{[\chi_i]_{i=1}^{N_p}\},
\end{equation}
where $N_u\ll N_u^h$ and $N_p\ll N_p^h$, and  $[\boldsymbol{\phi}_i]_{i=1}^{N_u}$ and $[\chi_i]_{i=1}^{N_p}$ are the velocity and pressure POD modes, respectively. $[\mathbf{s}(\chi_i)]_{i=1}^{N_{sup}}$ are the velocity supremizer modes, which are additional modes introduced in order to fulfill  the \emph{inf-sup} condition \cite{stabile2018finite}.

For each pressure basis function, the corresponding supremizer element can be found by solving the following problem:
\begin{equation}
\begin{cases}
 \Delta \mathbf{s}_i = - \nabla p_i \text{ in } \Omega , \\
 \mathbf{s}_i = 0 \text{ on } \partial \Omega .
\end{cases}
\label{sup_condition}
\end{equation}

\RA{
We 
note that each supremizer mode $\mathbf{s}_i$ 
corresponds to the Riesz representation in $H_0^1(\Omega)^{d}$ of the linear, continuous functional associated with $\nabla p$, 
where $d$ is the spatial dimension.}

The enrichment of the velocity POD space can be addressed either with an \emph{exact} or with an \emph{approximated} approach \cite{ballarin2015supremizer}. In the exact approach, the problem \eqref{sup_condition} is solved for each pressure basis function $\chi_i$ and each solution is added to the velocity space. In the approximated approach, the problem \eqref{sup_condition} is solved for each pressure
snapshot $p(\mathbf{x}, t_i)_{i=1}^{N_T}$, 
which yields the following supremizer snapshot matrix: 
\begin{equation*}
\mathcal{S}_{sup}=\{\mathbf{s}(\mathbf{x},t_1),...,\mathbf{s}(\mathbf{x},t_{N_T})\}  \in \mathbb{R}^{N_u^h \times N_T}.
\end{equation*}
A POD modal decomposition is then applied to the snapshot matrix in order to obtain the supremizer POD modes $(\boldsymbol{\eta}_i)_{i=1}^{N_{sup}}$ \cite{stabile2018finite}.
In this paper, we adopt the approximated procedure, since it 
significantly reduces the computational cost of the offline phase. 
We note, however, that when an approximated approach is used, it is not possible to rigorously demonstrate the inf-sup condition, which 
can lead to stability issues, as we will see in section \ref{sup_enrich_nocorr}.

For the sake of simplicity, we 
denote $(\boldsymbol{\eta}_i)_{i=1}^{N_{sup}}=(\boldsymbol{\phi}_i)_{i=N_u+1}^{N_u+N_{sup}}$ and\RA{, considering a generic time instant $t_j$, we 
}
express the approximated velocity and pressure fields 
as follows:

\begin{equation}
\mathbf{u}(\mathbf{x},t\RA{_j}) \approx \mathbf{u}_r(\mathbf{x},t\RA{_j})=\sum_{i=1}^{N_u+N_{sup}} a_i(t\RA{_j})\boldsymbol{\phi}_i(\mathbf{x}), \quad
p(\mathbf{x},t\RA{_j})\approx p_r(\mathbf{x},t\RA{_j})=\sum_{i=1}^{N_p} b_i(t\RA{_j})\chi_i(\mathbf{x}). 
\label{appfields}
\end{equation}
Performing a Galerkin projection of the momentum equation \eqref{mom_NSE} onto the velocity modes and of the continuity equation \eqref{cont_NSE} onto the pressure modes, the following reduced system is obtained: 
\begin{equation}
    \begin{cases}
    \mathbf{M} \dot{\mathbf{a}}=\nu(\mathbf{B}+\mathbf{B_T})\mathbf{a}-\mathbf{a}^T \mathbf{C} \mathbf{a}-\mathbf{H}\mathbf{b}+\tau \left( \sum_{k=1}^{N_{\text{BC}}}(U_{\text{BC},k}\mathbf{D}^k-\mathbf{E}^k \mathbf{a})\right)\, ,\\
    \mathbf{P}\mathbf{a}=\mathbf{0}\, ,
    \end{cases}
    \label{reduced_system}
\end{equation}
where $\mathbf{a}$ and $\mathbf{b}$ are the vectors of the coefficients associated with the velocity and pressure modes, respectively.

The matrices appearing in the system are defined 
as follows: 
\[
\begin{split}
&(\mathbf{M})_{ij}=(\boldsymbol{\phi}_i,\boldsymbol{\phi}_j)_{L^2(\Omega)}, \quad (\mathbf{P})_{ij}=(\chi_i,\nabla \cdot \boldsymbol{\phi}_j)_{L^2(\Omega)}\, ,\quad (\mathbf{B})_{ij}=(\boldsymbol{\phi}_i,\nabla \cdot \nabla \boldsymbol{\phi}_j)_{L^2(\Omega)}, \\
&(\mathbf{B_T})_{ij}=(\boldsymbol{\phi}_i,\nabla \cdot (\nabla \boldsymbol{\phi}_j)^T)_{L^2(\Omega)},\quad (\mathbf{C})_{ijk}=(\boldsymbol{\phi}_i,\nabla \cdot (\boldsymbol{\phi}_j \otimes \boldsymbol{\phi}_k))_{L^2(\Omega)}, \quad (\mathbf{H})_{ij}=(\boldsymbol{\phi}_i,\nabla \chi_j)_{L^2(\Omega)}\, .
\end{split}
\]
The term $\tau \left( \sum_{k=1}^{N_{\text{BC}}}(U_{\text{BC},k}\mathbf{D}^k-\mathbf{E}^k \mathbf{a})\right)$ in \eqref{reduced_system} is a penalization term used to enforce the Dirichlet boundary conditions at the reduced order level \cite{hijazi2020data, star2019extension}. 
In the penalization term formula,  $N_{\text{BC}}$ is the number of velocity boundary conditions on the $k $ different parts of the Dirichlet boundary,  $U_{\text{BC},k}$ is the velocity nonzero component at the $k$-th part of the Dirichlet boundary, $\tau$ is a penalization factor, and the matrices $\mathbf{E}^k$ and vectors $\mathbf{D}^k$ are defined as follows:
\[(\mathbf{E}^k)_{ij}=(\boldsymbol{\phi}_i, \boldsymbol{\phi}_j)_{L^2(\Gamma_{D_k})}, \quad (\mathbf{D}^k)_{i}=(\boldsymbol{\phi}_i)_{\Gamma_{D_k}}, \text{ for all }k=1,...,N_{\text{BC}}.\]

\subsection{The PPE-ROM approach} \label{ppe_sec}
The second stabilization technique for the reduced system is the PPE-ROM, which was 
proposed in \cite{noack2005need} and then used and extended 
in \cite{ akhtar2009stability, Stabile2017CAIM, stabile2018finite}.
In this formulation, the continuity equation of the reduced system \eqref{reduced_system} is replaced by the pressure Poisson equation, obtained by taking the divergence of the momentum equation and 
leveraging the fact that the velocity field is divergence-free:  
\begin{equation}
    \begin{cases}
    \mathbf{M} \dot{\mathbf{a}}=\nu(\mathbf{B}+\mathbf{B_T})\mathbf{a}-\mathbf{a}^T \mathbf{C} \mathbf{a}-\mathbf{H}\mathbf{b}+\tau \left( \sum_{k=1}^{N_{\text{BC}}}(U_{BC,k}\mathbf{D}^k-\mathbf{E}^k \mathbf{a})\right)\, ,\\
    \mathbf{D}\mathbf{b}+ \mathbf{a}^T \mathbf{G} \mathbf{a} - \nu \mathbf{N} \mathbf{a}- \mathbf{L}=\mathbf{0}\, .
    \end{cases}
    \label{reduced_systemPPE}
\end{equation}

The additional matrices appearing in system \eqref{reduced_systemPPE} are defined by Galerkin projection 
as follows:
\[
\begin{split}
&(\mathbf{D})_{ij}=(\nabla \chi_i,\nabla \chi_j)_{L^2(\Omega)}, \quad 
(\mathbf{G})_{ijk}=(\nabla \chi_i,\nabla \cdot (\boldsymbol{\phi}_j \otimes \boldsymbol{\phi}_k))_{L^2(\Omega)}, \\ &(\mathbf{N})_{ij}=(\mathbf{n} \times \nabla \chi_i,\nabla \boldsymbol{\phi}_j)_{\Gamma}, \quad (\mathbf{L})_{ij}=(\chi_i,\mathbf{n} \cdot \boldsymbol{R}_t)_{\Gamma}\, .
\end{split}
\]

In system \eqref{reduced_systemPPE}, the vector $\mathbf{R}$ is such that 
\begin{equation*}
\mathbf{u}(\mathbf{x},t)=\mathbf{R}(\mathbf{x}) \text{ on } \Gamma_{\text{Inlet}}\, ,
\end{equation*}
where $\Gamma_{\text{Inlet}}$ is the 
inlet boundary of the domain.
In the 
test case 
considered in the numerical investigation, the velocity conditions at the inlet do not change as time evolves and the term $(\mathbf{n} \cdot \mathbf{R}_t, \chi_i)_{\Gamma}$ is identically zero.

The penalty term $\tau \left( \sum_{k=1}^{N_{\text{BC}}}(U_{BC,k}\mathbf{D}^k-\mathbf{E}^k \mathbf{a})\right)$ in \eqref{reduced_systemPPE} is a constraint expressing the Dirichlet non-homogeneous boundary conditions, just as in the SUP-ROM described in section \ref{sup}.



%% file: sections_paper1/section_datasup.tex
\section{Data-driven corrections for the SUP-ROM approach}
\label{datasup}

To construct the new correction terms in this section (and in section~\ref{datappe}), we leverage the \emph{data-driven variational multiscale ROM} (\emph{DD-VMS-ROM}) framework~\cite{mou2021data,xie2018data}. This framework centers around the following principle, which is also used in LES~\cite{BIL05,sagaut2006large}:
\begin{itemize}
    \item[(i)] Filter the NSE with a spatial filter;
    \item[(ii)] Solve the filtered NSE obtained in step (i) to approximate the filtered flow variables.
\end{itemize}
The motivation for using this approach is that the filtered flow variables can be represented in the under-resolved (or marginally-resolved) regime, whereas the unfiltered flow variables cannot.
In the DD-VMS-ROM framework, the spatial filter used is the ROM projection~\cite{mou2021data,xie2018data}, i.e., the $L^2$ projection of the FOM variables onto the ROM space.

Specifically, for fixed $r$ and $q$, given the FOM velocity $\mathbf{u}$ 
and the FOM pressure $p$, the ROM filter \RA{
is defined as follows:} 
\begin{equation*}
(\bar{\mathbf{u}}^r,\boldsymbol{\phi}_i)=(\mathbf{u},\boldsymbol{\phi}_i)  \ \  \forall i=1,...,r, \quad
(\bar{p}^r,\chi_i)=(p,\chi_i)  \ \ \forall i=1,...,q,
\end{equation*}
where, in our setting, $r=N_u+N_{sup}$ and $q=N_p$.
We emphasize that, to our knowledge, using the ROM projection to define the spatially filtered pressure is novel.
Next, we use the ROM filters defined above and the DD-VMS-ROM framework to construct the new correction terms.

The reduced system obtained by projecting the NSE 
onto the POD modes can be reformulated as a system of spatially filtered NSE:
\begin{equation}
\begin{cases}
\left( -\frac{\partial \bar{\mathbf{u}}_d^r}{\partial t},\boldsymbol{\phi}_i \right) +\nu \left( \nabla \cdot \left( \nabla \bar{\mathbf{u}}_d^r - (\nabla \bar{\mathbf{u}}_d^r)^T \right), \boldsymbol{\phi}_i \right)-\left( (\bar{\mathbf{u}}_d^r \cdot \nabla) \bar{\mathbf{u}}_d^r,\boldsymbol{\phi}_i \right)-(\nabla\bar{p}_d^r,\boldsymbol{\phi}_i)+\\
\quad \quad+\left( \boldsymbol{\tau}_u^{SFS},\boldsymbol{\phi}_i\right)+\left( \boldsymbol{\tau}_{p(1)}^{SFS},\boldsymbol{\phi}_i\right) =0 \ \ \text{ for } i=1,...,r ,\\
\left(\nabla \cdot \bar{\mathbf{u}}_d^r,\chi_i \right)+\left( \tau_{p(2)}^{SFS},\chi_i\right)=0 \ \ \text{ for } i=1,...,q  ,
\end{cases}
\label{filtered_NSE}
\end{equation}
where the correction terms for the velocity and pressure 
are defined as follows:
\begin{equation}
\boldsymbol{\tau}_u^{SFS}=-\left(\overline{\left(\mathbf{u}_d \cdot \nabla \right)\mathbf{u}_d}^r-(\bar{\mathbf{u}}_d^r \cdot \nabla) \bar{\mathbf{u}}_d^r \right)  ,\quad
\boldsymbol{\tau}_{p(1)}^{SFS}=-\left(\overline{\nabla p_d}^r-\nabla \bar{p}_d^r\right),\quad
\tau_{p(2)}^{SFS}=\overline{\nabla \cdot \mathbf{u}_d}^r-\nabla \cdot \bar{\mathbf{u}}_d^r .
\label{closure_terms}
\end{equation}
\RA{We 
note that $\boldsymbol{\tau}_u^{SFS} \in \mathbb{R}^{N_{cells}}$, $\boldsymbol{\tau}_{p(1)}^{SFS} \in \mathbb{R}^{3\times N_{cells}}$, and $\boldsymbol{\tau}_{p(2)}^{SFS} \in \mathbb{R}^{N_{cells}}$.}
In order to reduce the computational effort, we assume that $\mathbf{u}_d \approx \mathbf{u}$ and $p_d \approx p$. 
Thus, the projection in the filtered NSE \eqref{filtered_NSE} is carried out not starting from the FOM fields, but starting from the fields reconstructed using a number $d$ and $d_p$ of modes, where $d$ and $d_p$ are smaller than the rank of the snapshot matrix. In \cite{xie2018data}, it was shown that this approach reduces the computational cost of the DD-VMS-ROM framework without significantly degrading its accuracy.

We note that, since the ROM projection is used as a spatial filter, the filtered variables satisfy the following formulas:
\[
\bar{\mathbf{u}}_d^r=\mathbf{u}_r \text{ and } \bar{p}_d^r=p_r.\]

Next, the following terms are introduced:
\begin{equation}
\boldsymbol{\tau}^u \text{ s.t. }\tau^u_i= \left( \boldsymbol{\tau}_u^{SFS},\boldsymbol{\phi}_i\right) , \quad
\boldsymbol{\tau}^{p(1)} \text{ s.t. }\tau^{p(1)}_i= \left( \boldsymbol{\tau}_{p(1)}^{SFS},\boldsymbol{\phi}_i\right) ,\quad
\boldsymbol{\tau}^{p(2)} \text{ s.t. }\tau^{p(2)}_i= \left( \tau_{p(2)}^{SFS},\chi_i\right).
\label{tau_all}
\end{equation}
\RA{We specify that $\boldsymbol{\tau}^u \in \mathbb{R}^{r}$, $\boldsymbol{\tau}^{p(1)} \in \mathbb{R}^{r}$, and $\boldsymbol{\tau}^{p(2)} \in \mathbb{R}^{q}$.}

Adding the 
correction terms, the dynamical system \eqref{reduced_system} becomes:
\begin{equation}
\begin{cases}
\mathbf{M} \dot{\mathbf{a}}=\nu(\mathbf{B}+\mathbf{B_T})\mathbf{a}-\mathbf{a}^T \mathbf{C} \mathbf{a}-\mathbf{H}\mathbf{b}+\tau \left( \sum_{k=1}^{N_{\text{BC}}}(U_{BC,k}\mathbf{D}^k-\mathbf{E}^k \mathbf{a})\right)+\boldsymbol{\tau}^u+\boldsymbol{\tau}^{p(1)} ,\\
    \mathbf{P}\mathbf{a}+\boldsymbol{\tau}^{p(2)}=\mathbf{0} .
\end{cases}
\label{new_reduced_system}
\end{equation}
The new system obtained in \eqref{new_reduced_system} is not a closed system because the correction terms depend on the fields $\mathbf{u}_d$ and $p_d$. In order to close the 
system, a data-driven modeling is adopted, as in \cite{xie2018data} and \cite{mohebujjaman2019physically}. The key problem is to find the approximated expressions $\boldsymbol{\tau}^u \approx \boldsymbol{\tau}^u (\mathbf{a}) $,  $\boldsymbol{\tau}^{p(1)}\approx \boldsymbol{\tau}^{p(1)} (\mathbf{b})$, and  $\boldsymbol{\tau}^{p(2)}\approx \boldsymbol{\tau}^{p(2)} (\mathbf{a})$. 

\subsection{Data-driven correction for velocity} \label{velocity}
In this section, only the correction term for velocity $\boldsymbol{\tau}^u$ is 
considered in system \eqref{new_reduced_system} and its effect on the dynamical system is 
evaluated. This technique, 
which was introduced in~\cite{xie2018data, mohebujjaman2019physically, mou2020data}, 
is applied to a 
new reduced formulation.

The correction term for velocity is modeled as in \cite{xie2018data} by using the ansatz 
\begin{equation}
\boldsymbol{\tau}^u (\mathbf{a})=\tilde{A} \mathbf{a}+\mathbf{a}^T \tilde{B} \mathbf{a} ,
\label{ansatz1}
\end{equation}
where $\tilde{B}$ is a three-dimensional tensor.
To find $\tilde{A}$ and $\tilde{B}$, the following optimization 
problem is solved:
\begin{equation}
\min_{\substack{\tilde{A} \in \mathbb{R}^{r \times r},\\ \tilde{B} \in \mathbb{R}^{r \times r \times r}}}{\sum_{j=1}^M || \boldsymbol{\tau_u}^{\text{exact}}(t_j)-\boldsymbol{\tau_u}^{\text{ansatz}}(t_j)||_{L^2(\Omega)}^2} ,
\label{opt_problem}
\end{equation}
where $M$ time instances are considered to build the correction term, and the term $\boldsymbol{\tau}^{\text{exact}}(t_j)$ is computed 
from the snapshot vectors $\mathbf{a}_d^{\text{snap}}(t_j)$, which satisfy 
the conditions 
\[a_{d_i}^{\text{snap}}(t_j)=\left(\mathbf{u}_d(t_j), \boldsymbol{\phi}_i\right) \quad \forall i=1,...,d.\]
The exact 
correction term is evaluated as follows:
\[ \boldsymbol{\tau_u}^{\text{exact}}(t_j)=\left(- \overline{(\mathbf{a}_d^{\text{snap}}(t_j))^T\mathbf{C_d}\mathbf{a}_d^{\text{snap}}(t_j)}^r\right) -\left(-(\mathbf{a}_r^{\text{snap}}(t_j))^T\mathbf{C}\mathbf{a}_r^{\text{snap}}(t_j) \right) , \]
where the tensor $\mathbf{C_d} \in \mathbb{R}^{d \times d \times d}$ is defined in the following way:
\[
\mathbf{C_d}_{ijk}=\left(\boldsymbol{\phi}_i, \nabla \cdot (\boldsymbol{\phi}_j \otimes \boldsymbol{\phi}_k) \right).
\]
At each time step $j$, the approximated correction term is evaluated as in \eqref{ansatz1}, but starting from $\mathbf{a}_r^{\text{snap}}(t_j)$: 
\begin{equation}
\boldsymbol{\tau_u}^{\text{ansatz}}(t_j)=\tilde{A}\mathbf{a}_r^{\text{snap}}(t_j)+(\mathbf{a}_r^{\text{snap}}(t_j))^T\tilde{B}\mathbf{a}_r^{\text{snap}}(t_j).
\end{equation}
The optimization problem \eqref{opt_problem} is rewritten as a least squares problem following a procedure similar to that used in \cite{peherstorfer2016data}. In particular, the following terms are defined:
\begin{itemize}
    \item the snapshot matrix $\hat{X} \in \mathbb{R}^{M \times r}$. 
    Denoting with $\hat{X}_{j,\cdot}$ the $j$-th row of the matrix, we have 
    \begin{equation}
    \hat{X}_{j,\cdot}=\mathbf{a}_r^{\text{snap}}(t_j);\end{equation}
    \item $r$ vectors $\mathbf{a}^{(1)}(t_j)$,...,$\mathbf{a}^{(r)}(t_j)$ for each time step, such that 
    \begin{equation}\mathbf{a}^{(i)}(t_j)=a^{\text{snap}}_i(t_j) \begin{bmatrix}a_0^{\text{snap}}(t_j)\\a_1^{\text{snap}}(t_j)\\...\\a_i^{\text{snap}}(t_j)
    \end{bmatrix} \in \mathbb{R}^i \text{ for }i=1,...,r,\end{equation} where $a^{\text{snap}}_i(t_j)$ is the $i$-th component of the snapshot vector at time step $j$;
    \item $r$ different matrices $\hat{X}^{(1)}$,...,$\hat{X}^{(r)}$, with $\hat{X}^{(i)}\in \mathbb{R}^{M \times i}$ such that 
    \begin{equation}
    \hat{X}^{(i)}_{j,\cdot}=\mathbf{a}^{(i)}(t_j) ;
    \end{equation}
    \item the matrix $\mathbf{R} \in \mathbb{R}^{M \times r}$ such that 
    \begin{equation}
        \mathbf{R}_{j,\cdot}=\mathbf{\tau_u}^{\text{exact}}(t_j) \quad \forall j=1,..., M  .
    \end{equation}
\end{itemize}
The optimization problem \eqref{opt_problem} can be expressed in the following way:
\begin{equation}
\min_{\substack{\tilde{A} \in \mathbb{R}^{r \times r},\\ \tilde{B} \in \mathbb{R}^{r \times r \times r} }}{||\mathbf{R} - \hat{X} \tilde{A}^T - \sum_{i=1}^r \hat{X}^{(i)} (\tilde{B}^{(i)})^T||_F^2} ,
\label{new_opt_problem}
\end{equation}
where $\tilde{B}^{(i)}$ are blocks of the tensor $\tilde{B}$ of dimension $i \times i$, and the norm considered in the minimization is the 
Frobenius 
norm.
\RA{We now define
}
\[\mathbf{D}=[\hat{X}, \hat{X}^{(1)}, \hat{X}^{(2)},..., \hat{X}^{(r)}], \quad \mathbf{O}=[\tilde{A},\tilde{B}^{(1)}, \tilde{B}^{(2)},...,\tilde{B}^{(r)}] .\]
\RA{We 
denote by $\ell=r+(1+2+\dots+r)$ 
the global number of columns of matrices $\mathbf{D}$ and $\mathbf{O}$ defined above.}
In a more compact form, the optimization problem \eqref{new_opt_problem} can be written as follows:
\begin{equation}
    \min_{\mathbf{O} \RA{\in \mathbb{R}^{r \times \ell}}}{||\mathbf{R}-\mathbf{D}\mathbf{O}^T||^2_F} ,
    \label{final_opt_problem}
\end{equation}
Problem \eqref{final_opt_problem} can be also seen as a set of $r$ optimization problems 
\begin{equation}
\min_{\mathbf{o}_i, i=1,...,r}{||\mathbf{r}_i-\mathbf{D}\mathbf{o}_i||^2_{L^2(\Omega)}} ,
\end{equation}
where $\mathbf{o}_i$ is the $i$-th row of matrix $\mathbf{O}$ and $\mathbf{r}_i$ is the $i$-th column of the matrix $\mathbf{R}$.
As in \cite{xie2018data} and \cite{peherstorfer2016data}, the problem \eqref{final_opt_problem} is ill-conditioned since the matrix $\mathbf{D}$ has a large condition number. In order to solve the least squares problem, a truncated singular value decomposition is applied to matrix $\mathbf{D}$, just as in step 6 of Algorithm 1 in \cite{xie2018data}. 

After $\tilde{A}$ and $\tilde{B}$ are found from the least squares problem, the dynamical system to be solved is 
\begin{equation}
\begin{cases}
\mathbf{M} \dot{\mathbf{a}}=\nu(\mathbf{B}+\mathbf{B_T})\mathbf{a}-\mathbf{a}^T \mathbf{C} \mathbf{a}-\mathbf{H}\mathbf{b}+\tau \left( \sum_{k=1}^{N_{\text{BC}}}(U_{BC,k}\mathbf{D}^k-\mathbf{E}^k \mathbf{a})\right)+\tilde{A}\mathbf{a}+\mathbf{a}^T \tilde{B} \mathbf{a} ,\\
    \mathbf{P}\mathbf{a}=\mathbf{0} .
\end{cases}
\label{system_corr_U}
\end{equation}
The ill-conditioning of the least squares problem leads to an ill-conditioning of the dynamical system \eqref{system_corr_U}. In order to 
mitigate this ill-conditioning, step 5 of Algorithm 1 of 
\cite{mou2021data} is applied: the number $R$ of singular values retained 
in matrix $\mathbf{D}$ is the optimal one, i.e., the 
number that minimizes the error metric 
\begin{equation}
    \varepsilon_u (L^2)= \sum_{j=1}^M ||\mathbf{u}_{\text{sol}}(t_j)-\mathbf{u}_r(t_j)||_{L^2(\Omega)} ,
    \label{error_metric}
\end{equation}
where, at each time step, $\mathbf{u}_{\text{sol}}(t_j)=\sum_{i=1}^r a_i(t_j) \boldsymbol{\phi}_i$ is found from the solution of the dynamical system \eqref{system_corr_U}. 

\RA{\remark{
We note that the correction term 
depends on the POD basis considered. Indeed, the 
correction term depends on the operators $\mathbf{C}$ and $\mathbf{C}_d$, which depend on the POD basis. 
Furthermore, the correction term depends on the time window used to perform the POD, and on the time window used to build the correction term. We 
denote these two time intervals as $[0,T_{\text{offline}}]$ and $[0,T_{\text{correction}}]$, respectively.}}

\subsection{Constrained data-driven correction for velocity}
\label{constrainedvel}
To increase the ROM accuracy, physical constraints can be added to the least squares problem for velocity~\cite{mohebujjaman2019physically}.
In this approach, the matrices $\tilde{A}$ and $\tilde{B}$ 
are endowed with the following physical properties:
\begin{itemize}
    \item $\mathbf{a}^T \tilde{A} \mathbf{a} \leq 0$, i.e., $\tilde{A}$ is negative semi-definite;
    \item $\mathbf{a}^T(\mathbf{a}^T \tilde{B} \mathbf{a})=0$, i.e., $\tilde{B}$ is skew-symmetric.
\end{itemize}
The optimization problem becomes 
\begin{equation}
    \min_{\substack{\tilde{A} \in \mathbb{R}^{r \times r},\\ \tilde{B} \in \mathbb{R}^{r \times r \times r},\\ \mathbf{a} \tilde{A} \mathbf{a} \leq 0,\\\mathbf{a}^T(\mathbf{a}^T \tilde{B} \mathbf{a})=0}}{\sum_{j=1}^M || \boldsymbol{\tau}^{\text{exact}}(t_j)-\boldsymbol{\tau}^{\text{ansatz}}(t_j)||_{L^2(\Omega)}^2} .
    \label{opt_problem2}
\end{equation}
As shown in \cite{mohebujjaman2019physically} and 
and confirmed in the numerical investigation in section~\ref{results}, 
the constrained method is more accurate than the unconstrained one when the number of modes 
used for the velocity and pressure 
yields a \emph{marginally-resolved} regime. 
However, as the number of modes increases, the unconstrained method seems to produce 
more accurate results. 

\RA{\remark{As extensively explained in \cite{mohebujjaman2019physically}, the constraints imposed on $\tilde{A}$ and $\tilde{B}$ aim 
at reproducing the constraints of the operators 
to which our correction is applied. In 
\cite{mohebujjaman2019physically}, the correction term $\mathbf{a}^T \mathbf{\tilde{B}} \mathbf{a}$ acts on the nonlinear term $\mathbf{a}^T\mathbf{C}\mathbf{a}$}, so $\tilde{B}$ inherits the skew-symmetric property of $\mathbf{C}$. On the other hand, the term $\tilde{A} \mathbf{a}$ is treated in \cite{mohebujjaman2019physically} as a separate 
correction 
to the diffusive term $\nu(\mathbf{B}+\mathbf{B}_T)\mathbf{a}$, and inherits the property of that operator. Since our test case has a 
high Reynolds number, the diffusive term 
does not have a significant impact on the system, and thus we do not consider it in the exact correction $\boldsymbol{\tau_u}^{\text{exact}}$}.

\subsection{Data-driven corrections for pressure}
\label{presSUP}
This section introduces the novel 
correction terms for pressure in the reduced equations of the supremizer approach. The system \eqref{new_reduced_system} is considered and now the focus is the evaluation of the closure expressions for $\boldsymbol{\tau}^{p(1)}$ and $\boldsymbol{\tau}^{p(2)}$.
From the expressions \eqref{tau_all} and following a procedure similar to that 
used in section \ref{velocity}, the ansatzes' expressions for the corrective terms are written as 
\begin{equation}
\boldsymbol{\tau}^{p(1)}(\mathbf{b})=\tilde{H} \mathbf{b},\quad
\boldsymbol{\tau}^{p(2)}(\mathbf{a})=\tilde{P} \mathbf{a}.
\end{equation}
We consider two 
approaches to find the matrices $\tilde{H}$ and $\tilde{P}$.
\begin{itemize}
    \item[1.] Solve two different optimization problems, one for each 
    correction term:
\begin{equation}
\min_{\tilde{H} \in \mathbb{R}^{r \times q}} ||\boldsymbol{\tau_{p(1)}}^{\text{exact}}-\boldsymbol{\tau_{p(1)}}^{\text{ansatz}}||_{L^2(\Omega)}^2,\quad
\min_{\tilde{P} \in \mathbb{R}^{q \times r}} ||\boldsymbol{\tau_{p(2)}}^{\text{exact}}-\boldsymbol{\tau_{p(2)}}^{\text{ansatz}}||_{L^2(\Omega)}^2,
        \label{opt_P_SUP}
\end{equation}
        where 
        \begin{equation}
        \begin{split} 
        &\boldsymbol{\tau_{p(1)}}^{\text{exact}}(t_j)
        = \left( -\overline{\mathbf{H_d} \mathbf{b}_{d_p}^{\text{snap}}(t_j)}^r\right)-\left(- (\mathbf{H} \mathbf{b}_{q}^{\text{snap}}(t_j)) \right),\quad
        \boldsymbol{\tau_{p(1)}}^{\text{ansatz}}
        =\tilde{H} \mathbf{b}_{q}^{\text{snap}}(t_j),\\
        &\boldsymbol{\tau_{p(2)}}^{\text{exact}}(t_j)=\left( \overline{\mathbf{P_d} \mathbf{a}_{d}^{\text{snap}}(t_j)}^r\right)-\left( \mathbf{P} \mathbf{a}_{r}^{\text{snap}}(t_j) \right), \quad
        \boldsymbol{\tau_{p(2)}}^{\text{ansatz}}=\tilde{P} \mathbf{a}_{r}^{\text{snap}}(t_j).
        \end{split}
        \end{equation}
    \item[2.] Solve a unique optimization problem to find both correction terms:
    \begin{equation}
    \min_{\substack{\tilde{H}\in \mathbb{R}^{r \times q},\\\tilde{P}\in \mathbb{R}^{q \times r}}}{||\boldsymbol{\tau_{p}}^{\text{exact}}-\boldsymbol{\tau_{p}}^{\text{ansatz}}||_{L^2(\Omega)}^2},
    \label{opt_P_SUP_2}
    \end{equation} 
    where \RA{we consider the following compact notation to vertically stack the vectors:}
    \begin{equation}
        \boldsymbol{\tau_{p}}^{\text{exact}}=[\boldsymbol{\tau_{p(1)}}^{\text{exact}}, \boldsymbol{\tau_{p(2)}}^{\text{exact}}], \quad \boldsymbol{\tau_{p}}^{\text{ansatz}}=[\boldsymbol{\tau_{p(1)}}^{\text{ansatz}}, \boldsymbol{\tau_{p(2)}}^{\text{ansatz}}].
    \label{opt_p2}
    \end{equation}
  
\end{itemize}

In both cases, the matrices 
in the least squares problems are ill-conditioned and the truncated singular value decomposition 
is applied to 
mitigate this issue.  
As in section \ref{velocity}, the number of singular values retained in each optimization problem is chosen in order to minimize the error metric $\varepsilon_p(L^2)$, defined as 
\begin{equation}
    \varepsilon_p (L^2)= \sum_{j=1}^M ||p_{\text{sol}}(t_j)-p_{q}(t_j)||_{L^2(\Omega)},
    \label{error_metricP}
\end{equation}
where $p_{\text{sol}}(t_j)=\sum_{i=1}^{q} b_i(t_j) \chi_i$ is found from the solution of the dynamical system at each time step, and $p_q(t_j)$ is the projection of the full order pressure on the space generated by the first $q$ modes. The reason for choosing a different error metric involving the pressure field is that, in this section, the pressure corrections are 
introduced to improve the accuracy of the pressure field.

The dynamical system obtained 
by adding the new correction terms is the following: 
\begin{equation}
\begin{cases}
\mathbf{M} \dot{\mathbf{a}}=\nu(\mathbf{B}+\mathbf{B_T})\mathbf{a}-\mathbf{a}^T \mathbf{C} \mathbf{a}-\mathbf{H}\mathbf{b}+\tilde{H}\mathbf{b}+\tau \left( \sum_{k=1}^{N_{\text{BC}}}(U_{BC,k}\mathbf{D}^k-\mathbf{E}^k \mathbf{a})\right)+\tilde{A}\mathbf{a}+\mathbf{a}^T \tilde{B} \mathbf{a},\\
    \mathbf{P}\mathbf{a}+\tilde{P}\mathbf{a}=\mathbf{0}.
\end{cases}
\label{system_corr_U_P}
\end{equation}
%
In the construction of the 
correction terms for pressure and velocity, the supremizer modes are not 
used since they induce numerical instability. 
We believe that the reason for this numerical instabilty is that the supremizer modes are fictitious modes without a physical significance, introduced to fulfill the 
inf-sup condition. Since the aim of the 
correction terms is to model the contribution of the 
neglected physical modes, 
retaining a large number of supremizer modes in the correction terms can adversely affect the 
data-driven procedure.
Therefore, in the least squares problems \eqref{opt_problem}, \eqref{opt_problem2}, \eqref{opt_P_SUP}, and \eqref{opt_P_SUP_2}, the number of velocity modes 
used is $r=N_u$ instead of $r=N_u+N_{sup}$.

%% file: sections_paper1/section_datappe.tex
\section{Data-driven corrections for the PPE-ROM approach}
\label{datappe}
In this section, the data-driven techniques used in 
section \ref{datasup} for the supremizer approach are adapted to the PPE-ROM approach. This approach 
constructs new correction 
terms 
that model the contribution of the neglected modes inside the Poisson equation, as showed in \eqref{filtered_NSEPPE}: 
\begin{equation}
\begin{cases}
\left( -\frac{\partial \bar{\mathbf{u}}_d^r}{\partial t},\boldsymbol{\phi}_i \right) +\nu \left( \nabla \cdot \left( \nabla \bar{\mathbf{u}}_d^r + (\nabla \bar{\mathbf{u}}_d^r)^T \right), \boldsymbol{\phi}_i \right)-\left( (\bar{\mathbf{u}}_d^r \cdot \nabla) \bar{\mathbf{u}}_d^r,\boldsymbol{\phi}_i \right)-(\nabla\bar{p}_d^r,\boldsymbol{\phi}_i)+\\
\quad \quad+\left( \boldsymbol{\tau}_u^{SFS},\boldsymbol{\phi}_i\right)+\left( \boldsymbol{\tau}_{p(1)}^{SFS},\boldsymbol{\phi}_i\right) =0 \ \ \text{ for } i=1,...,r,\\
\left(\nabla  \bar{p}_d^r,\nabla \chi_i \right) + \left(\nabla  \cdot (\bar{\mathbf{u}}_d^r \otimes \bar{\mathbf{u}}_d^r),\nabla \chi_i \right) - \nu \left( \nabla \times \bar{\mathbf{u}}_d^r, \mathbf{n} \times \nabla \chi_i \right)_{\Gamma} - \left(\mathbf{n} \cdot \mathbf{R}_{d_t}^r, \chi_i \right)_{\Gamma}+\\
\quad \quad +\left( \tau_{D}^{SFS}, \nabla \chi_i\right) + \left( \tau_{G}^{SFS}, \nabla \chi_i\right) =0 \ \ \text{ for } i=1,...,q.
\end{cases}
\label{filtered_NSEPPE}
\end{equation}
As pointed out in section \ref{presSUP}, the number of velocity modes used to build the correction terms is $r=N_u$.
The new exact closure terms introduced in system \eqref{filtered_NSEPPE} are 
\begin{equation}
    \boldsymbol{\tau}_{D}^{SFS}= \overline{ \nabla p_d}^r - \nabla p_d^r, \quad
    \boldsymbol{\tau}_{G}^{SFS}= \overline{\nabla \cdot (\mathbf{u}_d \otimes \mathbf{u}_d)}^r - \nabla \cdot(\mathbf{u}_d^r \otimes \mathbf{u}_d^r).
\end{equation}
Next, the following terms are introduced:
\begin{equation}
\boldsymbol{\tau}^D \RA{\in \mathbb{R}^{q}} \text{ such that }\tau^D_i= \left( \boldsymbol{\tau}_D^{SFS}, \nabla \chi_i\right) ,\quad
\boldsymbol{\tau}^{G} \RA{\in \mathbb{R}^{q}} \text{ such that }\tau^{G}_i= \left( \boldsymbol{\tau}_{G}^{SFS},\nabla \chi_i\right).
\label{tau_allPPE}
\end{equation}
Adding the new correction terms, the dynamical system \eqref{reduced_systemPPE} becomes:
\begin{equation}
\begin{cases}
\mathbf{M} \dot{\mathbf{a}}=\nu(\mathbf{B}+\mathbf{B_T})\mathbf{a}-\mathbf{a}^T \mathbf{C} \mathbf{a}-\mathbf{H}\mathbf{b}+\tau \left( \sum_{k=1}^{N_{\text{BC}}}(U_{BC,k}\mathbf{D}^k-\mathbf{E}^k \mathbf{a})\right)+\boldsymbol{\tau}^u+\boldsymbol{\tau}^{p(1)},\\
    \mathbf{D}\mathbf{b}+\mathbf{a}^T \mathbf{G} \mathbf{a} - \nu \mathbf{N} \mathbf{a} - \mathbf{L}+\boldsymbol{\tau}^{D}+\boldsymbol{\tau}^{G}=\mathbf{0}.
\end{cases}
\label{new_reduced_systemPPE}
\end{equation}

In order to close the system \eqref{new_reduced_systemPPE}, the 
following approximations
need to be found: 
$\boldsymbol{\tau}^u \approx \boldsymbol{\tau}^u (\mathbf{a}) $,  $\boldsymbol{\tau}^{p(1)}\approx \boldsymbol{\tau}^{p(1)} (\mathbf{b})$,  $\boldsymbol{\tau}^{D}\approx \boldsymbol{\tau}^{D} (\mathbf{b})$, $\boldsymbol{\tau}^{G}\approx \boldsymbol{\tau}^{G} (\mathbf{a},\mathbf{b})$. 

\subsection{Data-driven correction for the term $\mathbf{D} \mathbf{b}$ in \eqref{new_reduced_systemPPE}}
\label{Dcorr_ref}
The first proposed correction term 
is that related to the matrix $\mathbf{D}$.
For the term $\boldsymbol{\tau}^D$, the following ansatzes 
are proposed and 
tested:
    \begin{enumerate}
        \item A linear ansatz,
        $\boldsymbol{\tau}^{D} (\mathbf{b})=\tilde{D} \mathbf{b}$.
        
        The procedure followed to evaluate matrix $\tilde{D}$ is similar to 
        that used to find $\tilde{H}$ or $\tilde{P}$ in section \ref{presSUP}. The detailed implementation is described in \cite{ivagnes2021data}.
        After $\tilde{D}$ is found from the least squares problem, the Poisson equation in the dynamical system is 
\begin{equation}
\mathbf{D}\mathbf{b} + \mathbf{a}^T \mathbf{G} \mathbf{a} + \tilde{D} \mathbf{b} -\nu \mathbf{N} \mathbf{a} -\mathbf{L}=\mathbf{0}  .
\label{system_PPEcorr_D}
\end{equation}

As in the previous cases, the 
ill-conditioning of the least squares problem is 
mitigated by applying a truncated singular value decomposition in which the error metric \eqref{error_metricP} is 
minimized.
We note that the 
error metric to be minimized for the correction 
terms in the Poisson equation concerns the pressure field. The reason for this choice is that the correction term $\tilde{A} \mathbf{a} +\mathbf{a}^T \tilde{B} \mathbf{a}$ 
aims at improving the results mainly for the velocity field. However, it is 
numerically shown that the correction terms added in the Poisson equation have no effect on the velocity field, but they significantly improve the pressure results.
        
        \item A quadratic ansatz, 
        $\boldsymbol{\tau}^{D} (\mathbf{b})=\tilde{D} \mathbf{b} + \mathbf{b}^T\tilde{B}_P \mathbf{b}$, where $\tilde{B}_P$ is a 
        three-dimensional tensor. 
        This choice is motivated by the velocity correction term in section \ref{velocity}.
        
        The Poisson equation of the dynamical system in this case becomes 
\begin{equation}
\mathbf{D}\mathbf{b} + \mathbf{a}^T \mathbf{G} \mathbf{a} + \tilde{D} \mathbf{b} + \mathbf{b}^T \tilde{B}_P \mathbf{b} -\nu \mathbf{N} \mathbf{a} -\mathbf{L}=\mathbf{0}.
\label{system_PPEcorr_D2}
\end{equation}
In 
section~\ref{results}, the quadratic ansatz $2$ is chosen to approximate the exact correction term, since it provides the best results among the 
proposed choices. 
    \end{enumerate}
\subsection{Data-driven correction for the term $\mathbf{a}^T \mathbf{G} \mathbf{a}$ in \eqref{new_reduced_systemPPE}}
\label{Gterm}
The second proposed correction term 
is that depending on the matrix $\mathbf{G}$. The ansatz 
for the term $\boldsymbol{\tau}^G$ is the following:
\[\boldsymbol{\tau}^{G} (\mathbf{a})=\tilde{G}_A \mathbf{a}+ \mathbf{a}^T \tilde{G}_B \mathbf{a}.\]
This ansatz is similar to that proposed to approximate the velocity correction in the SUP-ROM approach in section \ref{velocity}.
The procedure followed to find the matrices $\tilde{G}_A$ and $\tilde{G}_B$ is exactly the same 
as that used
in section \ref{velocity}, but considering the following exact term:
\begin{equation}
    \boldsymbol{\tau}^{\text{exact}}_G = \overline{(\mathbf{a}_d^{\text{snap}}(t_j))^T \mathbf{G}_d \mathbf{a}_d^{\text{snap}}(t_j)}^r - (\mathbf{a}_r^{\text{snap}}(t_j))^T \mathbf{G} \mathbf{a}_r^{\text{snap}}(t_j),
\end{equation}
where the tensor $\mathbf{G}_d \in \mathbb{R}^{d \times d \times d}$ is defined as 
\[\mathbf{G}_{d_{ijk}}= (\nabla \chi_i, \nabla \cdot (\boldsymbol{\phi_j} \otimes \boldsymbol{\phi}_k)) , 
\ \text{ with } i,j,k=1,...,d.\]
The Poisson equation 
using just the correction for the term containing the tensor $\mathbf{G}$ can be rewritten as 
\begin{equation}
    \mathbf{D}\mathbf{b} + \mathbf{a}^T \mathbf{G} \mathbf{a} + \tilde{G}_A \mathbf{a} + \mathbf{a}^T \mathbf{G}_B \mathbf{a} -\nu \mathbf{N} \mathbf{a} -\mathbf{L}=\mathbf{0}.
\label{system_PPEcorr_G}
\end{equation}

\subsection{Combined data-driven corrections proposals}
\label{proposals}
We 
investigate other approaches
in order to provide a more compact ansatz in the reduced order systems. In particular, we 
approximate more than one correction term with a unique ansatz. The following ansatzes 
are proposed and investigated:
\begin{enumerate}
    \item \label{jointDG} Joint correction to $\mathbf{D} \mathbf{b} + \mathbf{a}^T \mathbf{G} \mathbf{a}$: $\tau^{\text{joint}}_{DG} (\mathbf{a}, \mathbf{b})=\tilde{D}_{pg} \mathbf{b}+ \mathbf{a}^T \tilde{B}_{pg}\mathbf{a}$. In this case, the following least squares problem 
    is solved:
\begin{equation}
\min_{\substack{\tilde{D}_{pg} \in \mathbb{R}^{q \times q}; \\ \tilde{B}_{pg} \in \mathbb{R}^{q \times r \times r}}}{\sum_{j=1}^M || \boldsymbol{\tau}_{DG}^{\text{exact}}(t_j)-\boldsymbol{\tau}_{DG}^{\text{ansatz}}(t_j)||_{L^2(\Omega)}^2},
\label{opt_problemDG}
\end{equation}
where the exact correction term is 
\begin{equation}
    \boldsymbol{\tau}_{DG}^{\text{exact}}(t_j) = \boldsymbol{\tau}_D^{\text{exact}}(t_j)+ \boldsymbol{\tau}_G^{\text{exact}}(t_j)  \ \ \forall j=1,...,M. 
\end{equation}
    \item \label{jointDGab} Joint correction to $\mathbf{D} \mathbf{b} + \mathbf{a}^T \mathbf{G} \mathbf{a}$ as a function of the vector of coefficients $\mathbf{ab}=(\mathbf{a}, \mathbf{b}) \in \mathbb{R}^{N_u+N_p}$: $\tau^{\text{joint,  ab}}_{DG} (\mathbf{ab})= \tilde{I}_A \mathbf{ab} + \mathbf{ab}^T \tilde{I}_B \mathbf{ab}$.
    
    Since the term $\mathbf{D} \mathbf{b}$ depends just on the pressure modes, whereas the term $\mathbf{a}^T \mathbf{G} \mathbf{a}$ depends on both the velocity and the pressure modes, the two terms can be merged into a unique least squares problem involving the total vector of coefficients $\mathbf{ab}=(\mathbf{a}, \mathbf{b})$.
Specifically, denoting $r_{\text{tot}}=r+q$, the least squares problem is 
\begin{equation}
    \min_{\substack{\tilde{I}_A \in \mathbb{R}^{q \times r_{\text{tot}}}; \\ \tilde{I}_B \in \mathbb{R}^{q \times r_{\text{tot}} \times r_{\text{tot}}}}}{\sum_{j=1}^M || \boldsymbol{\tau}_{DG}^{\text{exact}}(t_j)-\boldsymbol{\tau}_{joint, ab}^{\text{ansatz}}(t_j)||_{L^2(\Omega)}^2}. 
\label{opt_problem_ab}
\end{equation}
    \item \label{jointDGCab} Joint correction to $\mathbf{D} \mathbf{b} + \mathbf{a}^T \mathbf{G} \mathbf{a}$ in the Poisson equation and $\mathbf{a}^T \mathbf{C} \mathbf{a}$ in the momentum equation: $\tau^{joint,  ab}_{DCG} (\mathbf{ab})= \tilde{J}_A \mathbf{ab} + \mathbf{ab}^T \tilde{J}_B \mathbf{ab}$.
    
    In 
    this case, a single least squares problem is solved in order to find three 
    correction terms. The optimization problem is 
\begin{equation}
    \min_{\substack{\tilde{J}_A \in \mathbb{R}^{r_{\text{tot}} \times r_{\text{tot}}}; \\ \tilde{J}_B \in \mathbb{R}^{r_{\text{tot}} \times r_{\text{tot}} \times r_{\text{tot}}}}}{\sum_{j=1}^M || \boldsymbol{\tau}_{DCG, ab }^{\text{exact}}(t_j)-\boldsymbol{\tau}_{DCG, ab}^{\text{ansatz}}(t_j)||_{L^2(\Omega)}^2}.
\label{opt_problem_abC}
\end{equation}
The exact term is 
\begin{equation}
    \boldsymbol{\tau}_{DCG,  ab}^{\text{exact}}(t_j) = \left(\boldsymbol{\tau}_u^{\text{exact}}(t_j),\boldsymbol{\tau}_D^{\text{exact}}(t_j)+ \boldsymbol{\tau}_G^{\text{exact}}(t_j)\right)  \ \ \forall j=1,...,M.
\end{equation}
The matrices $\tilde{J}_A \in \mathbb{R}^{r_{\text{tot}} \times r_{\text{tot}}}$ and $\tilde{J}_B \in \mathbb{R}^{r_{\text{tot}} \times r_{\text{tot}}\times r_{\text{tot}}}$ are computed through a 
procedure similar to that used in case 2. The final correction is divided into two vectors:
\[\tilde{J}_A \mathbf{ab} + \mathbf{ab}^T \tilde{J}_B \mathbf{ab}=\left( \mathbf{J}_1, \mathbf{J}_2\right), \ \text{ where } \mathbf{J}_1 \in \mathbb{R}^{N_u}, \mathbf{J}_2 \in \mathbb{R}^{N_p}.\]
The dynamical system with the novel data-driven correction terms becomes 
\begin{equation}
\begin{cases}
\mathbf{M} \dot{\mathbf{a}}=\nu(\mathbf{B}+\mathbf{B_T})\mathbf{a}-\mathbf{a}^T \mathbf{C} \mathbf{a}-\mathbf{H}\mathbf{b}+\tau \left( \sum_{k=1}^{N_{\text{BC}}}(U_{BC,k}\mathbf{D}^k-\mathbf{E}^k \mathbf{a})\right)+ \mathbf{J}_1,\\
    \mathbf{D}\mathbf{b}+\mathbf{a}^T \mathbf{G} \mathbf{a} -\nu \mathbf{N} \mathbf{a}-\mathbf{L}+ \mathbf{J}_2=\mathbf{0}.
\end{cases}
\label{system_corr_all_PPE}
\end{equation}
\end{enumerate}
Details on all the compact cases 
are provided in~\cite{ivagnes2021data}.

%% file: sections_paper1/section_numresults.tex
\section{Numerical Results}
\label{results}
This section presents the numerical results obtained with the 
ROMs described in section \ref{sec:g-rom} combined with the data-driven techniques introduced in sections \ref{datasup} and \ref{datappe}. The test case considered is 
the unsteady turbulent flow around a circular cylinder, which is a typical case study in the field of fluid dynamics.

\begin{figure}[h!]
\centering
\subfloat[]{\includegraphics[width=0.6 \textwidth]{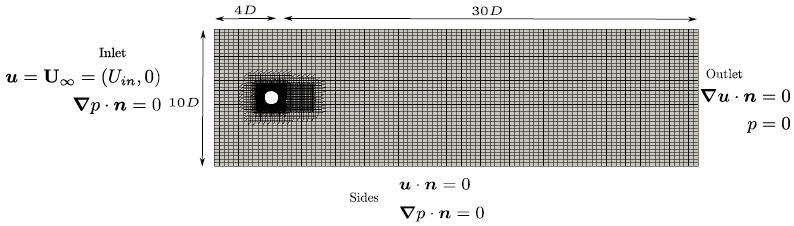}\label{1}}
\subfloat[]{\includegraphics[ width=0.4 \textwidth]{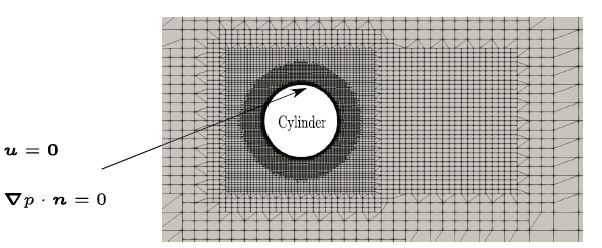}\label{2}}
\caption{(a) The mesh used in simulations.
(b) The mesh zoomed in around the cylinder. Image 
adapted from \cite{hijazi2020data}.}
\label{mesh}
\end{figure}

Given the two-dimensional nature of the vortex shedding phenomenon, the numerical test case is set up in two dimensions and the mesh used is composed 
of 11644 polygonal cells. The mesh and 
boundary conditions set for the velocity and pressure are presented in Figure \ref{mesh} \cite{hijazi2020data}. The diameter of the cylinder is $D= \SI{1}{\metre}$, the fluid kinematic viscosity is $\nu=\SI{1e-4}{\metre \squared \per \second}$, and the velocity at the inlet is horizontal and fixed at $U_{in}=\SI{5}{\metre \per \second}$. 
These parameters yield a Reynolds number Re=$\SI{5e4}{}$.

The open source software \emph{OpenFOAM} is used in the offline stage to 
generate the high-fidelity fields. In particular, we 
use the 
unsteady solver \emph{pimpleFoam} and 
the $\kappa-\omega$ model for pressure-velocity coupling and turbulence treatment, respectively.
\RA{
In the offline stage, the high-fidelity snapshots are collected 
every $0.004$ s.}

The POD algorithm is then e\RA{m}ployed to obtain the POD modes for the velocity and  pressure from the respective snapshot matrices, along with the supremizer modes. 
The following subsections will 
discuss the results obtained from simulations. Simulations are computed using both ITHACA-FV \cite{stabile2018finite,stabile2017pod} and in-house developed \emph{Python} scripts. Both the 
SUP-ROM and PPE-ROM approaches are considered.

In all the simulations reported, the $d$ value used for the construction of the closure terms is $d=50$.
The results of the FOM simulations are compared 
with the results obtained by solving the reduced order dynamical systems with or without the data-driven terms and considering the parameter $\tau=1000$ in~\eqref{system_corr_U_P} and in ~\eqref{new_reduced_systemPPE}.

\RA{
The SUP-ROM dynamical system is written in a more general form:}
\begin{equation}
\begin{cases}
\mathbf{M} \dot{\mathbf{a}}^{i}=\nu(\mathbf{B}+\mathbf{B_T})\mathbf{a}^i-(\mathbf{a}^{i})^T \mathbf{C} \mathbf{a}^i-\mathbf{H}\mathbf{b}^i\RA{+\tau \left( \sum_{k=1}^{N_{\text{BC}}}(U_{BC,k}\mathbf{D}^k-\mathbf{E}^k \mathbf{a}^{i})\right)}+  c_{u} \boldsymbol{\tau}_{u}(\mathbf{a}^i, \mathbf{b}^i)&+  \\ + c_{p(1)} \boldsymbol{\tau}_{p(1)}(\mathbf{a}^i, \mathbf{b}^i)  &\text{ at each }i=1,...,M ,\\
\mathbf{P}\mathbf{a}^i + c_{p(2)} \boldsymbol{\tau}_{p(2)}(\mathbf{a}^i, \mathbf{b}^i)=\mathbf{0} &\text{ at each }i=1,...,M ,
\end{cases}
\label{supgen}
\end{equation}
where $M$ is the total number of time steps in the online phase.
The matrices and tensors in \eqref{supgen} are defined in section \ref{sup}.

In the PPE-ROM framework, the reduced system supplemented with
the \RA{correction} terms is the following:
\begin{equation}
\begin{cases}
\mathbf{M} \dot{\mathbf{a}}^{i}=\nu(\mathbf{B}+\mathbf{B_T})\mathbf{a}^i-(\mathbf{a}^{i})^T \mathbf{C} \mathbf{a}^i-\mathbf{H}\mathbf{b}^i \RA{+\tau \left( \sum_{k=1}^{N_{\text{BC}}}(U_{BC,k}\mathbf{D}^k-\mathbf{E}^k \mathbf{a}^{i})\right)}+ c_u \boldsymbol{\tau}_u(\mathbf{a}^i, \mathbf{b}^i) &\text{ at each }i=1,...,M ,\\
\mathbf{D}\mathbf{b}^i+ (\mathbf{a}^i)^T \mathbf{G} \mathbf{a}^i - \nu \mathbf{N} \mathbf{a}^i- \mathbf{L}+c_D \boldsymbol{\tau}_{D}(\mathbf{a}^i, \mathbf{b}^i) +c_G \boldsymbol{\tau}_{G}(\mathbf{a}^i, \mathbf{b}^i)=\mathbf{0} &\text{ at each }i=1,...,M  .
\end{cases}
\label{ppegen}
\end{equation}
The matrices and tensors in \eqref{ppegen} are defined in section \ref{ppe_sec}.
\RA{In systems \eqref{supgen} and \eqref{ppegen}, the additional parameters $c_u, c_{p(i)}, c_D$, and $c_G$ 
are set to $1$ if the corresponding correction term is included in the system, and $0$ otherwise.
}

To better guide the reader through all the numerical tests carried out and their specific purpose, a brief summary of their evolution 
in this section is provided below.
\begin{enumerate}[label=\textnormal{(\arabic*)}]
\item \label{p1} A necessary preliminary step is that of analyzing the results of the POD modal decomposition. An analysis of the eigenvalue decay for the pressure and velocity is carried out in order 
to identify the \emph{under-resolved} and \emph{marginally-resolved} regimes. 
\item \label{p2} The solution of the standard SUP-ROM approach \RA{(system \eqref{supgen}, with $c_u=c_{p(i)}=0$, $i=1,2$)} is then studied. In particular, the stability issues associated with the number of supremizer modes 
are discussed.
\item \label{p3} The data-driven methods proposed in \cite{xie2018data,mohebujjaman2019physically,mou2020data} are then extended to the SUP-ROM approach. 
The influence of the term $\boldsymbol{\tau}_{u}(\mathbf{a}^i)$, which we refer to as correction term for velocity, is analyzed. As will be shown, the dynamical system \RA{\eqref{supgen}} with $c_u=1$ and $c_{p(1)}=c_{p(2)}=0$ produces a better approximation of the velocity field than the formulation without any correction term. However, the approximation of the pressure field is not significantly improved with respect to the solution of the system without the velocity correction term.
\item \label{p4} Since in important applications it is necessary to obtain an accurate pressure field prediction, a further test 
considers new pressure corrections for the SUP-ROM approach. In particular, the dynamical system \RA{\eqref{supgen}} with $c_u=c_{p(1)}=c_{p(2)}=1$ is solved making use of both the velocity and the pressure correction/closure terms. As will be 
shown, the results still 
display no meaningful improvement of the pressure field. A possible reason for this 
behavior is that the correction terms are mainly designed to model effects of nonlinear terms
, as is the case with its velocity counterpart. Since no nonlinear terms involve pressure, the corresponding correction/closure term is less effective than the velocity correction. Another possible reason for the inefficiency of the pressure corrections is that the supremizer model does not  consider a dedicated equation for the pressure field, therefore no corrections can directly affect the pressure approximation.

\item \label{p5} To investigate effective strategies for the reduced pressure field prediction, the next step is to consider a different approach where a dedicated equation for pressure is employed. 
To this end, we resort to the PPE-ROM. In the PPE formulation, novel correction terms are added to the pressure Poisson equation and different ansatzes are proposed. The results lead to an evident improvement of the reduced pressure field.
\RA{\item \label{p6} 
Steps \ref{p1}-\ref{p5} of our analysis serve as a \emph{validation} of the methods and ansatzes proposed in this paper, always considering the same time windows for both the POD construction ($T_{\text{offline}}=\SI{20}{\second}$) and the correction terms ($T_{\text{correction}}=\SI{2}{\second}$).

In this final 
step, we \emph{test} the methods considering different time windows, which are selected according to the Strouhal number, i.e., the frequency of vortex shedding. In our case, the Strouhal number is defined as follows:
\[
St = \dfrac{fD}{U\infty}=\frac{D}{T U\infty}=0.224.
\]
}
\end{enumerate}

For clarity of presentation, the models tested in this section are summarized in Table \ref{table1}.

Finally, two different time integration schemes 
are considered in the numerical tests. The time derivative $\dot{\mathbf{a}}^i$ appearing in the momentum equation in both formulations \eqref{supgen} and \eqref{ppegen} is either computed by making use of an implicit Euler time 
discretization or with an implicit second-order time 
discretization. It is worth remarking that the second-order time 
discretization corresponds to the scheme implemented in \emph{OpenFOAM} and used to solve the full order problem. For the standard ROM approaches, the results obtained with both \RA{discretizations are compared in order to assess if the 
FOM-ROM 
consistency of the second-order time 
discretization provides increased accuracy. 
For clarity of presentation, for the data-driven techniques, our analysis will focus only on the first-order time 
discretization}.
The quantitative evaluation of the ROM 
accuracy is obtained by means of the 
relative $L^2(\Omega)$ error of velocity and pressure. The errors are evaluated with respect to the full order fields and typically compared with the reconstruction errors, i.e., 
the errors yielded by projecting the full order fields on a subspace generated by a given number of reduced modes. The projection of the full order solution is in fact the best possible result which can be achieved with a given 
number of modes. Thus, the solution of the reduced system cannot generally lead to an error lower than the projection error. 
The 
relative errors with respect to the full order fields at the $j$-th time step 
are defined as
\begin{equation}
    \varepsilon_u\RA{^{full}}(t_j)=\dfrac{||\mathbf{u}_r^{\text{abs}}(\mathbf{x}, t_j)-\mathbf{u}_d^{\text{abs}} (\mathbf{x}, t_j)||_{L^2(\Omega)}}{||\mathbf{u}^{\text{abs}}_d (\mathbf{x}, t_j)||_{L^2(\Omega)}} , \quad \varepsilon_p\RA{^{full}}(t_j)=\dfrac{||p_r(\mathbf{x}, t_j)-p_d(\mathbf{x}, t_j)||_{L^2(\Omega)}}{||p_d(\mathbf{x}, t_j)||_{L^2(\Omega)}} ,
    \label{errors1}
\end{equation}
where we 
used the following quantities:
\begin{itemize}
    \item the reduced order fields of velocity 
    and pressure, namely
    $$\mathbf{u}_r(\mathbf{x}, t_j)=\sum_{i=1}^{r} a_i(t_j) \boldsymbol{\phi}_i(\mathbf{x}) , \quad p_r(\mathbf{x},t_j)=\sum_{i=1}^{q} b_i(t_j) \chi_i(\mathbf{x})  ,$$
    where the coefficients $a_i(t_j)$ and $b_i(t_j)$ are the solutions of the dynamical systems \eqref{supgen} (in the supremizer approach) and \eqref{ppegen} (in the Poisson approach);
    \item the approximated full order field of velocity 
    and pressure, 
    which are evaluated starting from the first $d$ modes, where $d=100$ when a supremizer approach is considered, and $d=50$ when a pressure Poisson approach is considered. 
    These terms are written as follows:
    $$\mathbf{u}_d(\mathbf{x}, t_j)=\sum_{i=1}^{d} a^{\text{snap}}_i(t_j) \boldsymbol{\phi}_i (\mathbf{x}) , \quad p_d(\mathbf{x}, t_j)=\sum_{i=1}^{d_p} b^{\text{snap}}_i (t_j) \chi_i(\mathbf{x})  .$$
    We remark that when the supremizer approach is considered,  $\{\boldsymbol{\phi}_i\}_{i=51}^{100}=\{\mathbf{s}_i(\chi_i)\}_{i=1}^{50}$ are the supremizer modes.
\end{itemize}

\begin{table}[h!]\centering

\begin{tabular}{ p{0.25 \textwidth}p{0.2 \textwidth}p{0.45 \textwidth}  }
 \toprule[0.3ex]
 \textbf{Model} & \textbf{Corresponding section} & \textbf{Reduced order system} \\
 \midrule
 Standard SUP-ROM  & \ref{sup_enrich_nocorr}  & System \eqref{supgen} with $c_u = c_{p(1)} = c_{p(2)} = 0$ \\
 \midrule
 SUP-ROM with velocity correction  & \ref{velresults}  & System \eqref{supgen} with $c_u = 1$, $c_{p(1)} = c_{p(2)} = 0$\\
  \midrule
 SUP-ROM with pressure corrections &   \ref{presSUP_2} & System \eqref{supgen} with $c_u = 0$, considering the cases:
 \[
 \begin{cases}
 c_{p(1)} = c_{p(2)} = 1\\
 c_{p(1)}=1, c_{p(2)} = 0\\
 c_{p(1)}=0, c_{p(2)} = 1\\
 \end{cases}
 \] \\
 \hline
 Standard PPE-ROM  & \ref{PPE_standard_results}  & System \eqref{ppegen} with $c_u = c_D = c_G = 0$ \\
  \midrule
 PPE-ROM with pressure corrections & \ref{pres_comb_sec} & System \eqref{ppegen} with $c_u = 0$, considering the cases:
 \[
 \begin{cases}
 c_{D} = c_{G} = 1\\
 c_{D}=1, c_{G} = 0\\
 c_{D}=0, c_{G} = 1\\
 \end{cases}
 \] \\
 \midrule
 PPE-ROM with velocity and pressure corrections   & \ref{vel_pres_combined} & System \eqref{ppegen} with $c_u = c_D = c_G = 1$\\
 \bottomrule[0.3ex]
 \end{tabular}
 \caption{Summary of the models presented in the numerical results \RA{considering the time windows corresponding to $T_{\text{offline}}=\SI{20}{\second}$, $T_{\text{correction}}=\SI{2}{\second}$, and $T_{\text{online}}=\SI{2}{\second}$ or $\SI{8}{\second}$}.}
\label{table1}
\end{table}
\subsection{Eigenvalue decay for velocity, pressure and supremizer}

\begin{figure}[h!]
    \centering
    \includegraphics[scale=0.7]{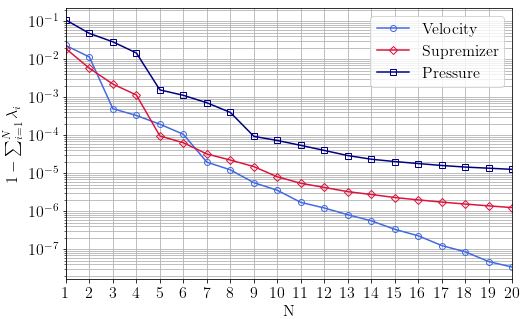}
    \caption{Cumulative ignored normalized eigenvalue 
    decay, for velocity, supremizer and pressure eigenvalues.}
    \label{fig:eigen}
\end{figure}
In Figure \ref{fig:eigen}, the normalized eigenvalue decay is 
displayed for the velocity, supremizer, and pressure modes. The plot shows that a relatively small number of modes is sufficient in order to retain most of the energetic information in the snapshots. The plot suggests that the marginally-resolved regime --- in which our numerical investigation will take place --- extends between 3 to 8 modes. 
However, it must be pointed out that the energetic interpretation for the eigenvalue decay of the supremizer modes is not as straightforward as that of the velocity and pressure modes. 
In particular, the eigenvalues of the supremizer modes cannot be interpreted as an energetic contribution since 
these modes are fictitious nonphysical modes added to fulfill a numerical stability condition. Thus, the number of supremizer modes included in the online simulations is not chosen based on energetic consideration. Instead, we impose that $N_{sup} \geq N_p$ in order to avoid stability issues \cite{stabile2018finite}. 

\subsection{Analysis of the SUP-ROM without corrections}
\label{sup_enrich_nocorr}
This section is dedicated to the preliminary analysis of the solutions of the dynamical system 
without data-driven corrections.

\begin{figure}[htpb!]
\subfloat[Velocity error ($1^{st}$-order time scheme) ]{\includegraphics[ scale=0.42]{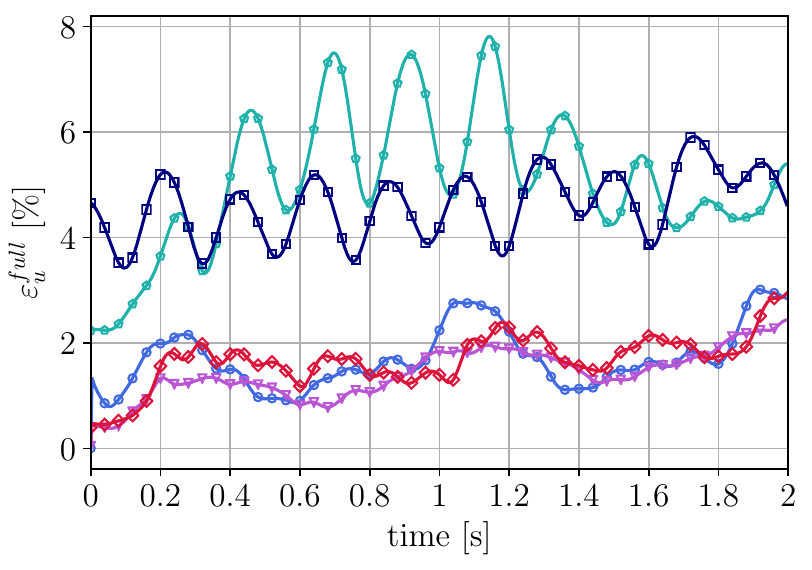}\label{vel_1_instab}}
\subfloat[Pressure error ($1^{st}$-order time scheme) ]{\includegraphics[ scale=0.42, trim=0 0 6.5cm 0]{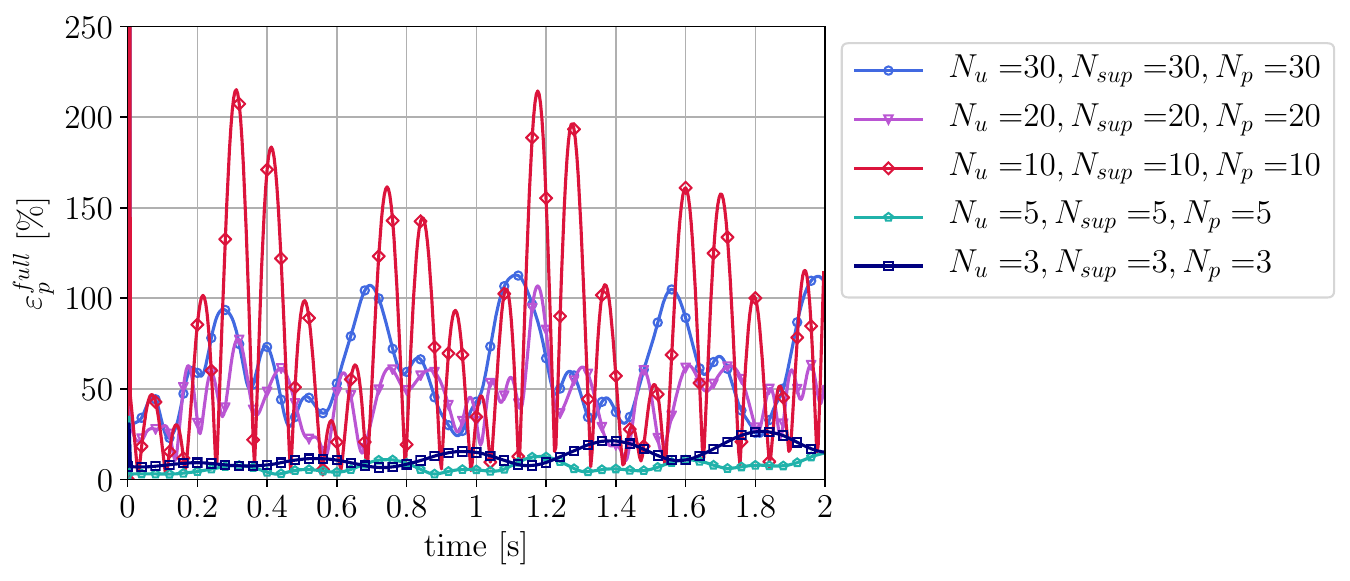}\label{pres_1_instab}}\\
\subfloat[Velocity error ($1^{st}$-order time scheme) ]{\includegraphics[ scale=0.42]{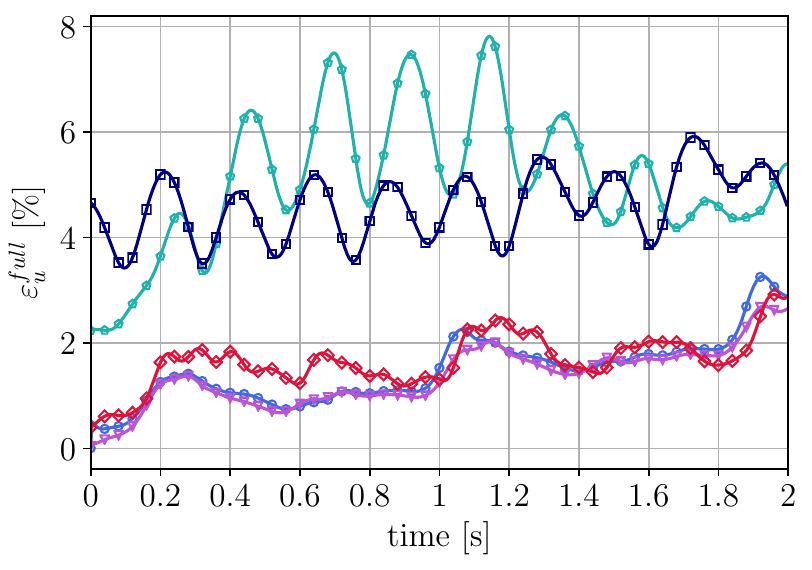}\label{vel_1_stab}}
\subfloat[Pressure error ($1^{st}$-order time scheme) ]{\includegraphics[ scale=0.42, trim=0 0 6.5cm 0]{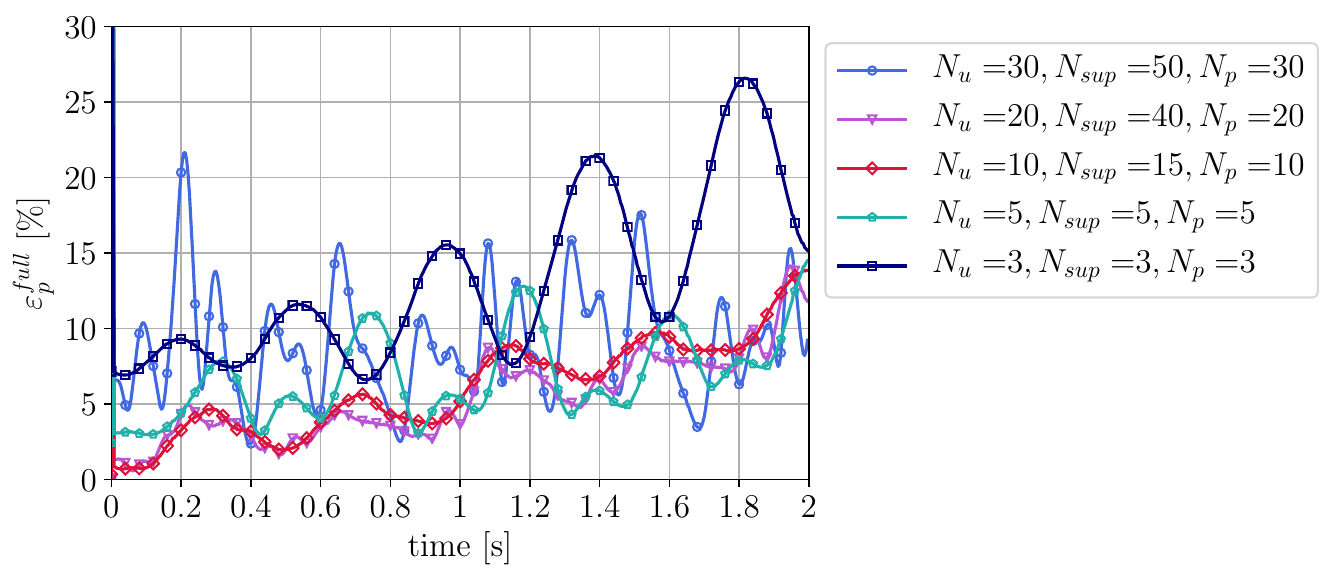}\label{pres_1_stab}}\\
\subfloat[Velocity error ($2^{nd}$-order time scheme) ]{\includegraphics[ scale=0.42]{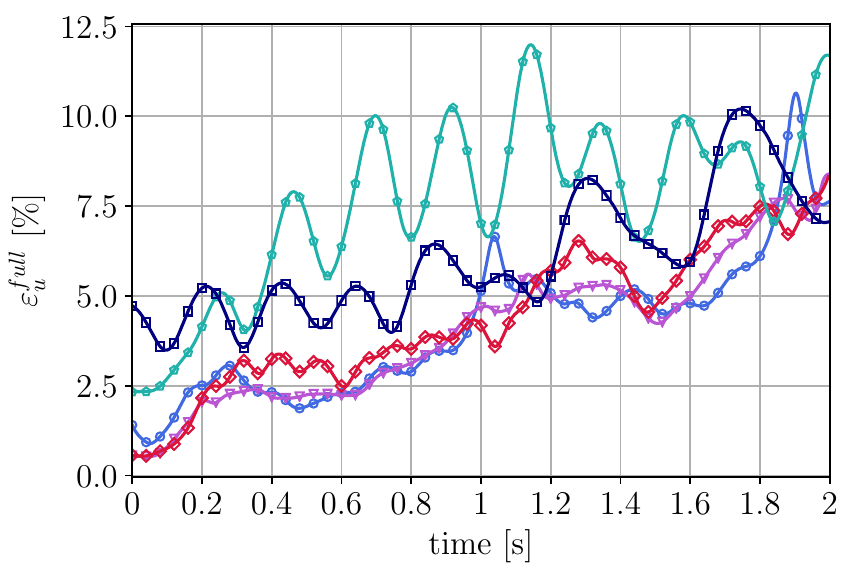}\label{vel_2_instab}}
\subfloat[Pressure error ($2^{nd}$-order time scheme) ]{\includegraphics[scale=0.42, trim=0 0 6.5cm 0]{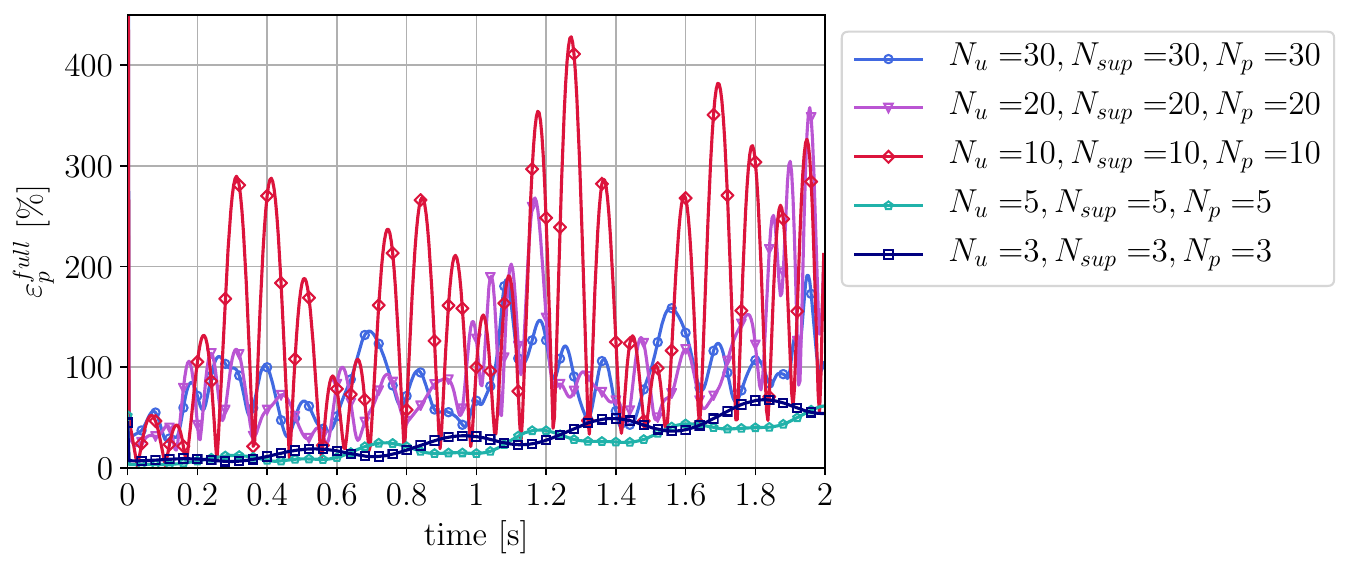}\label{pres_2_instab}}\\
\subfloat[Velocity error ($2^{nd}$-order time scheme) ]{\includegraphics[ scale=0.42]{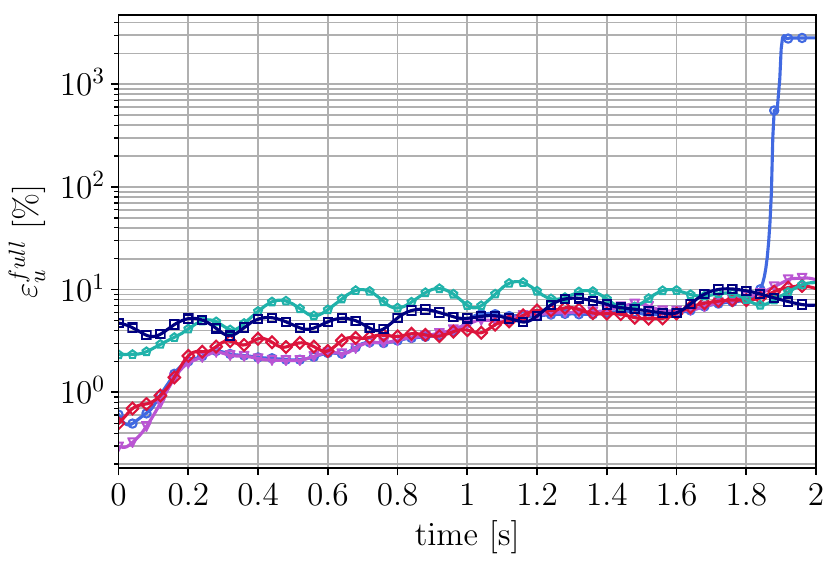}\label{vel_2_stab}}
\subfloat[Pressure error ($2^{nd}$-order time scheme) ]{\includegraphics[ scale=0.42, trim=0 0 6.5cm 0]{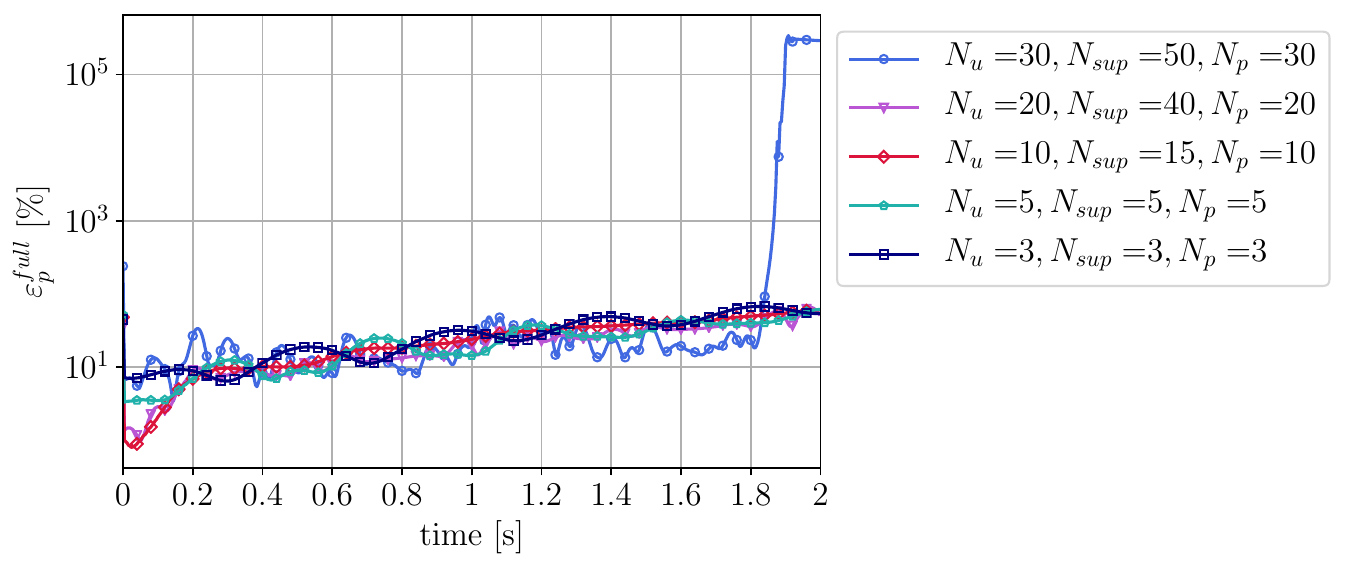}\label{pres_2_stab}}
\caption{
Relative errors represented at different time steps, when solving the reduced system with a first-order  (\protect \subref{vel_1_instab}, \protect \subref{pres_1_instab}, \protect \subref{vel_1_stab}, and \protect \subref{pres_1_stab}) or second-order time discretization 
(\protect \subref{vel_2_instab}, \protect \subref{pres_2_instab},\protect \subref{vel_2_stab} and \protect \subref{pres_2_stab}). The first and third rows of plots 
correspond to the case with a number of supremizer modes $N_{sup}=N_p=N_u$, whereas the second and fourth rows display results with the addition of extra supremizer modes.}
\label{modes_all}
\end{figure}

Plots \ref{modes_all}\protect \subref{vel_1_instab}, \protect \subref{pres_1_instab}, \protect \subref{vel_1_stab}, and \protect \subref{pres_1_stab} display the time evolution of the 
relative errors for reduced velocity and pressure fields. 
The diagrams are obtained making use, at the reduced order level, of a first-order time 
discretization scheme. Different combinations of truncation orders $N_u, N_{p}, N_{sup}$ for velocity, pressure, and supremizer modes, respectively, are considered; the expectation is that the accuracy of the standard Galerkin-ROM 
improves as the number of modes is increased. However, the \emph{approximated} supremizer approach is characterized by stability issues particularly affecting the pressure field, as specified in section \ref{sup}. 
As can be seen in Figure \ref{modes_all}\protect \subref{pres_1_instab}, when we are in the resolved regime and the number of supremizer modes is $N_{sup}=N_u=N_p$, the reduced pressure solution error significantly increases. These issues are not as severe when $N_{sup}>N_p$, as shown in Figure 
\ref{pres_1_stab}, where results obtained in the totally-resolved regime show an acceptable accuracy, except for the case $N_u=N_p=30$ and $N_{sup}=50$, where stability issues are again observed. Again, such stability problems are likely associated with the \emph{approximated} supremizer procedure adopted, and might be solved resorting to an \emph{exact} supremizer procedure as suggested in \cite{ballarin2015supremizer,stabile2018finite}. In this paper, we chose to adopt an approximated supremizer enrichment procedure, since it allows to reduce the online computational cost.
%
%
The corresponding results obtained when considering a second-order time discretization are presented in Figure \ref{modes_all}\protect \subref{vel_2_instab}, \protect \subref{pres_2_instab},\protect \subref{vel_2_stab}, and \protect \subref{pres_2_stab}. The plots suggest that the second-order time discretization suffers from the same stability issue observed 
for the first-order time discretization. 

The comparison between Figures \ref{modes_all}\protect \subref{vel_1_stab}, \protect \subref{pres_1_stab} and \ref{modes_all}\protect \subref{vel_2_stab}, \protect \subref{pres_2_stab} shows that when a first-order time discretization 
is considered the velocity field results are similar to the ones obtained with a second-order time discretization, 
except for 
the error blow up occurring just at the end of simulation for the case $N_u=N_p=30, N_{sup}=50$. 
This fact is likely due to the reduced numerical dissipation associated with the second-order time discretization. 
Such lower dissipation might in fact make the system more exposed to the instabilities associated with the approximated supremizer pressure treatment.
\medskip

The following sections 
will focus on the error analysis in the marginally-resolved regime, which is not affected by stability issues. However, the plots in Figure \ref{modes_all} also suggest that the standard SUP-ROM consistently leads to velocity field predictions that are significantly more accurate than their pressure counterparts.; the main goal of the following sections will be 
to evaluate strategies to improve the pressure accuracy.

\subsection{Effect of velocity correction in the SUP-ROM approach}
\label{velresults}
This section is dedicated to the analysis of the effect of the velocity correction terms discussed in section \ref{velocity}
on the solution of  the supremizer ROM dynamical system. 

The data-driven term is introduced in the supremizer formulation following the two approaches --- unconstrained and constrained --- described in section \ref{velocity}. In section \ref{constrained_and_not}, 
these alternative formulations are compared in terms of 
accuracy. In addition,  
for both formulations, we address the problem of choosing the optimal number of singular values $R$ and $R_c$ retained in the least squares problem appearing in the data driven 
optimization of the coefficients.
In section \ref{vel_res}, the prediction efficiency of the method is tested by building the correction terms from snapshots 
in a time interval that is smaller than the interval
used for the resolution at the ROM level.

\subsubsection{Velocity correction: constrained and unconstrained cases}
\label{constrained_and_not}

This section 
aims at making a comparison between the constrained and unconstrained corrections for the term $\boldsymbol{\tau}_u$. The parameters appearing in system \eqref{supgen} are set as follows: $c_u=1, c_{p(1)}=c_{p(2)}=0$. 
The number of modes 
is fixed at $N_u=N_p=N_{sup}=5$. Thus, we are in the marginally-resolved modal regime.

\begin{figure}[h!]
\subfloat[Relative error on velocity]{\includegraphics[ scale=0.5]{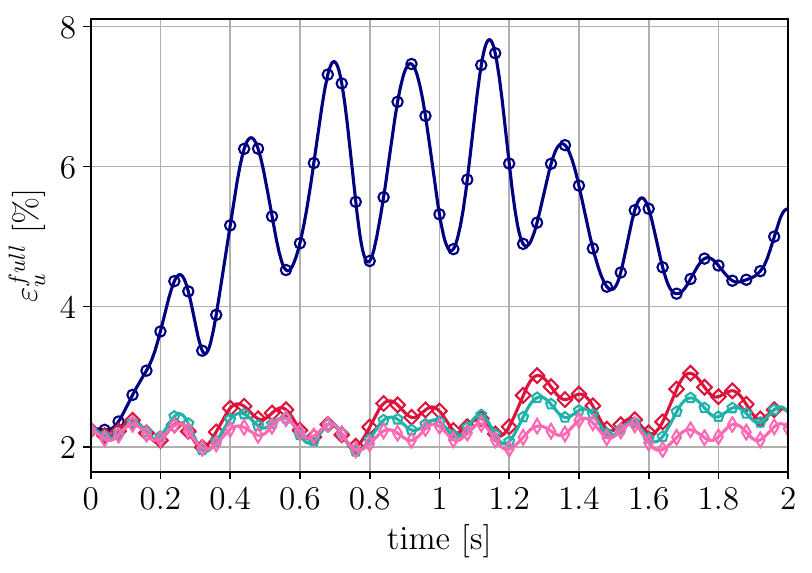}}
\subfloat[Relative error on pressure]{\includegraphics[ scale=0.5]{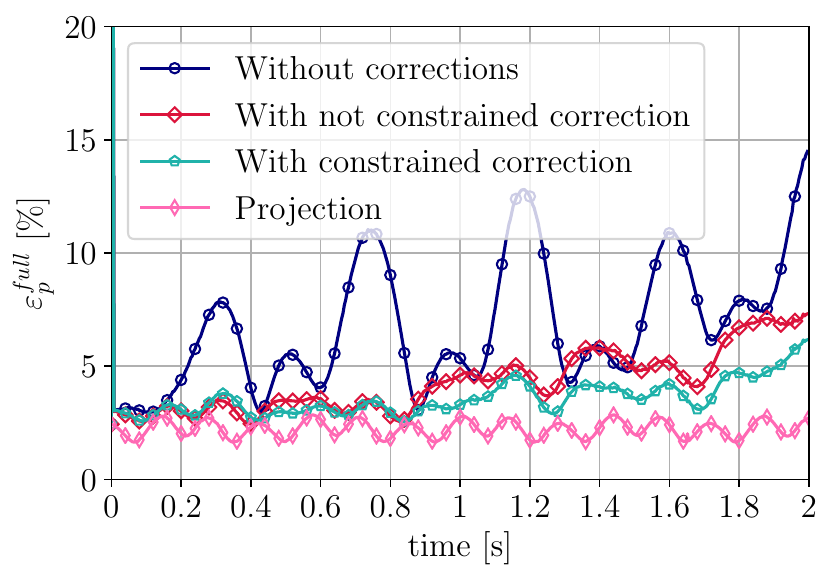}}
\caption{
Relative errors of the absolute value of velocity (a) and 
pressure (b) with respect to the full order simulations, considering $N_u=N_p=N_{sup}=5$. Results without any correction term, and with the unconstrained and constrained corrections, are displayed.}
\label{const_comp2}
\end{figure}

The plots in Figure \ref{const_comp2} clearly show that the velocity correction, both in its constrained and 
unconstrained variant, improves the approximation of both the velocity and the pressure field with respect to the standard SUP-ROM. 
Quite remarkably, the addition of the velocity correction is able to bring the velocity error to values that approach the  projection error. The reduction of the pressure error, albeit 
significant, is not as pronounced.

We recall that the constrained correction is derived by including physical conditions in the optimization problem, which are supposed to positively influence the velocity and pressure fields. However, the results in our numerical investigation suggest that the accuracy gain associated with the constrained method appears marginal. Moreover, the addition of constraints to the method is not sufficient to obtain pressure field error comparable with the projection error.

\noindent In order to further lower the pressure error, new correction terms are introduced and evaluated in the numerical simulations of the following sections.

\begin{figure}[h!]
    \centering
    \includegraphics[scale=0.55]{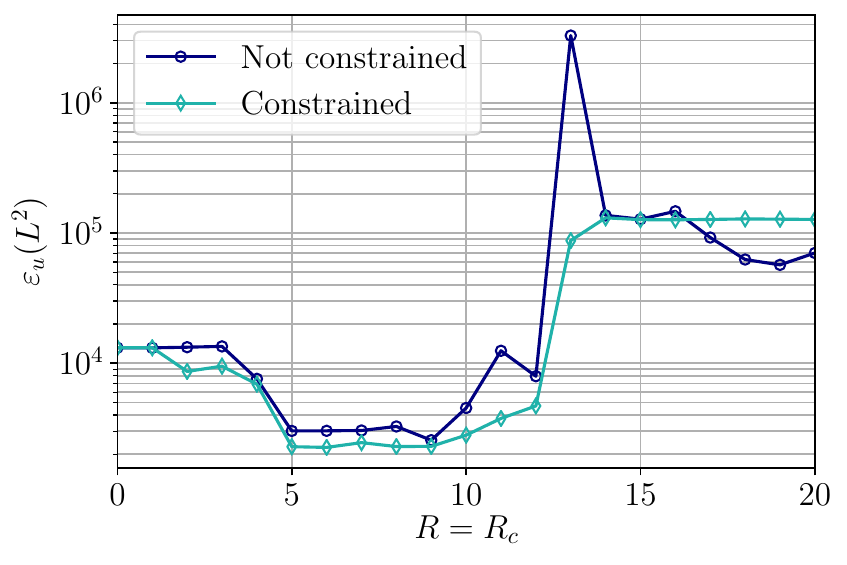}
    \caption{
    The metric $\varepsilon_u(L^2)$ 
    for different numbers of singular values retained in the singular value decomposition.}
    \label{metric_varyingR}
\end{figure}

We now assess how the accuracy of the solution for system \eqref{supgen} with $c_u=1$ and $c_{p(1)}=c_{p(2)}=0$ is affected by the number of singular values retained in the least square problem appearing in the optimization process 
used to obtain the data driven correction coefficients. In particular, we want to establish whether it is possible to identify optimal values of 
the truncation orders 
that can be conveniently used in the rest of the numerical 
investigation.
In the least squares problems \eqref{opt_problem} and \eqref{opt_problem2}, the number of singular values retained in the truncated SVD is optimized with respect to the error metric expressed in \eqref{error_metric}. The trend of the error metric for different values of $R$ and $R_c$ is 
displayed in Figure \ref{metric_varyingR}. For both the unconstrained and constrained cases the error metric has a similar trend with $R$ or $R_c$. Increasing values of $R$ improves the solution of the reduced system with respect to the standard case until a minimum in the error metric is reached; too large values of $R$ lead to a divergence of the reduced solution with respect to the high-fidelity result. In all the numerical simulations presented in the article, the optimal value of singular values is 
used to build the correction terms. 

Another important observation is that the selected $R$ value optimizes the error of the velocity field \RA{with respect to the projection of the full order field. Hence, $R$ is chosen \emph{a posteriori} such that the velocity solution of the dynamical system is as accurate as possible.} Thus, \RA{it does not optimize the error with respect to the exact correction and} it can provide better results than the exact correction. Using the terminology in section VI.A of \cite{ahmed2021closures}, 
\emph{model regression} is utilized to determine the model form of the correction terms, whereas 
\emph{trajectory regression} is used to determine the optimal parameters in these model forms.

\subsubsection{Velocity correction: 
predictive accuracy}
\label{vel_res}
This section evaluates the performance of the method proposed in sections \ref{velocity} and \ref{constrainedvel} when larger time windows are considered and time prediction is carried out. In particular, the matrices $\tilde{A}$ and $\tilde{B}$ are built using 
snapshots taken from a smaller interval that does not cover all the simulated time, in order to test the \emph{predictive} capability of the proposed method. 
Specifically, the time interval on which we collect the snapshots is $[0, \SI{20}{\second}] $, the time interval used to test the ROM is $[0, \SI{8}{\second}] $, and the 
velocity correction is built starting from data extracted from $[0, \SI{2}{\second}]$.
The modal regime considered is $N_u=N_p=N_{sup}=5$ and the time 
discretization 
is first-order.

\begin{figure}[h!]
\centering
\subfloat[Relative error on velocity]{\includegraphics[ scale=0.45,valign=t]{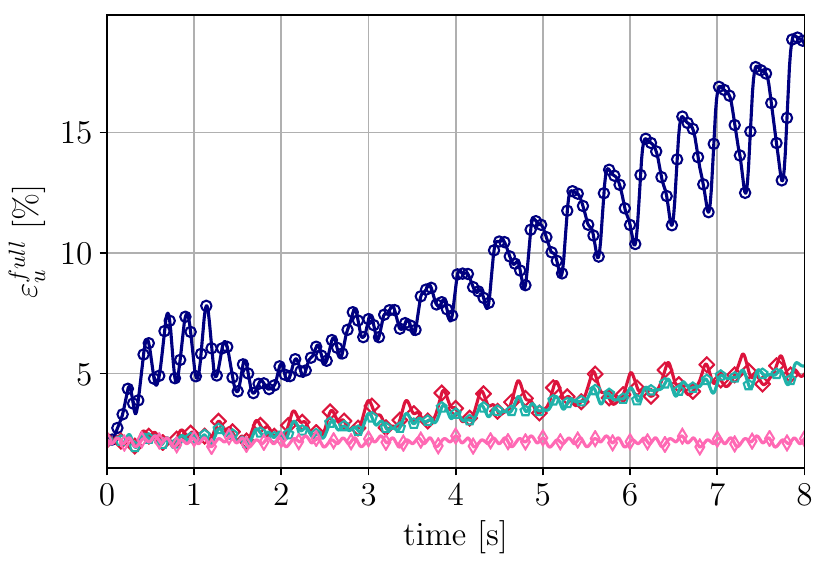}}
\subfloat[Relative error on pressure]{\includegraphics[ scale=0.45,  valign=t]{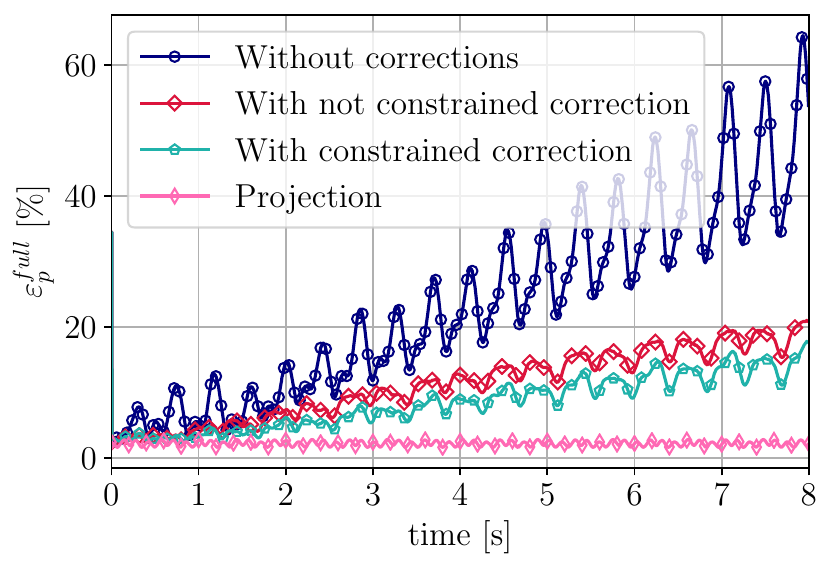}}
\caption{
Relative errors of the 
velocity (a) and 
pressure (b), 
for $N_u=N_p=N_{sup}=5$. Results without any correction term, 
with the constrained velocity correction term, and with the unconstrained velocity correction term.
The optimal values $R$ and $R_c$ are $9$ and $6$, respectively.}
\label{errlong2}
\end{figure}

The 
relative velocity and pressure errors  are displayed in Figure \ref{errlong2}.
The plots confirm that 
adding the correction term improves significantly both the velocity and the pressure field accuracy, especially in the interval $[2,8]$ seconds. This suggests that, provided that the 
training 
time interval (i.e., the time interval from which the snapshots were collected) contains all the significant frequencies characterizing the flow field time evolution, even a limited number of  
snapshots are sufficient to construct the matrices $\tilde{A}$ and $\tilde{B}$ 
in the correction and 
obtain a significant accuracy improvement 
with respect to the case without corrections.

In addition, the 
ROM approximations for 
all the test cases in which corrections have been applied appear to be significantly more stable than the test case without correction. This might once again suggest that if the training set contains a sufficient number of solution cycles, the correction will increase 
the ROM's ability to lead to stable solutions over time, and extended time integration and 
prediction will be possible.
However, in such a stable scenario, constraining the minimization for the calculation of $\tilde{A}$ and $\tilde{B}$ does not seem to lead to significant improvements, as can be seen from Figure \ref{errlong2}.

\subsection{Effect of pressure corrections in the SUP-ROM approach}
\label{presSUP_2}
The results shown in section \ref{velresults} extend to 
the RANS setting the data-driven correction terms tested so far only for LES simulations \cite{xie2018data, mohebujjaman2019physically}. In addition, 
these results show 
that the velocity correction significantly improves the 
accuracy of the velocity field approximations,
but the same cannot be said for the pressure field. 
We note that, to the best of our knowledge, there are no similar data-driven correction models that aim at improving the accuracy of the pressure field.
However, in our numerical investigation, the proposed data-driven corrections for the pressure did not yield completely satisfactory results. 

\begin{figure}[h!]
\centering
\subfloat[Relative error on velocity]{\includegraphics[ scale=0.45,valign=t]{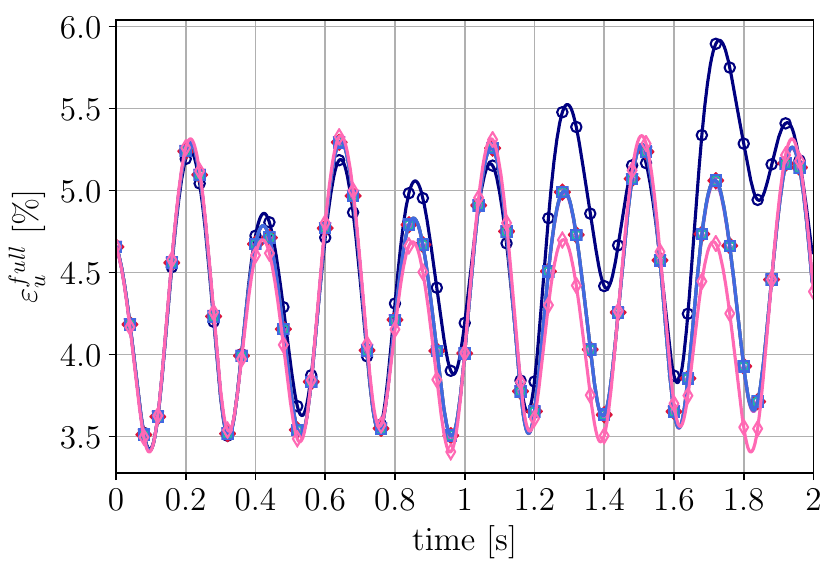}}
\subfloat[Relative error on pressure]{\includegraphics[ scale=0.45,  valign=t]{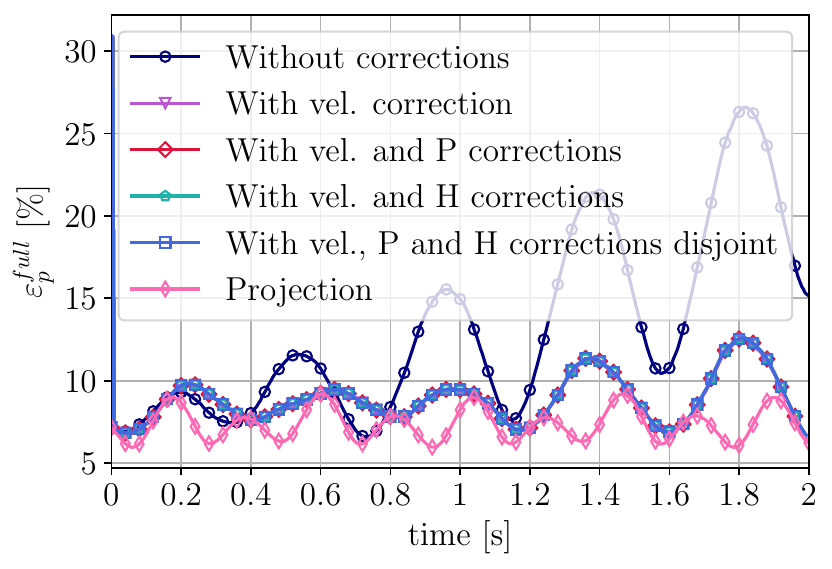}}
\caption{
Relative errors of the 
velocity (a) and 
pressure (b), 
using $N_u=N_p=N_{sup}=3$. Results without any correction term, with only the velocity correction term, 
and with both the velocity and the pressure correction terms
are displayed.}
\label{err1}
\end{figure}

We note that the pressure field is extremely important in 
many applications that require the computation of other flow quantities 
or output values such as forces. Thus, inaccurate pressure fields can lead to poor fluid 
flow force predictions, which can in turn make the proposed techniques 
ineffective in many applications. The present subsection investigates the effectiveness of the new methods developed in section \ref{presSUP} for increased pressure field accuracy.
The following cases for system \eqref{new_reduced_system} are 
used to produce the results in 
Figure \ref{err1}:
\begin{itemize}
    \item $c_u=c_{p(1)}=c_{p(2)}=0$, i.e., absence of any correction term;
    \item $c_u= 0,  c_{p(1)} =1,   c_{p(2)} = 0$, 
    i.e., the correction for the term $-\mathbf{H}\mathbf{b}$ is added in the momentum equation;
    \item $c_u =0,  c_{p(1)} = 0,   c_{p(2)}  =1$, 
    i.e., the correction for the term $\mathbf{P}\mathbf{a}$ is added in the continuity equation;
    \item $c_u = 0,  c_{p(1)} = c_{p(2)} =1$, 
    i.e., the pressure corrections are added to the original system.
\end{itemize}
The plots in Figure \ref{err1} represent the 
relative errors 
for the velocity and pressure fields predicted with each method tested, and 
using $N_u=N_p=N_{sup}=3$ modes for velocity, pressure, and supremizers. 
We 
note that for all the simulations 
using the velocity correction, the constrained data-driven correction method is selected, since it provides slightly better results --- as pointed out in section \ref{constrained_and_not}. In addition, the number of singular values retained for the matrices in the least squares problems is the optimal one. 

\begin{figure}[h!]
        \centering
        \subfloat[
        Relative velocity error]{\includegraphics[ scale=0.45,valign=t]{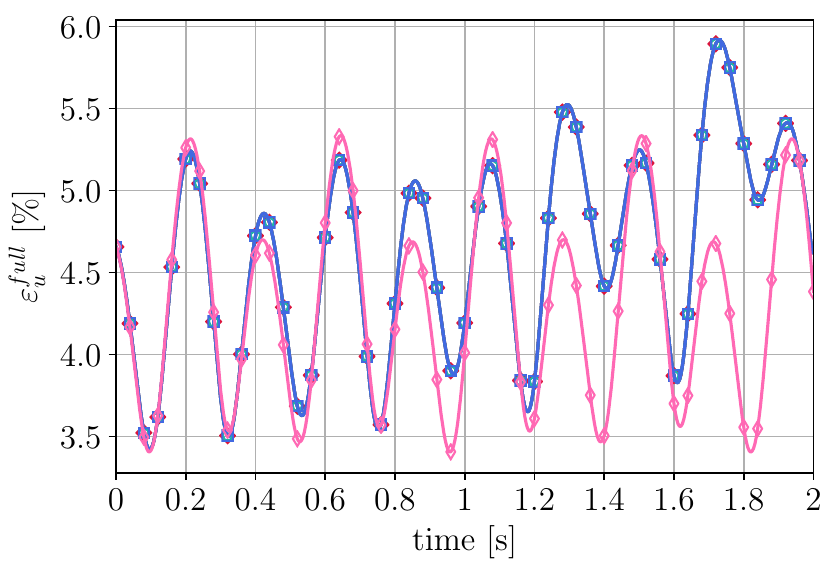}}
\subfloat[
Relative pressure error]{\includegraphics[ scale=0.45,  valign=t]{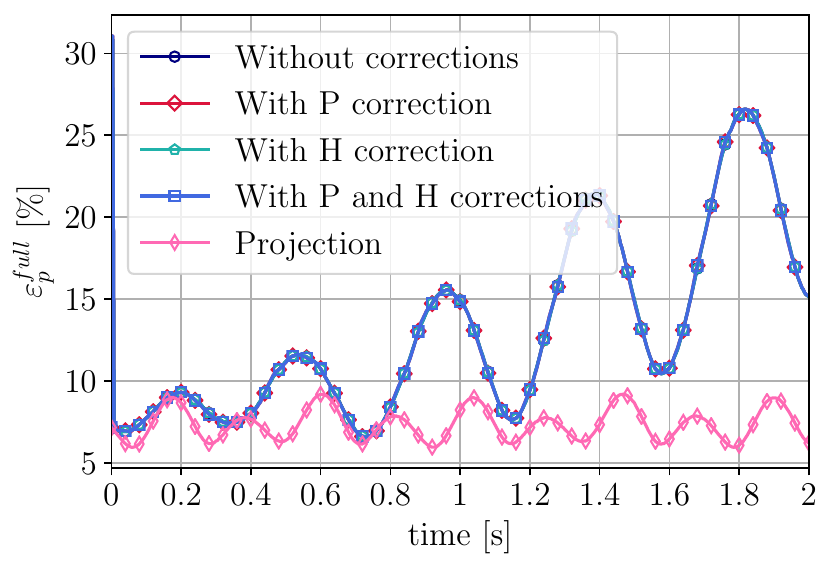}}
\caption{
Relative errors of the 
velocity (a) and 
pressure (b) 
for $N_u=N_p=N_{sup}=3$. Results with and without the exact pressure corrections are presented.}
\label{exactP}
\end{figure}

Figure \ref{err1} shows that the effect of the velocity data-driven correction is much more evident than the effect of the data-driven pressure corrections. 
Indeed, the lines corresponding to the simulations with only the velocity correction term overlap with 
the lines corresponding to the simulations with 
both velocity and pressure correction terms. In other words, once the velocity correction is activated, turning on or off all the pressure corrections 
yields insignificant changes to the ROM error plots.
The ineffectiveness of the pressure 
correction terms 
in the pressure field approximation is further investigated with 
the following tests: To investigate if the poor pressure reconstruction is a 
result of inaccuracies 
in the minimization problems  \eqref{opt_P_SUP} and \eqref{opt_P_SUP_2}, or if it is due to the inherent inability of the correction terms to approximate the pressure field, the performance of the \emph{exact} correction terms is evaluated (i.e., an ``ideal'' data-driven correction~\cite{mohebujjaman2019physically} is considered.). 
The rationale for this investigation is that removing the minimization error will provide insight into the ability of the pressure correction terms to 
improve the pressure accuracy.
In Figure \ref{exactP}, the effect of the addition of the exact pressure correction terms is displayed.
%
As can be seen in Figure \ref{exactP}, the pressure corrections have no effect on the dynamical system obtained with the supremizer enrichment. 
The plots in Figure \ref{exactP} show that the approximations $\tilde{H}\mathbf{b}$ and $\tilde{P} \mathbf{a}$ are ineffective because the exact terms $\mathbf{\tau}^{\text{exact}}_{p(1)}$ and $\mathbf{\tau}^{\text{exact}}_{p(2)}$ are also ineffective. 
We note that, if there is no improvement with the exact terms, we cannot hope to 
obtain an improvement with the approximated corrections. This clearly suggests that the structure of the corrections, rather than the accuracy of the minimization, must be improved.

We make the following remarks:
\begin{itemize}
\item The terms $\mathbf{H} \mathbf{b}$ and $\mathbf{P} \mathbf{a}$ are linear, 
whereas the velocity term $\mathbf{a}^T \mathbf{C} \mathbf{a}$ is 
nonlinear.
Thus, we expect the latter to have a more pronounced effect on the pressure and velocity accuracy. 
Similarly, we expect the correction terms for $\mathbf{a}^T \mathbf{C} \mathbf{a}$ to have a more significant effect on the velocity and pressure accuracy than the correction terms for $\mathbf{H} \mathbf{b}$ and $\mathbf{P} \mathbf{a}$.
(We note, however, that there can be cases at lower Reynolds numbers where the linear correction terms have a 
significant effect on results; for instance, in \cite{koc2019commutation} a linear correction term for the viscous term in the momentum equation is introduced, producing an improvement for 
$\nu$ values larger than $10^{-3}$.) 
\item The supremizer formulation \eqref{reduced_system} does not present a specific equation for the pressure variable. As a result, it is not possible to exploit a term exclusively based on the pressure modes and modal coefficients. 
For 
this reason, given the construction strategy of the data-driven model, there is no way to build a correction 
that is directly acting upon the pressure reduced coefficients through the pressure modes. Instead, the only pressure correction terms that can be included in the supremizer formulation either involve both pressure and velocity coefficients within the momentum equation, or 
the velocity and supremizer coefficients in the continuity equation. 
\end{itemize}
Thus, a different formulation is taken into account in the following section, in order to understand whether the presence of a specific equation for pressure offers the opportunity to devise a targeted correction term, which is able to affect the pressure field. The PPE model, which includes a pressure equation, naturally offers the opportunity to include different specific pressure corrections.
\subsection{Analysis of the PPE-ROM without corrections}
\label{PPE_standard_results}
\begin{figure}[h!]
\subfloat[Velocity error ($1^{st}$-order time scheme)]{\includegraphics[ width =0.419 \textwidth]{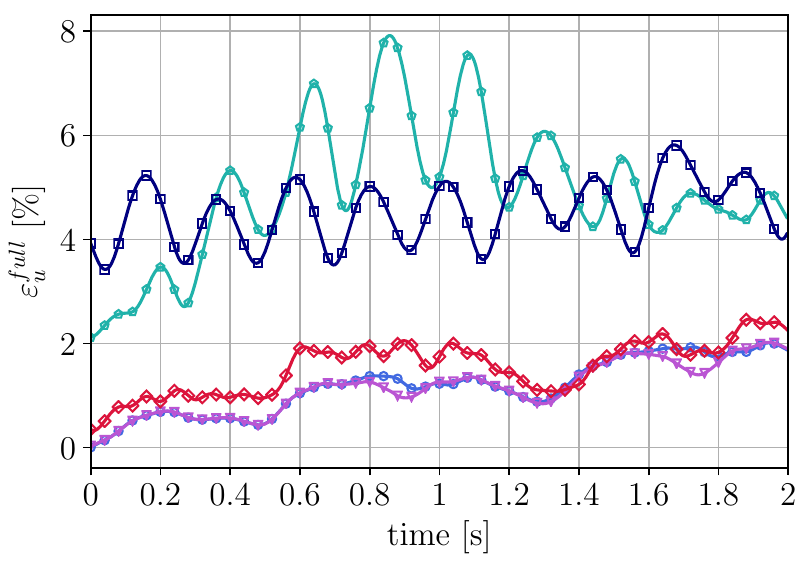} \label{PPE_vel_nocorr1ord} }
\subfloat[Pressure error ($1^{st}$-order time scheme)]{\includegraphics[ width =0.4 \textwidth, trim=0 0 6cm 0]{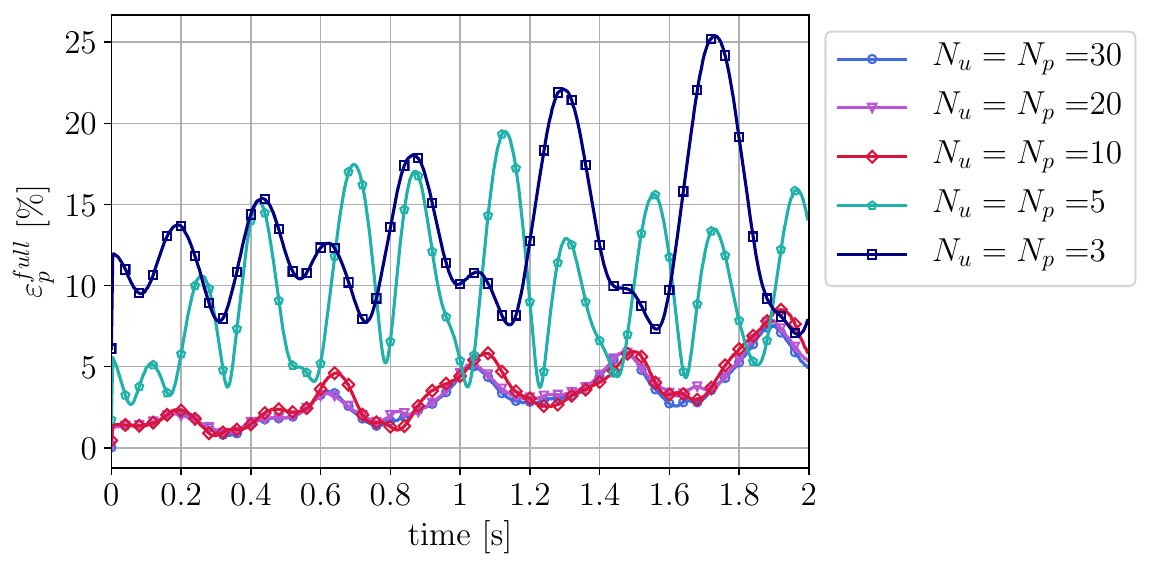} \label{PPE_pres_nocorr1ord}} \\
\subfloat[Velocity error ($2^{nd}$-order time scheme)]{\includegraphics[ width =0.425 \textwidth]{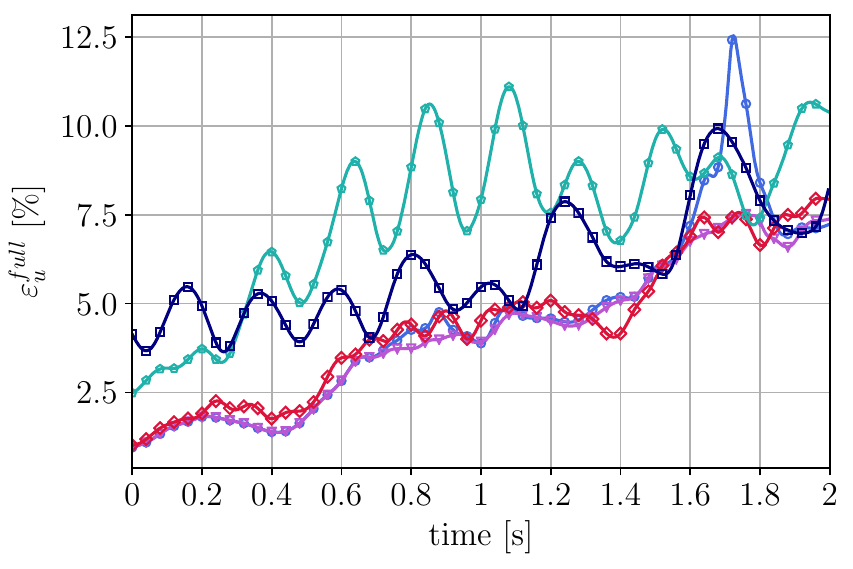} \label{PPE_vel_nocorr2ord}}
\subfloat[Pressure error ($2^{nd}$-order time scheme)]{\includegraphics[ width =0.4 \textwidth, trim=0 0 6cm 0]{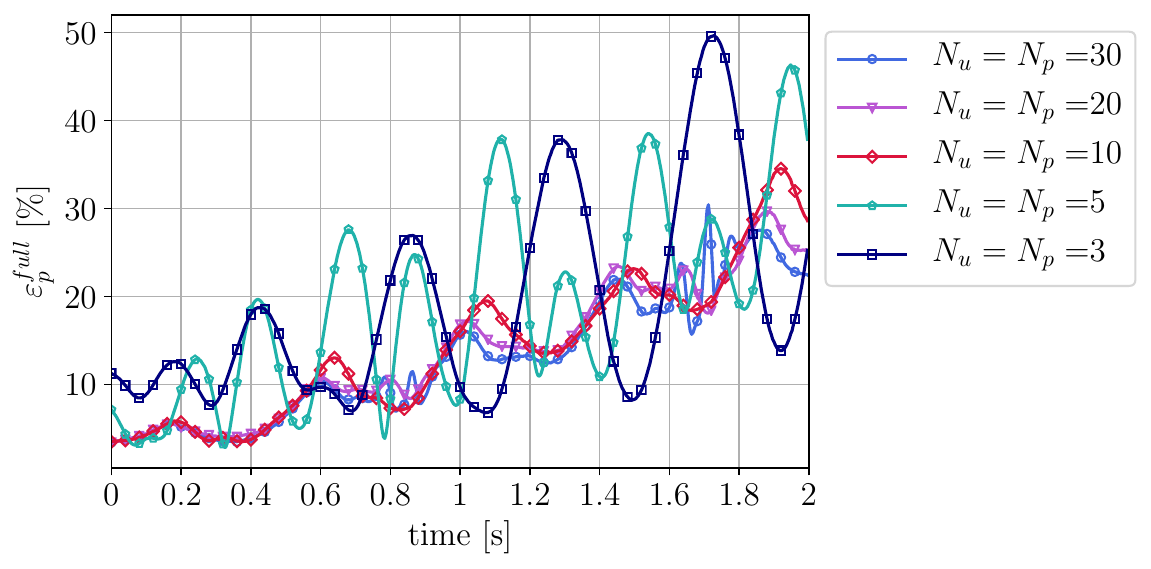} \label{PPE_pres_nocorr2ord}}
\caption{
Relative errors of the 
velocity and 
pressure, 
for different combinations of number of modes, with $N_u=N_p$.}
\label{PPE_nocorr}
\end{figure}

In this section, the solution of the dynamical system \eqref{ppegen} is 
investigated for the parameters $c_u=c_{D}=c_G=0$, i.e., the standard PPE approach without any correction terms is considered.
The relative velocity and pressure errors
are displayed in Figure \ref{PPE_nocorr}, where a first- and a second-order time scheme are used (Figure \ref{PPE_nocorr}\protect \subref{PPE_vel_nocorr1ord}, \protect \subref{PPE_pres_nocorr1ord} and \ref{PPE_nocorr}\protect \subref{PPE_vel_nocorr2ord}, \protect \subref{PPE_pres_nocorr2ord}, respectively).

From Figure \ref{PPE_nocorr}\protect \subref{PPE_vel_nocorr1ord} and \protect \subref{PPE_pres_nocorr1ord}, one can note that the accuracy increases as the number of modes increases, as expected. The stability issue exhibited by the supremizer approach 
and discussed in section \ref{sup_enrich_nocorr} is not observed in this case. However, when a second-order time scheme is used, the results in Figure \ref{PPE_nocorr}\protect \subref{PPE_vel_nocorr2ord} and \protect \subref{PPE_pres_nocorr2ord} appear more unstable, especially when a large number of pressure and velocity modes are 
used in the reduced simulations. The reason for this 
could be that already pointed out in section \ref{sup_enrich_nocorr}, 
i.e., the numerical damping of the 
time discretizations.

As was the case with the supremizer approach, 
in the investigation of the PPE approach we are also interested in the \emph{marginally-resolved} regime. Thus, most of the numerical tests 
in this section are carried out using $N_u=N_p=3$ or $N_u=N_p=5$.
The aim of the next sections will be 
to evaluate different strategies for the improvement of the pressure 
accuracy in the PPE framework. This will involve comparing 
three different approaches:
\begin{itemize}
    \item[(i)] the standard PPE method (i.e., without any correction terms);
    \item[(ii)] the PPE correction terms proposed in this section and presented in section \ref{datappe};
 \item[(iii)] the velocity correction terms developed for the supremizer approach (section \ref{velresults}).
\end{itemize}

\subsection{Effect of corrections in the PPE-ROM approach}
\label{PPE1}
The present section discusses the results obtained with the PPE approach enhanced with different data-driven corrections. The following studies are carried out:
\begin{itemize}
    \item a study of the separate and combined effects of different pressure corrections in the pressure Poisson 
    equation; 
    \item an evaluation of the combined influence of velocity 
    and pressure corrections, comparing the 
    approaches presented in section \ref{proposals}.
\end{itemize}
\subsubsection{Effect of pressure corrections in the PPE-ROM}
\label{pres_comb_sec}
\begin{figure}[h!]
\centering
\subfloat[Relative velocity error 
]{\includegraphics[  scale=0.45]{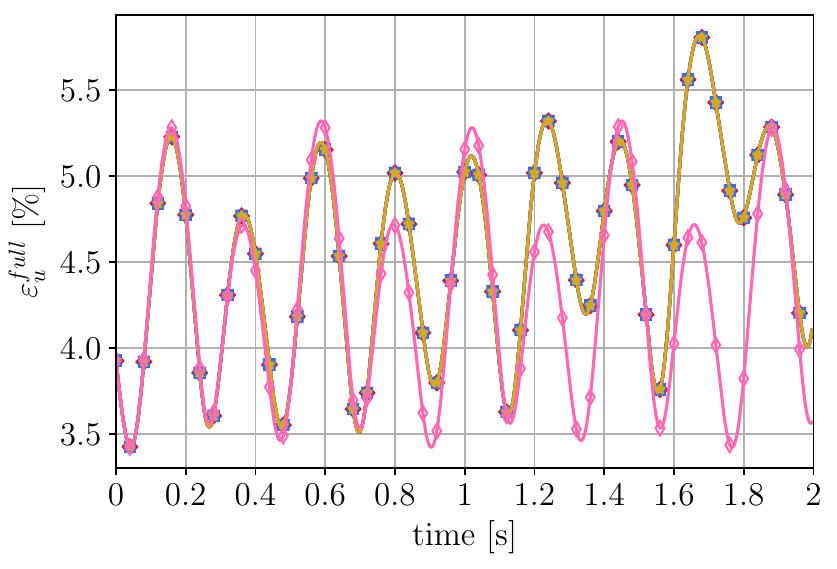}}
\subfloat[Relative pressure error 
]{\includegraphics[  scale=0.45]{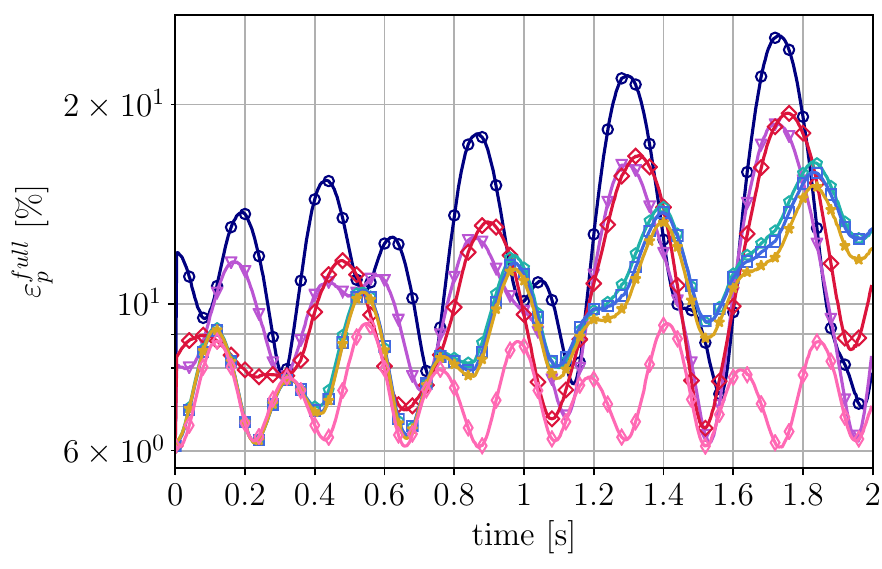}}
\caption{
Relative errors of the 
velocity (a) and pressure (b), considering $N_u=N_p=3$. Results in the following cases are displayed: without corrections (\protect\input{markers-plots/mark0}); with the correction $\boldsymbol{\tau}_D$ (\protect\input{markers-plots/mark1}); with the correction $\boldsymbol{\tau}_G$ (\protect\input{markers-plots/mark2}); and with both $\boldsymbol{\tau}_D$ and $\boldsymbol{\tau}_G$ found from two disjoint least squares problems (\protect\input{markers-plots/mark3}), in Cases 1 and 2 of section \ref{proposals} (\protect\input{markers-plots/mark4} and \protect\input{markers-plots/mark6}, respectively). Results are compared with the reconstruction errors, 
computed from the projected fields (\protect\input{markers-plots/mark5}).}

\label{errDG}
\end{figure}

In this subsection, the pressure correction terms for $\boldsymbol{\tau}_D$ and $\boldsymbol{\tau}_G$, presented in sections \ref{Dcorr_ref} and \ref{Gterm}, are introduced. We recall that 
these corrections, used in the pressure Poisson 
equation, are based on the reduced pressure and velocity vectors, respectively. 
Their separate 
and combined effects are 
evaluated in this section. In particular, the dynamical system \eqref{ppegen} is solved considering $c_u=0$. Figure \ref{errDG} displays the results for the following cases:
\begin{itemize}
    \item $c_D=0$, $c_G=1$;
    \item $c_D=1$, $c_G=0$, where the quadratic ansatz for $\boldsymbol{\tau}$ is computed;
    \item $c_D=c_G=1$, where the corrections are obtained by solving two different optimization problems for the pressure corrections presented in sections \ref{Dcorr_ref} and \ref{Gterm};
    \item $c_D=c_G=1$, where a unique optimization problem is solved, using an ansatz with a linear dependence on $\mathbf{a}$ and a quadratic dependence on $\mathbf{b}$ (section \ref{proposals} Case 1);
    \item $c_D=c_G=1$, in which a unique optimization problem is solved, considering an ansatz with a quadratic dependence on the vector $\mathbf{ab}$ (section \ref{proposals} Case 2).
\end{itemize}
The plots in Figure \ref{errDG} show that when only 
(one or more) pressure corrections are added to the system, the accuracy of the pressure field approximation improves, while the velocity approximation accuracy 
remains unchanged.
In addition, the plots suggest that when the system \eqref{ppegen} is solved including both pressure corrections at the same time, there is a significant improvement in the pressure accuracy. 
Finally, Figure \ref{errDG}(b) shows that the results achieved using the method presented in Case 2 (section \ref{proposals}) are slightly better than the results achieved following the method in Case 1.

An important observation is that the pressure corrections tested in the PPE-ROM approach significantly 
improve 
the pressure results, whereas in the SUP-ROM approach the results are unchanged when pressure corrections are added to the reduced order system (section \ref{presSUP_2}).
The conclusion of the present part is that both pressure corrections added in the PPE-ROM produce a significant improvement 
of the pressure accuracy. 
\subsubsection{Combined effect of velocity and pressure corrections in the PPE-ROM}
\label{vel_pres_combined}
\begin{figure}[h!]
\centering
\subfloat[Relative velocity error 
]{\includegraphics[ scale=0.45]{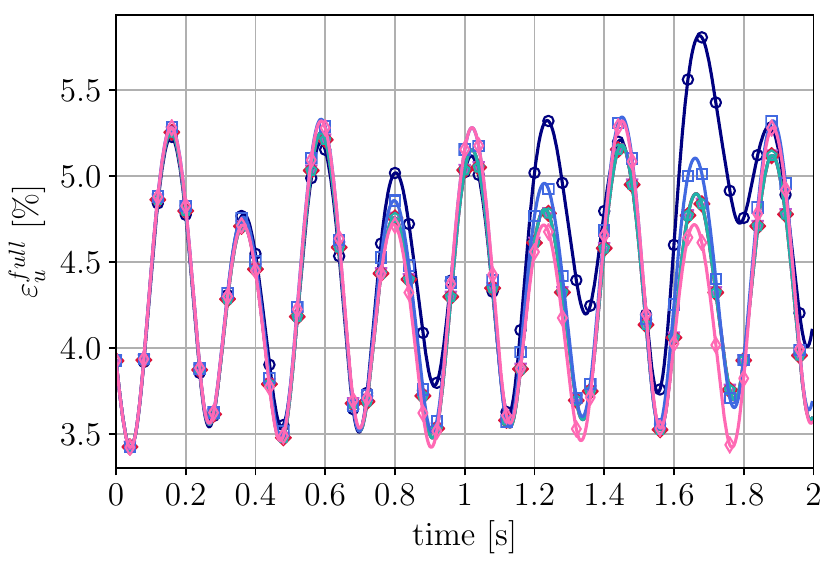}}
\subfloat[Relative pressure error 
]{\includegraphics[ scale=0.45 ]{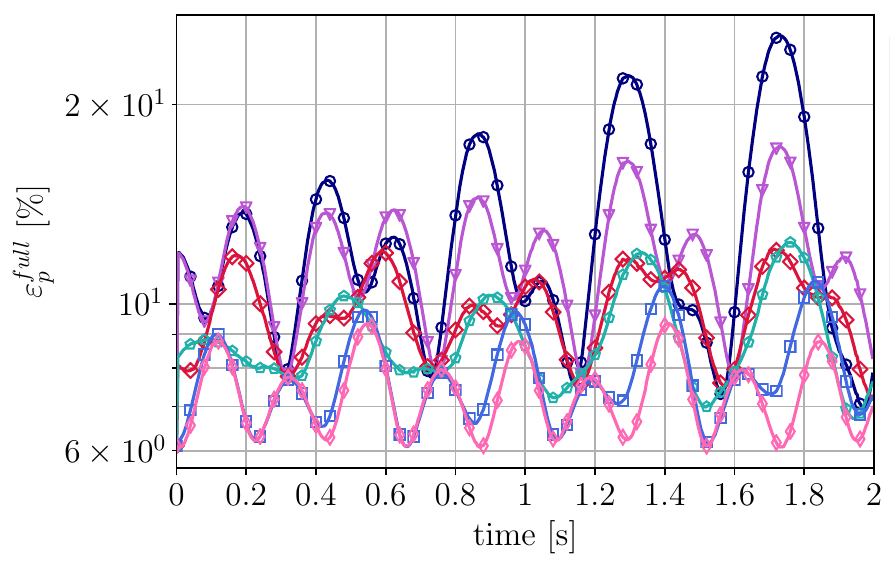}}
\caption{Relative errors 
of the velocity (a) and 
pressure (b), 
considering $N_u=N_p=3$. Results in the following cases are displayed: without corrections (\protect\input{markers-plots/mark0}); with the velocity correction $\boldsymbol{\tau}_u$ (\protect\input{markers-plots/mark1}); with the corrections $\boldsymbol{\tau}_u$ and $\boldsymbol{\tau}_D$ (\protect\input{markers-plots/mark2}); and with $\boldsymbol{\tau}_u$ and $\boldsymbol{\tau}_G$ (\protect\input{markers-plots/mark3}), with $\boldsymbol{\tau}_u$, $\boldsymbol{\tau}_D$, and $\boldsymbol{\tau}_G$ in Case 3 of section \ref{proposals} (\protect\input{markers-plots/mark4}). Results are compared with the reconstruction errors, 
computed from the projected fields (\protect\input{markers-plots/mark5}).}
\label{errCDG}
\end{figure}

The results in the previous subsection 
confirmed that adding data-driven corrections in the pressure Poisson equation 
leads to improvements in the ROM pressure 
accuracy. 
Next, we try to understand
if combining the corrections in the Poisson equation and momentum equation 
yields additional gains. 
Thus, in this section, the velocity correction term $\boldsymbol{\tau}_u$ is introduced in system \eqref{ppegen}.
Figure \ref{errCDG} displays results for the following cases:
\begin{itemize}
    \item $c_u=c_D=1$, $c_G=0$;
    \item $c_u=c_G=1$, $c_D=0$;
    \item $c_u=c_G=c_D=1$, where the method presented in \eqref{proposals} (Case 3) is 
    used and a unique least squares problem is solved to find all correction terms.
\end{itemize}
Figure \ref{errCDG} displays the results obtained when the velocity correction is combined with one or both pressure correction terms. As expected, the left plot confirms that the presence of data-driven velocity correction terms 
improves the ROM velocity approximation. In addition, the right diagram suggests that, as was the case for the supremizer approach, adding a velocity correction term in the momentum equation leads to more accurate pressure results. 
These accuracy gains 
are consistent with the 
gains obtained by 
using the pressure correction terms in the pressure Poisson equation. 
As a result, errors obtained when 
both the velocity and the pressure corrections are active 
are also very close to the reconstruction error for the velocity field --- which is the best result we can obtain --- as can be seen in the left diagram in Figure \ref{errCDG}.
\begin{figure}[h!]
\centering
\subfloat[Relative velocity error 
]{\includegraphics[ scale=0.45]{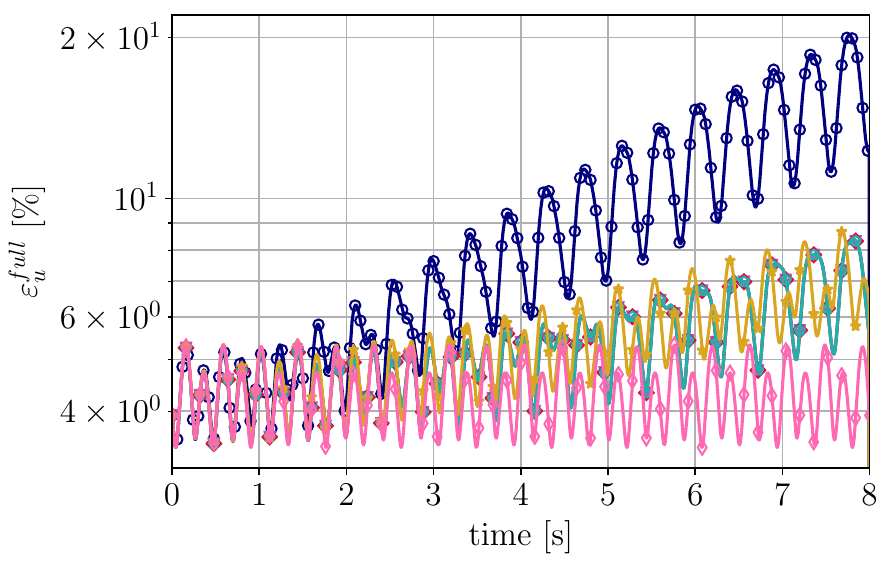}}
\subfloat[Relative pressure error 
]{\includegraphics[ scale=0.45 ]{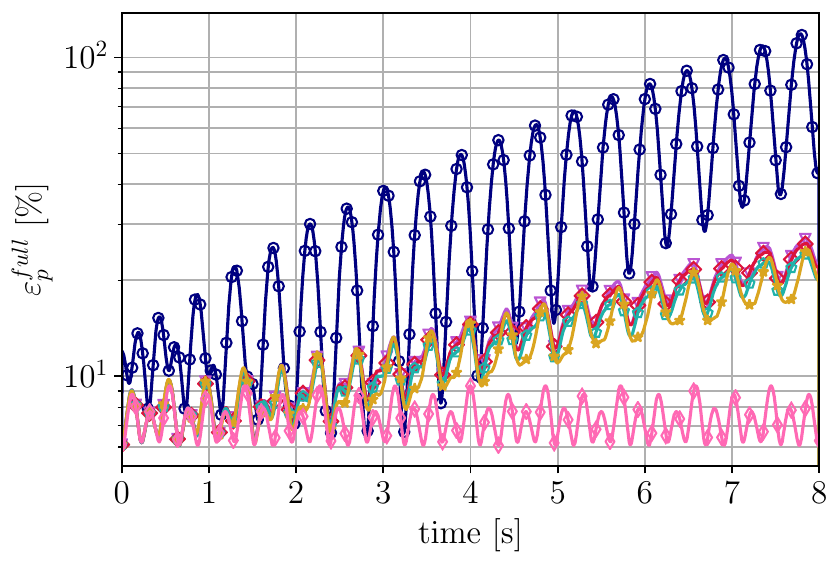}}
\caption{Relative errors of the 
velocity (a) and 
pressure (b), 
considering $N_u=N_p=3$. Results in the following cases are displayed: without corrections (\protect\input{markers-plots/mark0}); with $\boldsymbol{\tau}_u$, $\boldsymbol{\tau}_D$, and $\boldsymbol{\tau}_G$ disjoint (\protect\input{markers-plots/mark1}); and with the velocity and pressure corrections in Cases 1, 2, and 3 of section \ref{proposals} (\protect\input{markers-plots/mark2}, \protect\input{markers-plots/mark3}, and \protect\input{markers-plots/mark6}, respectively). Results are compared with the reconstruction errors (\protect\input{markers-plots/mark5}).}
\label{errCDGtot_long}
\end{figure}
In Figure \ref{errCDGtot_long},  the results obtained combining the momentum equation correction with all the Poisson corrections forms developed are 
shown.  In particular, the corrections proposed in section \ref{proposals} are compared; those data-driven terms are built from the first 2 seconds of the online simulation and in Figure \ref{errCDGtot_long} the results for 8 seconds of simulation are displayed. All the methods produce a pressure reduced solution that is very close to the projected pressure 
solution and the results of all methods look similar. 
However, the most accurate method is the one presented in section \ref{proposals} Case 3.

The following conclusions can be drawn from the analysis of the correction terms 
in the PPE-ROM approach:
\begin{itemize}
    \item the velocity correction reduces the error for both the velocity and the pressure fields, whereas the pressure corrections added in the Poisson equation only improve the pressure field;
    \item the most significant improvement in the accuracy of the reduced pressure field is reached when all data-driven corrections are added to the reduced system.
\end{itemize}
\subsection{Flow field qualitative inspection}
\label{graph_sec}
The effect of the novel correction terms in the reduced formulations is also examined through the observation of contour plots of the velocity and pressure fields obtained for both SUP-ROM and PPE-ROM approaches. 

\begin{figure}[h!]
\centering
    \subfloat[]{\includegraphics[width=0.46\textwidth]{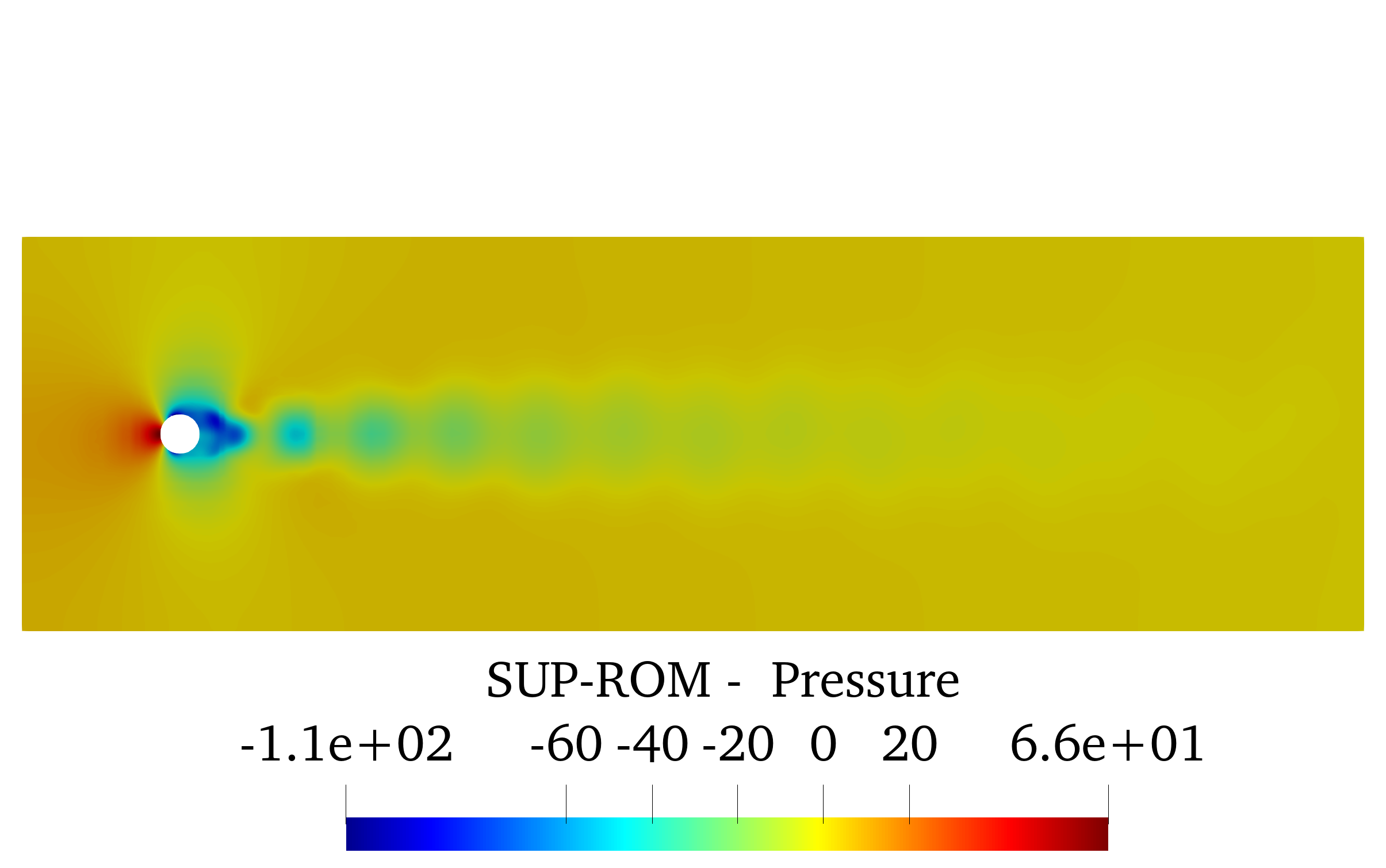}}
    \subfloat[]{\includegraphics[width=0.46\textwidth]{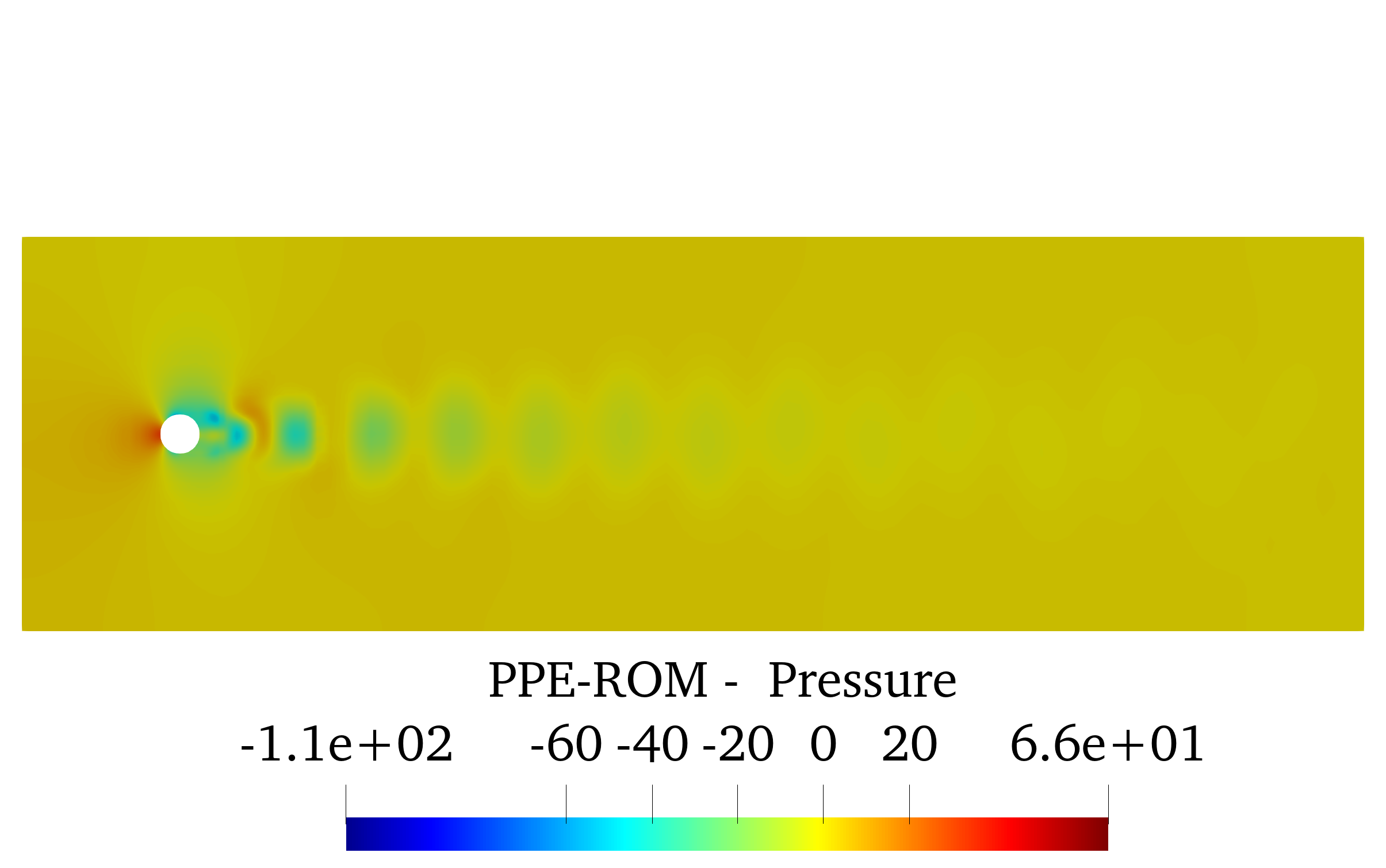}}
    \\ \vspace{-1.6cm}
    \subfloat[]{\includegraphics[width=0.46\textwidth]{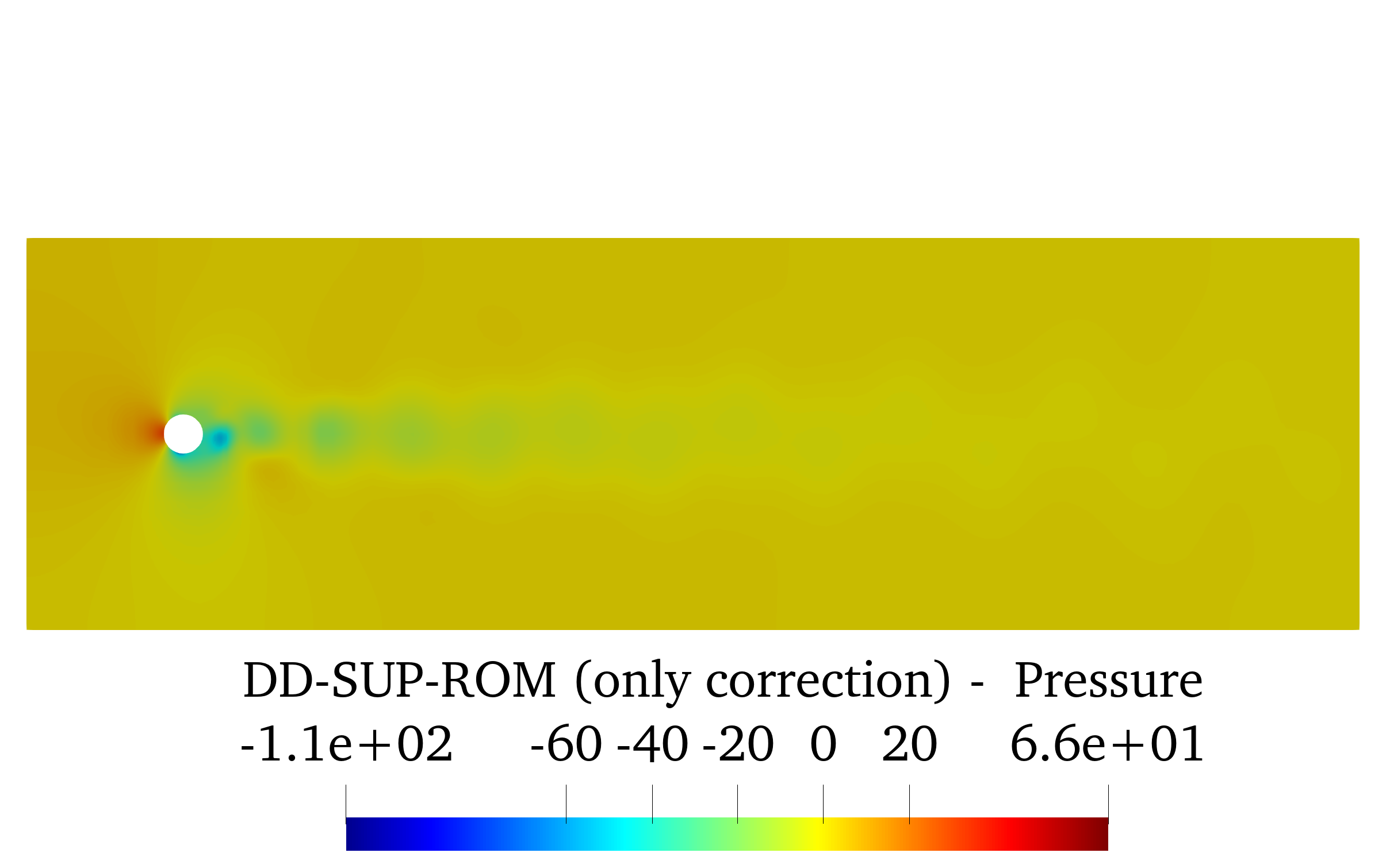}}
    \subfloat[]{\includegraphics[width=0.46\textwidth]{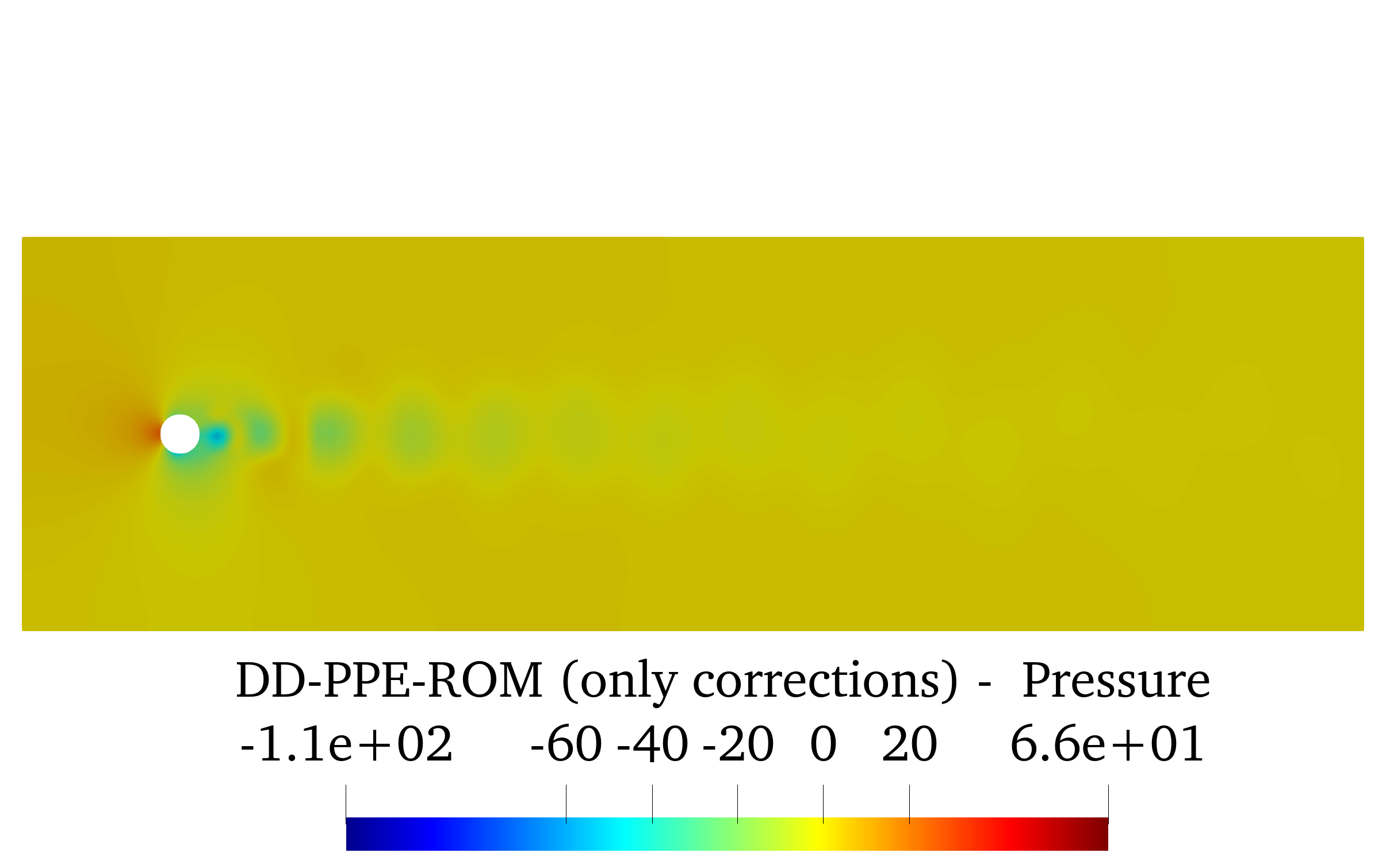}}\\ \vspace{-1.6cm}
    \subfloat[]{\includegraphics[width=0.46\textwidth]{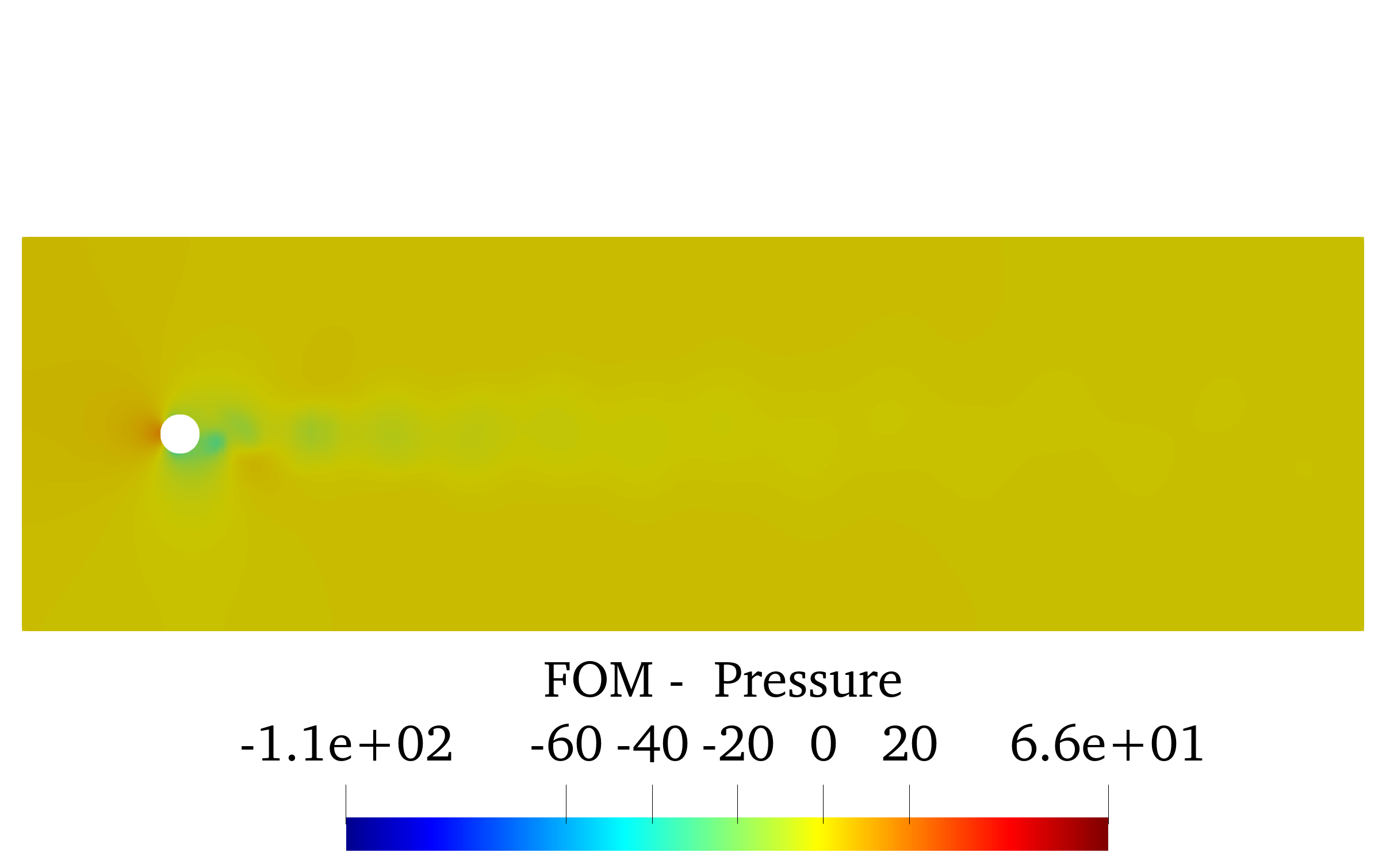}}
    \vspace{-0.4cm}
    \caption{Representation of the pressure field for the FOM, 
    SUP-ROM, and 
    PPE-ROM simulations with and without the correction terms.}
    \label{paraview1}
\end{figure}

\begin{figure}[h!]
\centering
    \subfloat[]{\includegraphics[width=0.46\textwidth]{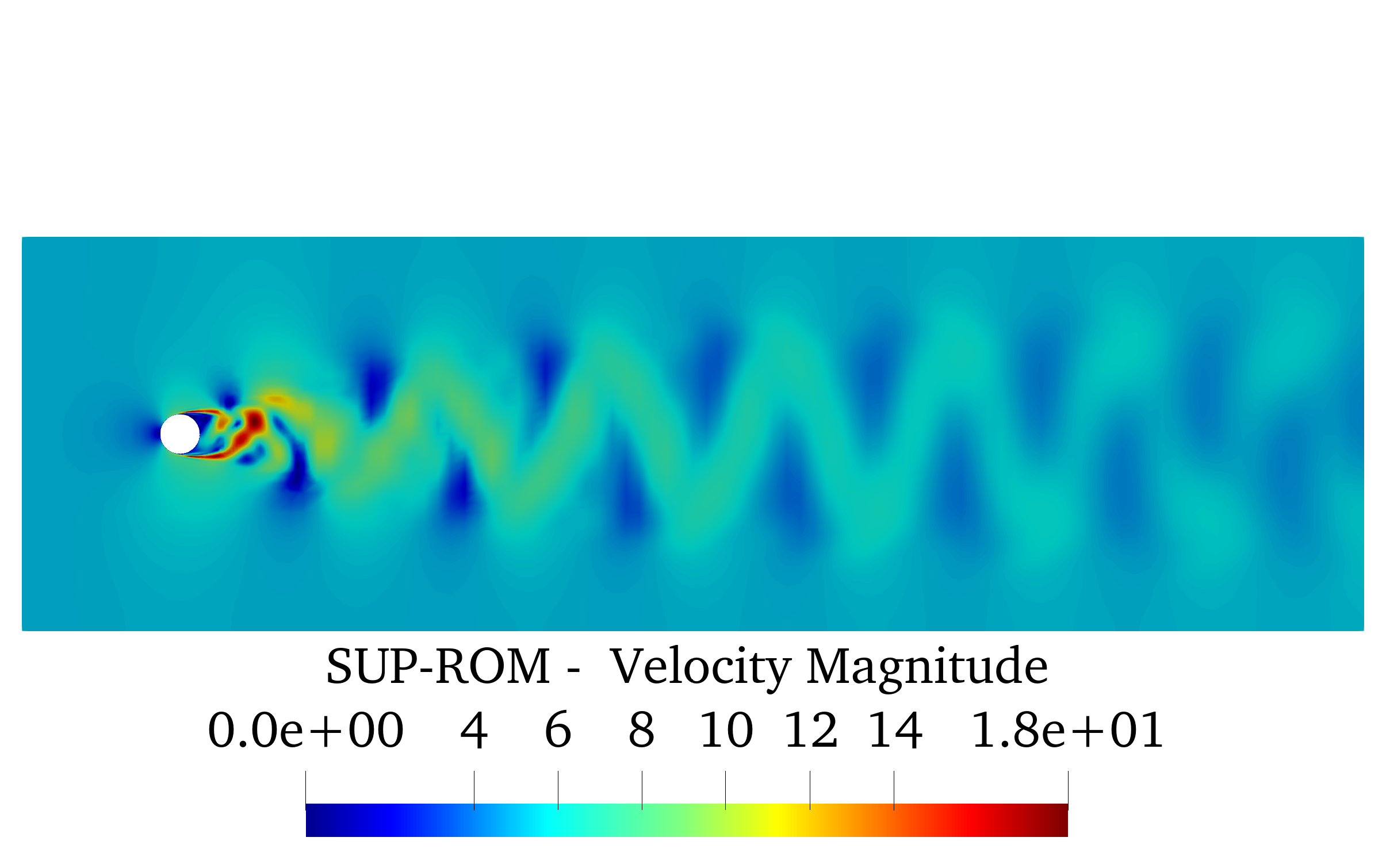}}
    \subfloat[]{\includegraphics[width=0.46\textwidth]{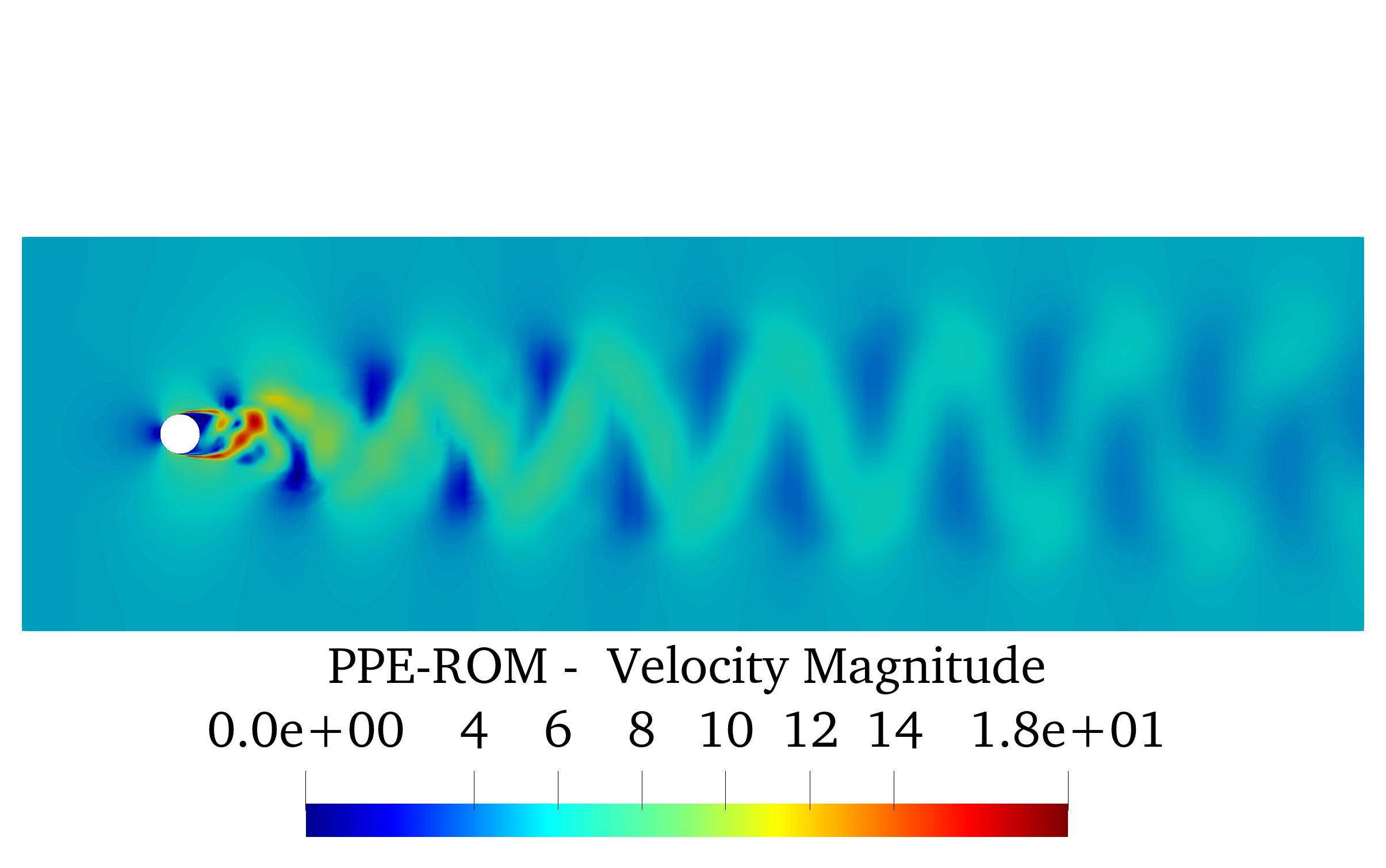}}
    \\ \vspace{-1.6cm}
    \subfloat[]{\includegraphics[width=0.46\textwidth]{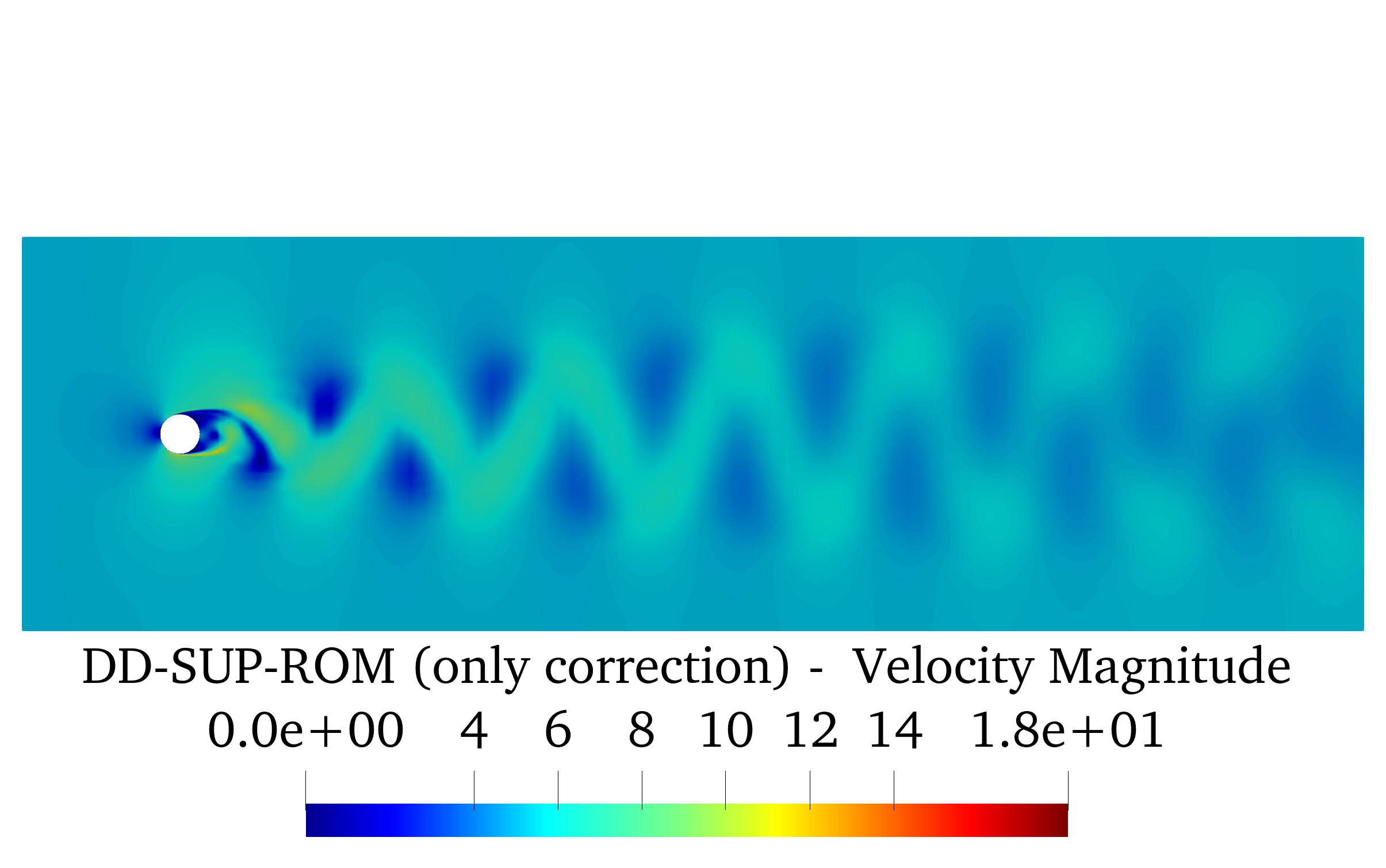}}
    \subfloat[]{\includegraphics[width=0.46\textwidth]{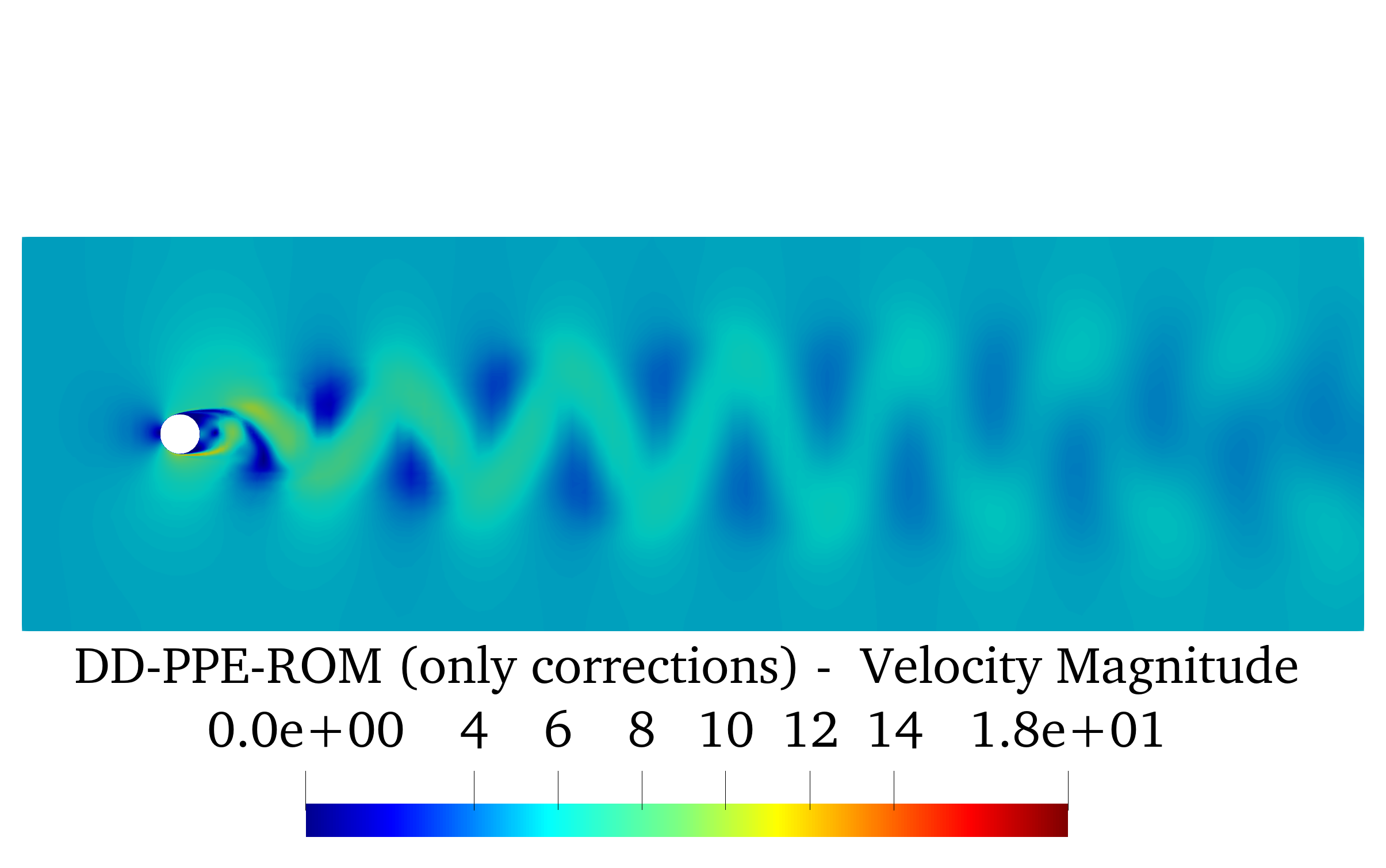}}\\ \vspace{-1.6cm}
    \subfloat[]{\includegraphics[width=0.46\textwidth]{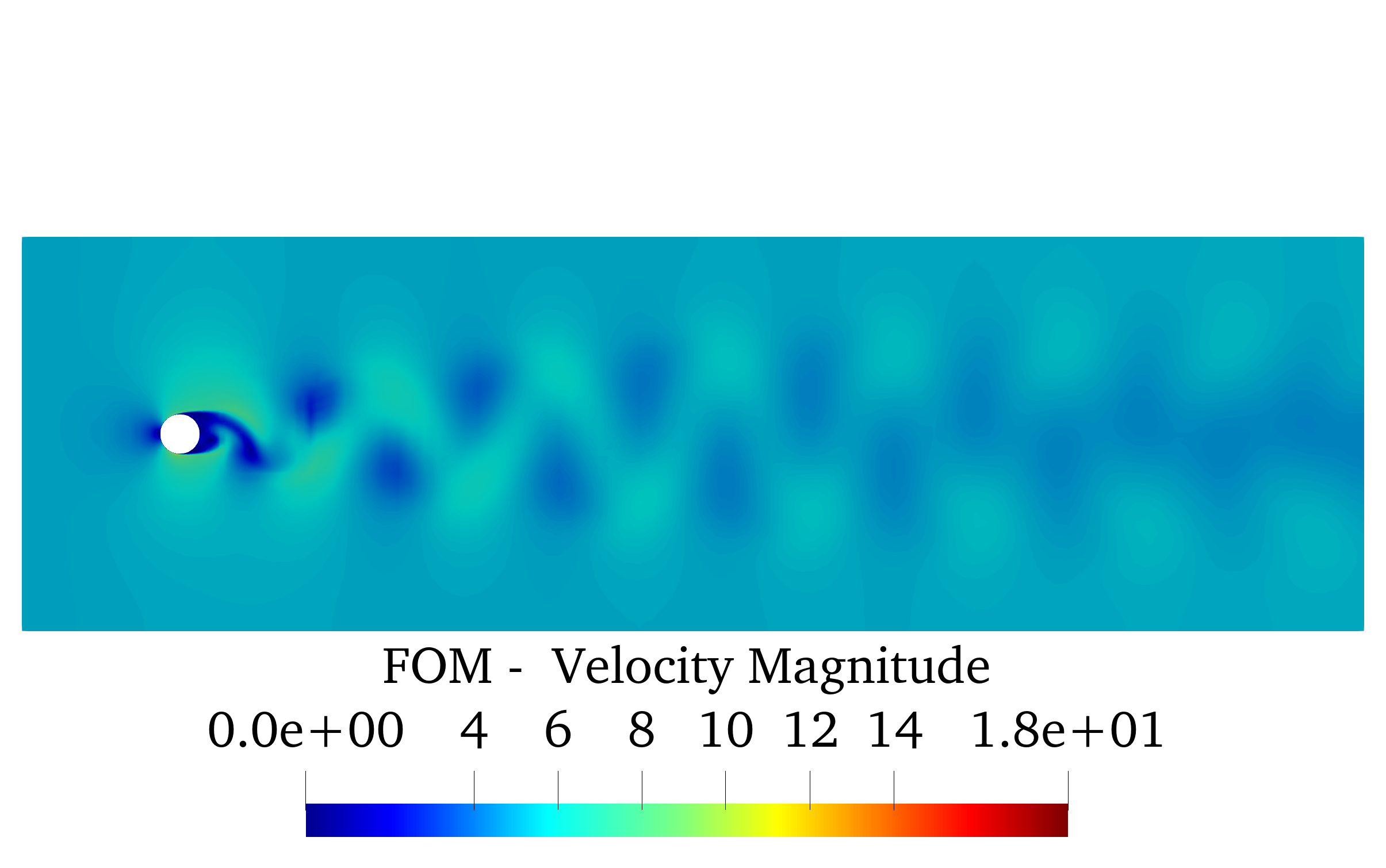}}
    \vspace{-0.4cm}
    \caption{Representation of the velocity magnitude field for the FOM, 
    SUP-ROM, and 
    PPE-ROM simulations with and without the correction terms.}
    \label{paraview2}
\end{figure}

The pressure and the velocity magnitude fields are 
displayed in Figures \ref{paraview1} and \ref{paraview2}, respectively, for different SUP-ROM and PPE-ROM simulations.
The corrections 
used are the constrained velocity correction examined in section \ref{constrainedvel} for the SUP-ROM, and the joint velocity and pressure correction presented in section \ref{proposals} (Case 3). 
The second-order time integration scheme is used to generate the plots.

The POD is performed on the time interval $[79.992,99.992]$ seconds and the reduced order systems \eqref{supgen} and \eqref{ppegen} are solved 
on the time interval $[79.992,87.992]$ seconds, since the maximum length of the online simulations carried out is $8$ seconds. 
All the plots display results
at the final time step of online simulations, which is second $87.992$.

In order to obtain the reduced fields, the reduced order systems with $N_u=N_p=5$, and $N_{sup}=5$ for the supremizer approach, are solved. The reduced fields are computed from the 
vectors of coefficients $\mathbf{a}$ and $\mathbf{b}$, and the POD modes $(\boldsymbol{\phi}_i)_{i=1}^{N_u+N_{sup}}$ and $(\chi_i)_{i=1}^{N_p}$, as in \eqref{appfields}. 

The contour plots confirm that 
there is a significant difference between the fields computed with the standard ROMs and those including the data-driven terms. 
Specifically, the fields in Figures \ref{paraview1} and \ref{paraview2}(c) and (d) are closer to the full order fields, especially in the region around the cylinder. 
We emphasize that the improvement of the error near the circular cylinder is an important gain as it 
leads to a better reconstruction of the ROM lift coefficient.

\newpage

\RA{
\subsection{Tests with different time windows}
\label{more_tests}
In subsections \ref{sup_enrich_nocorr}, \ref{velresults}, \ref{PPE1}, and \ref{graph_sec}, we considered 
the following time windows: 
$T_{\text{offline}}=\SI{20}{\second}$ and $T_{\text{correction}}=\SI{2}{\second}$ 
to construct the POD basis and the correction terms, respectively;
$\SI{2}{\second}$
to test the ROMs in the reconstructive regime; and
$\SI{8}{\second}$ 
to test the ROMs in the predictive regime.
However, a single test considering only one time window to build the POD basis does not guarantee the 
effectiveness of data-driven corrections.
Therefore, in this section we test our method 
on different datasets. 

We note that, for periodic and quasi-periodic flows, 
the time windows used for collecting snapshots ($T_{\text{offline}}$) and for building the correction terms ($T_{\text{correction}}$) are generally chosen according to the Strouhal number, i.e., the time period of vortex shedding.
In Table~\ref{table2} we summarize the tests presented in this subsection, where the time windows are related to the flow period $P\sim \SI{0.89}{\second}$.
In our numerical investigation, we apply the unconstrained and constrained correction terms for the SUP-ROM formulation, and both the velocity and the pressure correction in the PPE-ROM approach, following the ansatz in Case 3 in subsection \ref{proposals}.
We also emphasize that tests (b) and (d) in Table~\ref{table2} investigate the predictive accuracy of the numerical methods. 
Specifically, in tests (b) and (d) we test the new methods on time intervals that are longer than the time intervals used to construct the POD basis and correction terms.

From Figures \ref{2P_window}, \ref{2P_window_predictive}, \ref{P_window}, \ref{P_window_predictive}, the following conclusions can be drawn:
\begin{itemize}
    \item The data-driven techniques proposed in the paper lead to improvements in the results,  especially for the DD-PPE-ROM formulation and for the constrained DD-SUP-ROM, 
    when considering different time intervals for the 
    construction of the POD basis, construction of the correction terms, and ROM testing.
    \item The constrained correction in the supremizer approach outperforms the unconstrained one (Figures \ref{2P_window}\subref{sup_2P} and \ref{2P_window_predictive}\subref{sup_2P_pred}). In particular, the unconstrained correction 
    is less accurate than
    the standard SUP-ROM, especially outside 
    the time window used to compute the correction terms (i.e., in the time interval $[P,2P]$; see Figure \ref{2P_window}\subref{sup_2P}) and 
    outside the time window used to 
    build the POD basis (i.e., in the time interval $[2P,8P]$; see Figure \ref{2P_window_predictive}\subref{sup_2P_pred}).
    \item From Figure \ref{P_window}, we can conclude that a time window 
    $\bigl[0,\frac{P}{2}\bigr]$ does not provide enough information to build the data-driven corrections, since the predictive accuracy, especially for the pressure field, 
    decreases in the interval 
   $\bigl[\frac{P}{2},P\bigr]$.
    These results suggest
    that a time interval of at least one period is necessary to build the data-driven corrections in periodic and quasi-periodic flows.
    \item The predictive analysis with respect to the POD time interval (Figures \ref{2P_window_predictive} and \ref{P_window_predictive}) shows that the data-driven formulation leads to improvements with respect to the standard approaches both in the SUP-ROM (when considering the constrained correction) and in the PPE-ROM. 
\end{itemize}

\remark{
The proposed strategy significantly increases the standard approach accuracy in the predictive regime.
We emphasize, however, that the test case used in our numerical investigation is a quasi-periodic flow. 
Thus, higher Reynolds number test cases need to be carefully analyzed.
}

\begin{table}[h!]
    \centering
    \begin{tabular}{
    p{0.06\textwidth}
    p{0.2\textwidth}
    p{0.12\textwidth}
    p{0.12\textwidth}
    p{0.12\textwidth}}
    \toprule[0.3ex]
         \textbf{Tests}&\textbf{Corresponding Figure}& $T_{\text{offline}}$ & $T_{\text{correction}}$ & $T_{\text{online}}$ \\
         \midrule
         \textbf{(a)} &\ref{2P_window}& $\SI{1.78}{\second} \sim 2P$&  $\SI{0.89}{\second} \sim P$& $\SI{1.78}{\second} \sim 2P$ \\
         \midrule
       \textbf{(b)} &\ref{2P_window_predictive} & $\SI{1.78}{\second} \sim 2P$&  $\SI{0.89}{\second} \sim P$& $\SI{7.12}{\second} \sim 8P$ \\
       \midrule
       \textbf{(c)}&\ref{P_window} & $\SI{0.89}{\second} \sim P$&  $\SI{0.45}{\second} \sim \frac{P}{2}$& $\SI{0.89}{\second} \sim P$ \\
       \midrule
      \textbf{(d)} &\ref{P_window_predictive}& $\SI{0.89}{\second} \sim P$&  $\SI{0.89}{\second} \sim P$& $\SI{3.56}{\second} \sim 4P$ \\
         \bottomrule[0.3ex]
    \end{tabular}
    \caption{Time windows considered in the test cases of section \ref{more_tests}.}
    \label{table2}
\end{table}

\newpage

\begin{figure}[h!]
\centering
\subfloat[Relative errors for SUP-ROM]{\includegraphics[ width =0.85 \textwidth, trim={0 0 0 1.5cm}, clip]{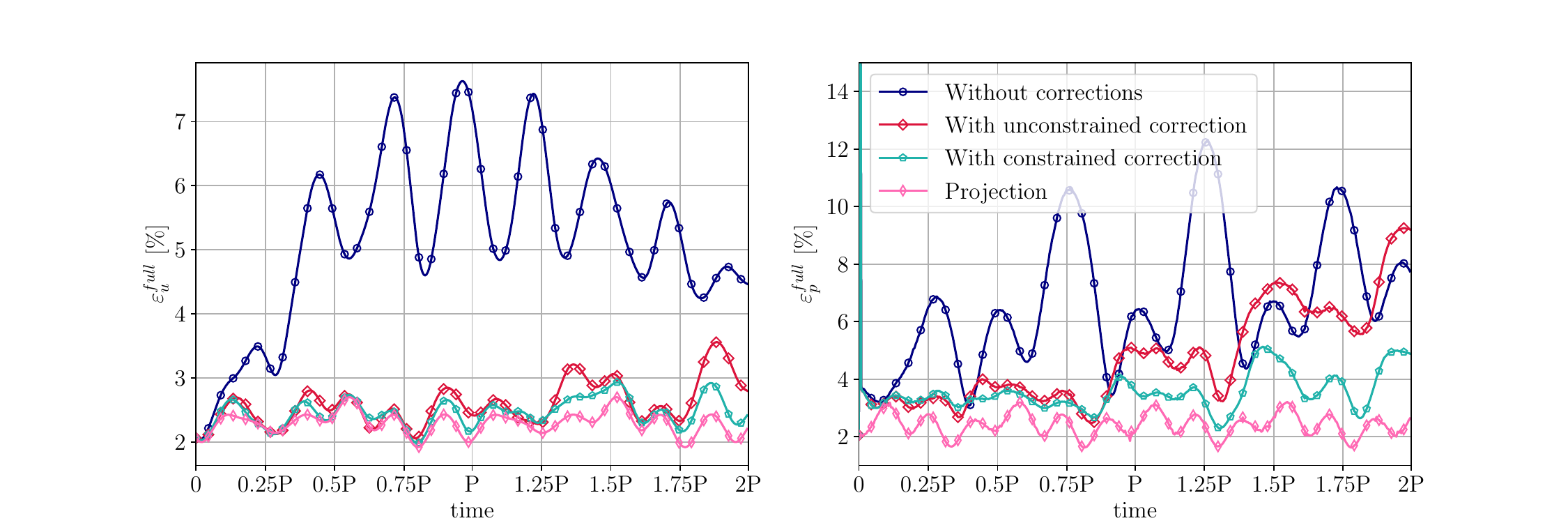} \label{sup_2P} }\\
\vspace{-0.3cm}
\subfloat[Relative errors for PPE-ROM]{\includegraphics[ width =0.85 \textwidth, trim={0 0 0 1.5cm}, clip]{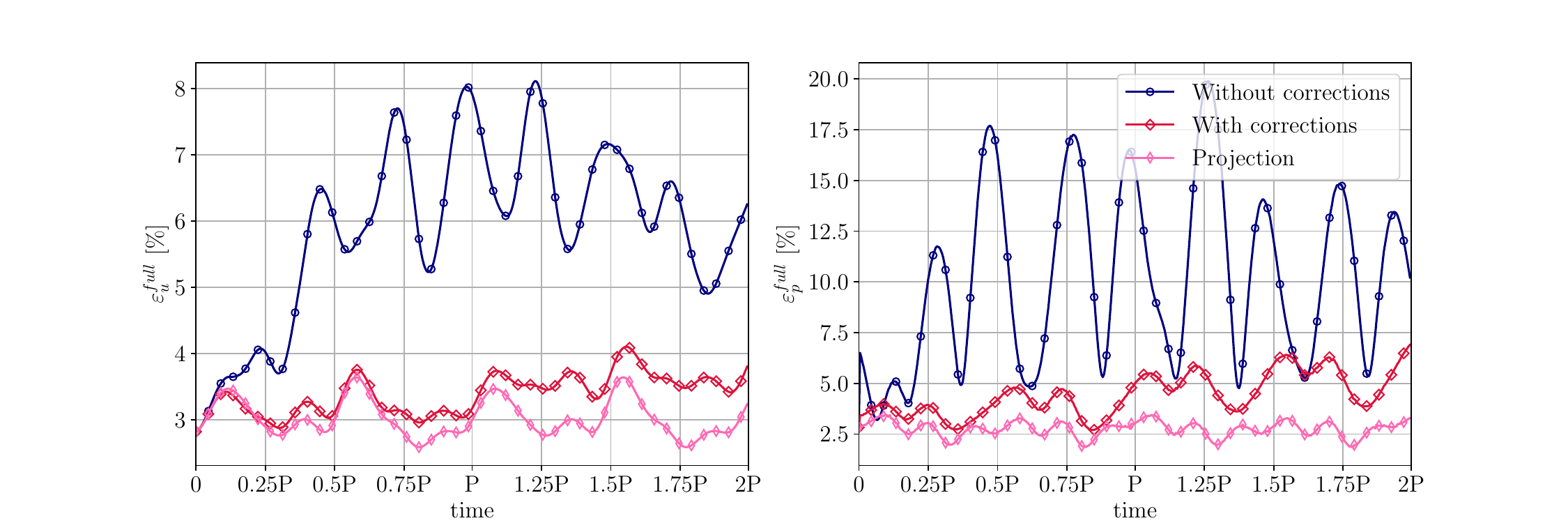} \label{ppe_2P} }
\caption{Relative errors of velocity and pressure for SUP-ROM ($N_u=N_{sup}=N_p=5$ modes) and PPE-ROM ($N_u=N_p=5$), considering $T_{offline}=2P$.}
\label{2P_window}
\end{figure}

\begin{figure}[h!]
\centering
\subfloat[Relative errors for SUP-ROM]{\includegraphics[ width =0.85 \textwidth, trim={0 0 0 1.5cm}, clip]{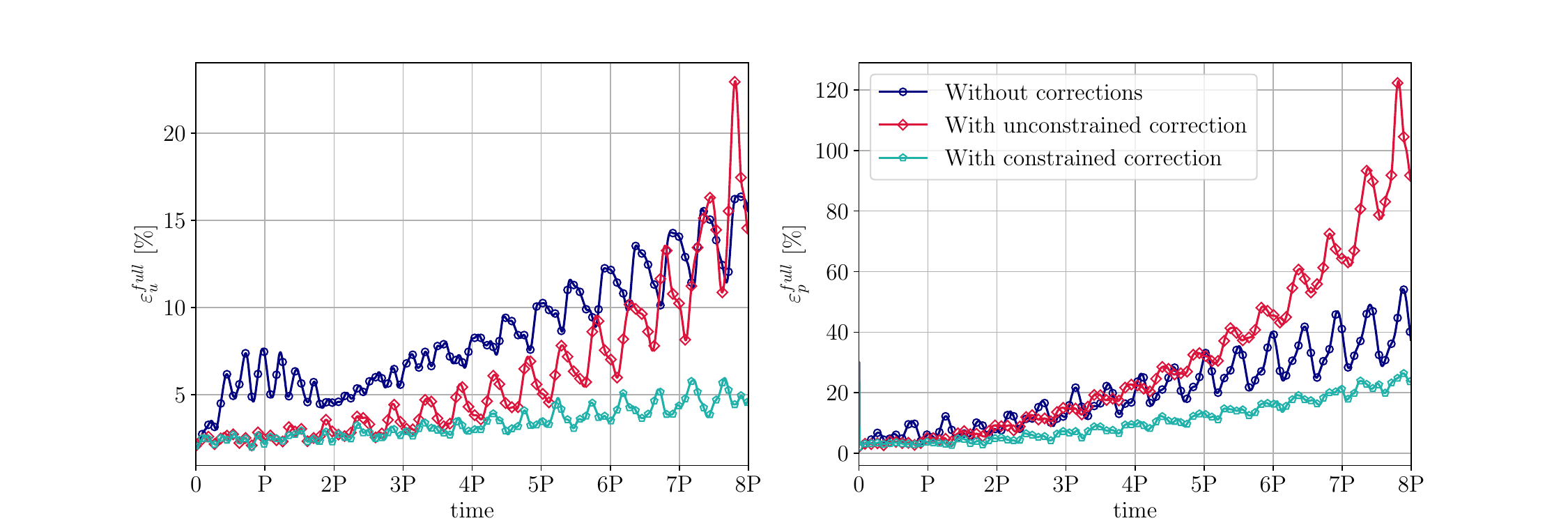} \label{sup_2P_pred} }\\
\vspace{-0.3cm}
\subfloat[Relative errors for PPE-ROM]{\includegraphics[ width =0.85 \textwidth, trim={0 0 0 1.5cm}, clip]{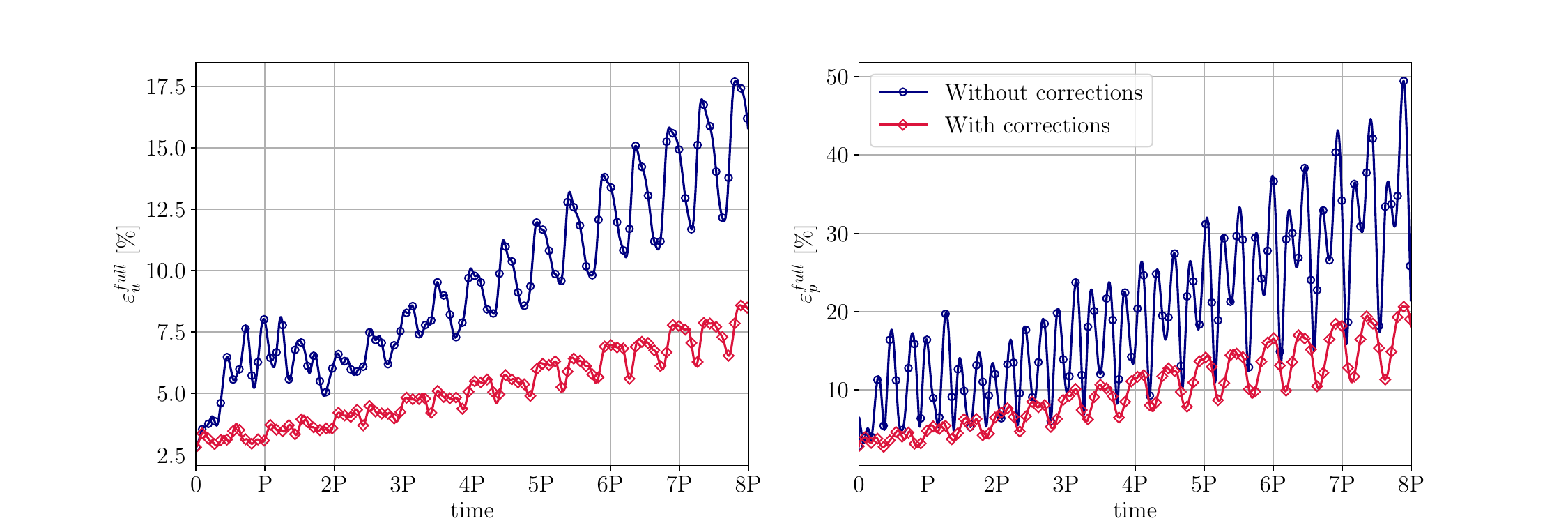} \label{ppe_2P_pred} }
\caption{
Relative errors of velocity and pressure for SUP-ROM ($N_u=N_{sup}=N_p=5$ modes) and PPE-ROM ($N_u=N_p=5$) considering a time window larger than the offline time interval, i.e.,  $T_{\text{online}}=8P$ and $T_{\text{offline}}=2P$.}
\label{2P_window_predictive}
\end{figure}

\begin{figure}[h!]
\centering
\subfloat[Relative errors for SUP-ROM]{\includegraphics[ width =0.85 \textwidth, trim={0 0 0 1.5cm}, clip]{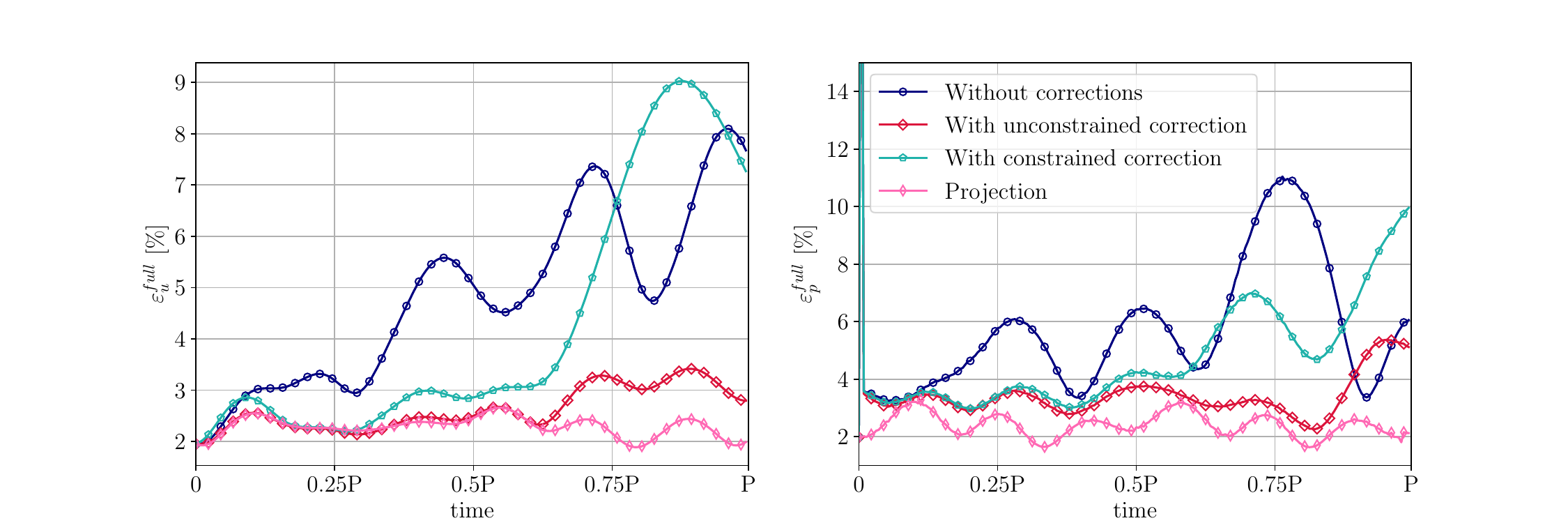} \label{sup_P} }\\
\vspace{-0.3cm}
\subfloat[Relative errors for PPE-ROM]{\includegraphics[ width =0.85 \textwidth, trim={0 0 0 1.5cm}, clip]{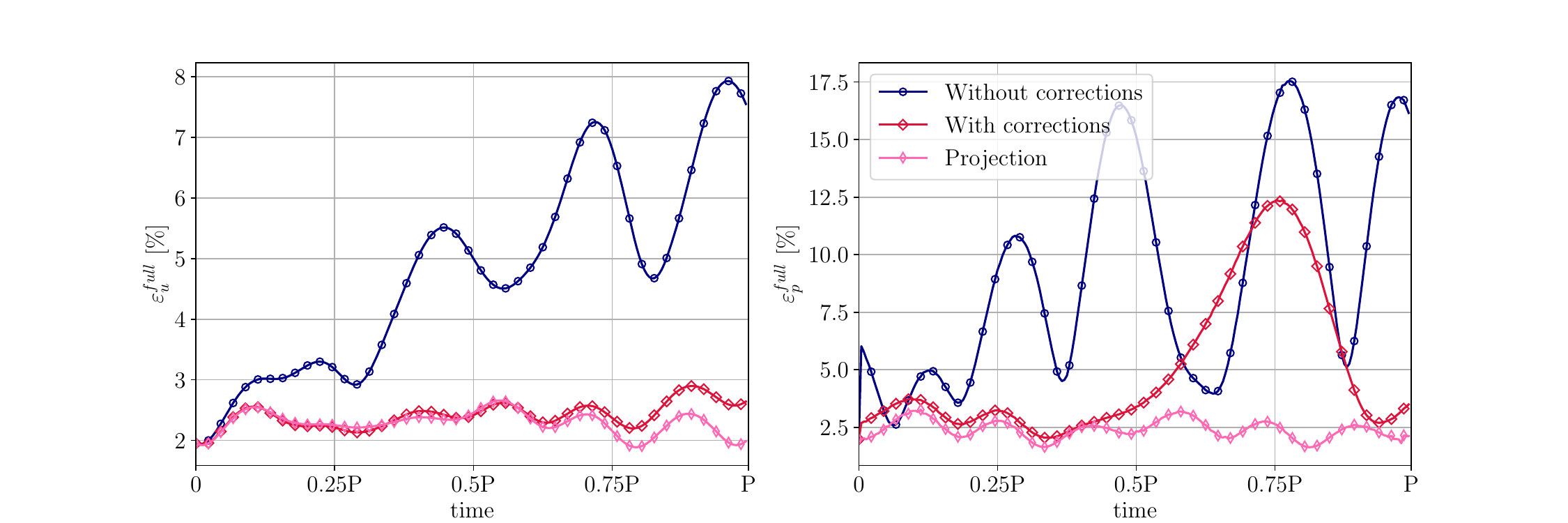} \label{ppe_P} }
\caption{
Relative errors of velocity and pressure for SUP-ROM ($N_u=N_{sup}=N_p=5$ modes) and PPE-ROM ($N_u=N_p=5$), considering $T_{offline}=P$.}
\label{P_window}
\end{figure}

\begin{figure}[h!]
\centering
\subfloat[Relative errors for SUP-ROM]{\includegraphics[ width =0.85 \textwidth, trim={0 0 0 1.5cm}, clip]{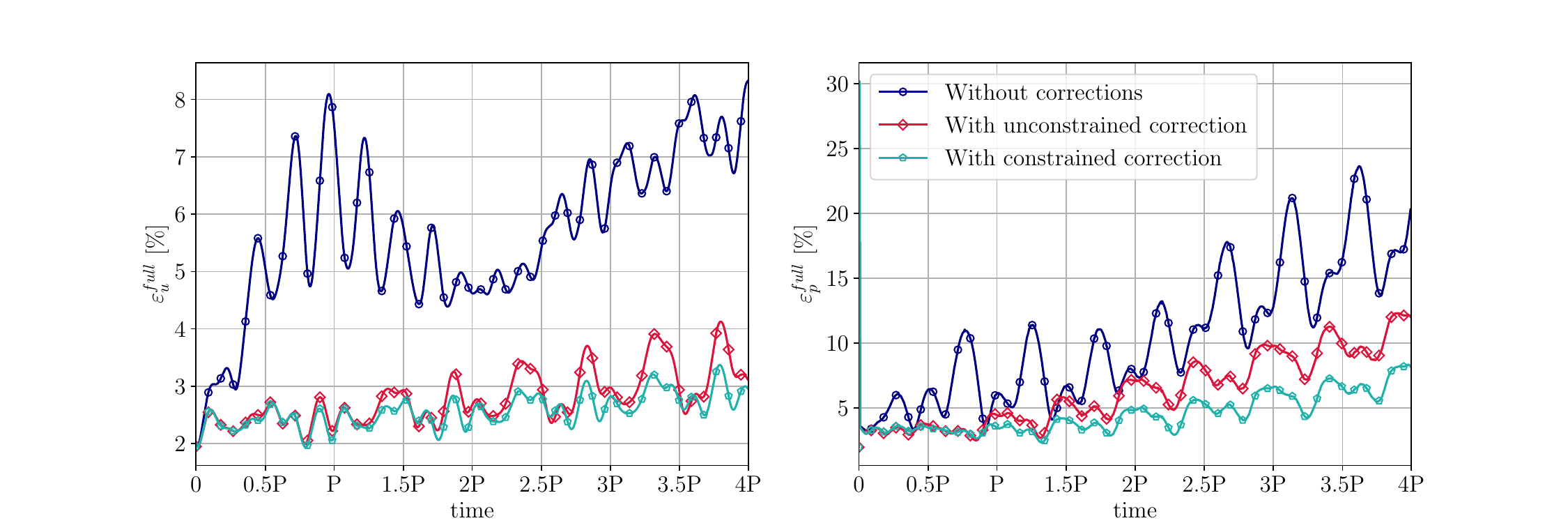} \label{sup_P_pred} }\\
\vspace{-0.3cm}
\subfloat[Relative errors for PPE-ROM]{\includegraphics[ width =0.85 \textwidth, trim={0 0 0 1.5cm}, clip]{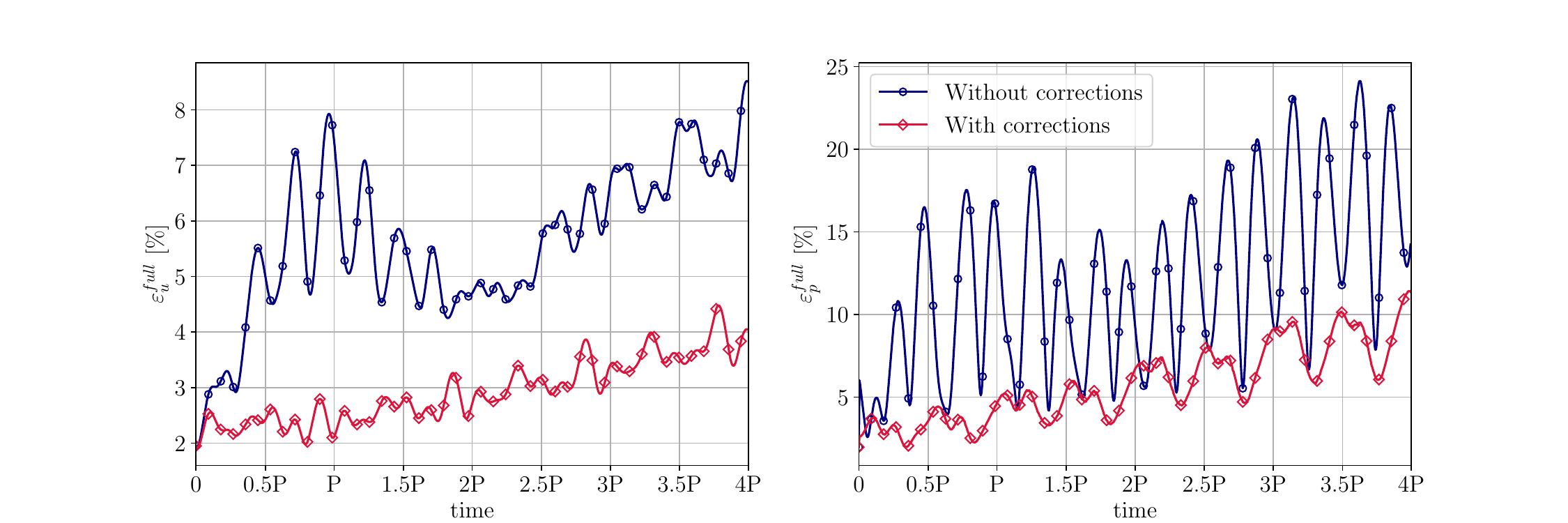} \label{ppe_P_pred} }
\caption{
Relative errors of velocity and pressure for SUP-ROM ($N_u=N_{sup}=N_p=5$ modes) and PPE-ROM ($N_u=N_p=5$), considering a time window larger than the offline time interval, i.e.,  $T_{\text{online}}=4P$ and $T_{\text{offline}}=P$.}
\label{P_window_predictive}
\end{figure}

}

%% file: markers-plots/mark0.tex
\begin{tikzpicture}
\draw[navyblue](0,0) -- (0.8,0);
\draw [navyblue] (0.4,0) circle (2pt);
\end{tikzpicture}

%% file: markers-plots/mark1.tex
\begin{tikzpicture}
\draw[mediumorchid](0,0) -- (0.8,0);
\node[regular polygon,draw, regular polygon sides = 3, rotate=180, color=mediumorchid, scale=0.3] (p) at (0.4,0) {};
\end{tikzpicture}

%% file: markers-plots/mark2.tex
\begin{tikzpicture}
\draw[crimson](0,0) -- (0.8,0);
\node[regular polygon, draw, color=crimson, scale=0.3, regular polygon sides=4, rotate=45] (p) at (0.4,0) {};
\end{tikzpicture}

%% file: markers-plots/mark3.tex
\begin{tikzpicture}
\draw[lightseagreen](0,0) -- (0.8,0);
\node[regular polygon,draw, color=lightseagreen, scale=0.4] (p) at (0.4,0) {};
\end{tikzpicture}

%% file: markers-plots/mark4.tex
\begin{tikzpicture}
\draw[royalblue](0,0) -- (0.8,0);
\node[regular polygon, draw, color=royalblue, regular polygon sides=4, scale=0.4] (p) at (0.4,0) {};
\end{tikzpicture}

%% file: markers-plots/mark6.tex
\begin{tikzpicture}
\draw[goldenrod] (0,0) -- (0.8,0); 
\color{goldenrod}
\pgfplothandlermark{\pgfuseplotmark{star}}
\pgfplotstreamstart
\pgfplotstreampoint{\pgfpoint{0.4cm}{0cm}}
\pgfplotstreamend
\end{tikzpicture}

%% file: markers-plots/mark5.tex
\begin{tikzpicture}
\draw[hotpink](0,0) -- (0.8,0);
\node[diamond, draw, color= hotpink, scale=0.3, aspect=0.5] (p) at (0.4,0) {};
\end{tikzpicture}

%% file: sections_paper1/section_conclusions.tex
\newpage

\section{Conclusions and Outlook}

One popular way to increase the accuracy of Galerkin ROMs in the under-resolved or marginally-resolved regimes is adding closure or correction terms, respectively.
In this paper, we investigated data-driven modeling of these closure and correction terms.
Specifically, we leveraged the data-driven variational multiscale ROM (DD-VMS-ROM) framework to construct for the first time closure and correction terms that significantly increase the pressure accuracy.
To this end, in the offline stage, we postulate an ansatz (i.e., a model form) for the closure and correction terms, and then solve a least squares problem to determine the ansatz parameters that yield the best fit between the exact (i.e., computed with FOM data) closure and correction terms and the ansatz.
We emphasize that developing accurate ROM pressure models is critical, e.g., in computing important engineering quantities (e.g., lift and drag) or when the snapshots used to construct the ROM basis are not weakly divergence-free (e.g., when FOMs enforce the continuity equation weakly).

In our numerical investigation of the novel pressure DD-VMS-ROM, we considered the two-dimensional flow past a circular cylinder at $Re= \num[group-separator={,}]{50000}$ in the marginally-resolved regime.
We also considered several model configurations.

First, we tested two fundamentally different ROM pressure formulations:

(i) The supremizer ROM (SUP-ROM), in which additional (supremizer) modes for the velocity approximation are used in order to satisfy the inf-sup condition.

(ii) The pressure Poisson equation (PPE-ROM), in which the pressure approximation is determined by solving a Poisson equation instead of the continuity equation.

In the SUP-ROM investigation, the introduction of the velocity correction improves the 
approximations of the velocity and pressure fields with respect to the standard SUP-ROM, as can be seen in section \ref{velresults}. However, the improvement 
of the approximation of the reduced pressure field is not as significant as that observed 
for the reduced velocity field. 

In order to further improve the pressure accuracy, pressure correction terms are proposed and added to the reduced system (section \ref{presSUP_2}). These terms do not appear to change the solution in a significant way and their effect is negligible with respect to the one of the velocity correction term.
One possible reason for this fact is that the SUP-ROM formulation does not include a dedicated pressure equation, in which a correction term can directly affect the pressure field.

Therefore, 
a different formulation, 
the PPE-ROM, is used, 
in which the Poisson pressure equation is 
employed. This formulation allows for the introduction of data-driven pressure correction terms, which 
leads to a significant improvement of the reduced pressure 
accuracy, as can be seen in section \ref{PPE1}.

\RA{Moreover, 
the proposed method is certified using different time intervals to collect data for the POD and for the correction terms (subsection \ref{more_tests}). We also 
study the existence of a correlation between the Strouhal number and the selection of the time intervals. 
This analysis leads to the conclusion that one flow period provides enough information to build the data-driven corrections.}
%

The methods developed in this paper can be applied and extended in different research directions. The pressure data-driven corrections developed in this work, when introduced in the SUP-ROM formulation, have not significantly improved the results of the standard formulation in terms of pressure accuracy. However, the supremizer approach, first introduced in \cite{ballarin2015supremizer} and explored in \cite{stabile2018finite}, has been a successful technique for the stabilization of the POD-Galerkin ROMs. Therefore, further data-driven terms including the reduced pressure coefficients should be explored and tested to identify an effective pressure data-driven correction for the SUP-ROM formulation. 

Moreover, the data-driven techniques developed in this paper are embedded in a fluid dynamics framework that does not include any turbulence modeling. In \cite{hijazi2020data}, instead, turbulence modelling is studied by including other data-driven terms in the reduced formulation. 
The present work can be extended by comparing and combining the two data-driven techniques in the reduced formulation.

%% file: sections_paper1/section_acknowledgements.tex
\section*{Acknowledgements}

We thank Prof. Claudio Canuto for his constant support.
We acknowledge the support by the European Commission H2020 ARIA (Accurate ROMs for Industrial Applications, GA 872442) project, by MIUR (Italian Ministry for Education University
and Research) and PRIN "Numerical Analysis for Full and Reduced
Order Methods for Partial Differential Equations" (NA-FROM-PDEs) project, and by the European
Research Council Consolidator Grant Advanced Reduced Order Methods with Applications in
Computational Fluid Dynamics-GA 681447, H2020-ERC COG 2015 AROMA-CFD. 
We also acknowledge support through National Science Foundation Grant Number DMS-2012253.
The main computations in this work were carried out by the usage of ITHACA-FV \cite{ithacasite}, a library maintained at
SISSA mathLab, an implementation in OpenFOAM \cite{ofsite} for reduced order modeling techniques. Its
developers and contributors are acknowledged.